%% file: main.tex
\documentclass[twocolumn,10pt]{IEEEtran}
\usepackage{ifpdf, flushend,subfigure}
\ifCLASSINFOpdf
  \usepackage[pdftex]{graphicx}
  \graphicspath{{../pdf/}{../jpeg/}}
  \DeclareGraphicsExtensions{.pdf,.jpeg,.png}
\else
  \usepackage[dvips]{graphicx}
  \graphicspath{{../eps/}}
  \DeclareGraphicsExtensions{.eps}
\fi

\usepackage{hyperref,multirow,tikz,enumitem}
\usepackage{xspace,amsmath,amsfonts,amssymb,tikz,multirow,float,newclude,hyperref,pgfplots, stmaryrd,lipsum,mathtools, tikz, amsthm,changepage }
\usepackage{xspace,amsmath,amsfonts,amssymb,hyperref,tikz,multirow,filecontents, pgfplots,pbox}
\usepackage{caption}
\usepackage{xspace}
\usepackage{latexsym}
\usepackage{xcolor}
\usepackage{amsmath}
\usepackage{hyperref,graphicx}
\usepackage{multirow,colortbl,wrapfig}
\usepackage{enumitem}
\usepackage{comment,verbatim}
\usepackage{array}
\usepackage{tabularx}
\usepackage{booktabs}
\usepackage{cite}
\usepackage{fancyhdr}
\usepackage[linesnumbered,ruled,vlined]{algorithm2e}
\usepackage{comment}
\usepackage{epstopdf}
\usepackage{float}

\newcommand{\argmin}{\operatornamewithlimits{argmin}}

\newcommand{\beq}{\begin{equation}}
\newcommand{\eeq}{\end{equation}}
\newcommand{\beqn}{\begin{eqnarray}}
\newcommand{\eeqn}{\end{eqnarray}}
\newcommand{\beqno}{\begin{eqnarray*}}
\newcommand{\eeqno}{\end{eqnarray*}}
\newcommand{\bma}{\begin{displaymath}}
\newcommand{\ema}{\end{displaymath}}
\newcommand{\bnu}{\begin{enumerate}}
\newcommand{\enu}{\end{enumerate}}
\newcommand{\bce}{\begin{center}}
\newcommand{\ece}{\end{center}}
\newcommand{\btb}{\begin{tabular}}
\newcommand{\etb}{\end{tabular}}

\newcommand{\bmip}{\textbf{BMIP}}
\newcommand{\bmipd}{\textbf{BMIP1}}
\newcommand{\ebmipd}{\textbf{BMIPe}}
\newcommand{\bmipl}{\textbf{SMIP}}

\def\bT{{\mathbf{T}}}
\def\bt{{\mathbf{t}}}
\def\bx{{\mathbf{x}}}
\def\by{{\mathbf{y}}}
\def\bz{{\mathbf{z}}}
\def\bp{{\mathbf{p}}}
\def\bq{{\mathbf{q}}}
\def\bu{{\mathbf{u}}}
\def\br{{\mathbf{r}}}
\def\zero{{\mathbf{0}}}

\newtheorem{theorem}{Theorem}[section]

\newtheorem{proposition}[theorem]{Proposition}

\newtheorem{remark}[theorem]{Remark}
\newtheorem{assumption}[theorem]{Assumption}

\begin{document}

\title{A Mixed-Integer Bi-level  Model for Joint Optimal Edge  Resource Pricing and Provisioning}

\author{\IEEEauthorblockN{Duong Thuy Anh Nguyen,~\IEEEmembership{Student Member,~IEEE}, Tarannum Nisha,~\IEEEmembership{Student Member,~IEEE}, \\  Ni Trieu,~\IEEEmembership{Member,~IEEE}, and Duong Tung Nguyen,~\IEEEmembership{Member,~IEEE}}  %\vspace{-0.2cm}
\thanks{D. T. A. Nguyen, Ni Trieu, and D. T. Nguyen are with the School of Electrical, Computer and Energy Engineering, Arizona State University, Tempe, AZ,  United States. Email: \{dtnguy52,nitrieu,duongnt\}@asu.edu.
T. Nisha is with  Theory+Practice, Vancouver, Canada.
Email: tarannum@ece.ubc.ca.
}
 %\thanks{Manuscript received April 19, 2005; revised August 26, 2015.}
 }

% make the title area

\maketitle

\begin{abstract}
This paper studies the joint optimization of edge node activation and resource pricing in edge computing, where an edge computing platform provides heterogeneous resources to accommodate multiple services with diverse preferences.
%This problem is formulated 
We cast this problem as a bi-level program, with the platform acting as the leader and the services as the followers.  
The platform aims to maximize net profit by optimizing edge resource prices and edge node activation, with the services' optimization problems acting as constraints.
%Specifically,
Based on the platform's decisions, each service aims to minimize its costs and enhance user experience through optimal service placement and resource procurement decisions. 
The presence of integer variables in both the upper and lower-level problems renders this problem particularly challenging. 
Traditional techniques for transforming bi-level problems into single-level formulations are inappropriate owing to the non-convex nature of the follower problems.
Drawing inspiration from the column-and-constraint generation method in robust optimization, we develop an efficient decomposition-based iterative algorithm to compute an exact optimal solution to the formulated bi-level problem. 
Extensive numerical results are presented to demonstrate the efficacy of the proposed model and technique. 
\end{abstract}

\begin{IEEEkeywords}
 Edge computing, service placement, decomposition, bi-level mixed-integer optimization, dynamic pricing.
\end{IEEEkeywords}
%\IEEEpeerreviewmaketitle
%\vspace{3pt}

\allowdisplaybreaks
\section{Introduction}
The proliferation of mobile devices and applications has led to an explosive increase in mobile data traffic. 
This growth is further amplified by the rise of new services such as augmented/virtual reality (AR/VR),  manufacturing automation, and autonomous driving. 
Therefore, it is imperative to develop innovative solutions that can effectively satisfy the stringent requirements of these evolving applications.
%In response to  these challenges,
To this end, edge computing (EC) has emerged as a pivotal technology, working in tandem with traditional cloud computing to deliver an exceptional user experience and enable %innovative %the seamless implementation of 
low-latency and ultra-reliable Internet of Things (IoT) applications. 
With EC, data processing and storage can be brought closer to the sources, reducing latency and enhancing overall performance \cite{wshi16}.

%Although EC  has experienced rapid growth and holds immense potential for the future,
Despite the rapid growth and immense potential, % of EC, it 
EC is still in its early developmental stages, % of development,
necessitating attention to various challenges. 
%and several challenges need to be addressed. 
One critical challenge is the issue of multi-tenancy concerning shared and heterogeneous edge resources, which is the main focus of this paper.
Specifically, we examine the interaction between an EC platform and multiple services (e.g., AR/VR, Uber, Google Maps). The EC platform, which may be operated by cloud providers, telecommunications companies, or third-party entities like Equinix\footnote{https://www.equinix.com}, manages a set of edge nodes (ENs). It can monetize the available edge resources by offering them for sale to the services. On the other hand, %by deploying %and running 
%the services at the ENs, service providers (SPs) can improve
%provide a more seamless 
%user experience. % by leveraging the proximity and lower latency provided by the EC infrastructure. 
deploying services at the ENs enables service providers (SPs) to provide a superior user experience.
While the placement and execution of services at the ENs enable SPs to improve service quality, %provide superior user experience, %offer their users enhanced performance and better responsiveness, 
it also %results in increased 
entails increased expenses due to the costs associated with storing and running the services. 
 
Our work aims to address two fundamental questions. The first question centers on \textit{determining which ENs should be active for cost reduction and how the platform can optimally set edge resource prices}, 
considering the diverse preferences and valuations of the services towards different ENs. The second question pertains to the \textit{strategic placement of services and the optimal allocation of resources} that SPs should procure from each EN. The complexity arises from the inherent interdependence between the decisions made by the platform and the services.  The service placement and resource procurement decisions of the services rely on the EN activation decision and resource prices established by the platform. Conversely, the EN activation and pricing decisions of the platform are influenced by the resource demands expressed by the services.

Furthermore, the heterogeneity of the ENs adds another layer of complexity, as services may have diverse preferences towards them. Typically, services prefer low-priced edge resources that are geographically close to their operations. However, this can lead to some ENs being over-demanded, particularly those located in or near high-demand areas, while other nodes may experience underutilization, resulting in suboptimal resource utilization. 
To address this issue, the platform can strategically adjust the resource prices of underutilized ENs to encourage load shifting from overloaded ENs.
This pricing strategy is crucial in % plays a crucial role in 
achieving workload balance and maximizing the platform's revenue, especially when multiple services compete for limited computing resources.

In pursuit of our research objectives,  we propose a novel approach that casts the joint edge resource pricing and allocation problem as a bi-level mixed-integer program (BMIP) \cite{asin18,jmore90,bzeng14,Saharidis2009,Fischetti2017,Huang2020,Kaiyi2023}, where the platform acts as the leader, and each service operates as a follower.  The proposed model aims to maximize the platform's profit while minimizing costs for individual services. In our model, the platform first optimizes EN activation, determining which ENs should be active or not %activated or deactivated 
to minimize its operational costs. 
It further computes the optimal resource prices to assign to different ENs, aiming to maximize its revenue, while anticipating the followers' reactions. % of the followers. 
Each service then utilizes this information to make  informed decisions regarding service placement, resource procurement, and workload allocation across different ENs to minimize costs while %simultaneously
enhancing the user experience.
Specifically, each service solves a follower problem
to identify the most suitable ENs for service installation (service placement) 
and determine the optimal amount of resources procured from each EN (resource procurement), considering its operational constraints.

This model builds upon our previous work \cite{tnduong21}, but there are significant differences in terms of both modeling and solution approaches.
In \cite{tnduong21}, the platform is assumed to be responsible for service placement, and the services do not incur costs for service installation %at the ENs
or storage expenses.
Consequently, the follower problems in \cite{tnduong21} were formulated as linear programs (LPs), facilitating the application of Karush–Kuhn–Tucker (KKT) or LP duality reformulation to convert them into sets of equivalent linear equations. This enabled the transformation of the bi-level problem into a single-level mixed-integer linear program (MILP), which could be efficiently solved using standard solvers such as  Mosek\footnote{https://www.mosek.com/} and Cplex\footnote{https://www.ibm.com/analytics/cplex-optimizer}.

However, in practice, services %are 
typically %required to 
bear the costs related to service installation and storage (i.e., service placement cost). This is crucial to ensure fairness among services, as without considering these costs,   services with high placement expenses (e.g., services demanding extensive storage, high bandwidth costs for data transfer, or those with expensive licensing fees per installation)  could utilize their budgets to acquire more resources, thereby gaining unfair advantages over services with lower placement costs.
To overcome this issue, we examine the case where services must cover the service placement costs, which introduces significant complexity to the follower problems.   Specifically, the follower problems become 
MILPs, rendering them non-convex. 
This non-convexity poses a challenge as traditional methods cannot be directly applied to solve the resulting bi-level problem that involves integer variables in the follower problems.
Therefore, a new approach is necessary to address this problem.
Furthermore, with the ability to handle integer variables in the follower problems,  our model goes beyond computing resource pricing and incorporates the design of storage pricing, which was not considered in our previous work  \cite{tnduong21}.

Our proposed model addresses the challenging bi-level mixed-integer program (BMIP) with mixed-integer structures in both upper and lower levels, in the context of EC, which to the best of our knowledge, has not been tackled before. 
%The system comprises two DMs, and due to the mixed-integer structures in both levels, conventional reformulation techniques that depend on lower-level convexity cannot be applied.
Finding a solution to the problem is extremely challenging.
Indeed, bi-level optimization problems, including even the simplest bi-level linear program, are known to be hard to solve   \cite{nphard90}. However, based on the following key observation, we can develop an iterative algorithm for computing an exact global optimal solution to the formulated BMIP. Specifically, it can be observed that if the service placement decisions are fixed (i.e., no lower-level integer variables), the follower problems become convex and we can employ traditional approaches to solve the underlying bi-level problem. Therefore, one approach to address the
presence of integer variables in the follower problems is to enumerate all possible values of the integer variables in those problems. % lower-level problems.
Unfortunately, this method leads to a large-scale MILP with an exponential number of constraints, rendering the problem intractable. 

To overcome this challenge, inspired by the Column-and-Constraint Generation (CCG) approach proposed for solving two-stage robust optimization (RO) \cite{zeng13}, we devise an iterative decomposition algorithm that generates solutions sequentially while ensuring that all possible integer variables are examined within a finite number of iterations. The core idea behind CCG is to use a partial enumeration over a collection of all extreme points of the uncertainty set in RO. Here, our algorithm exploits a similar idea by expanding a partial enumeration over the collection of all possible integer variables in the lower problems by identifying and gradually adding the most critical vector of lower-level integer variables.
The proposed algorithm operates in a master-subproblem framework where the master problem (MP) is a relaxed version of the original problem, thus producing an upper bound (UB), while the subproblem produces a lower bound (LB) for the original problem. 
This process of alternating between the upper and lower-level problems continues until convergence. 
Our contributions can be summarized as follows:

\begin{itemize}
    \item \textit{Modeling:}We formulate a novel bi-level mixed-integer optimization model for joint edge resource pricing and provisioning, employing a Stackelberg game framework to capture the strategic interactions between the EC platform and SPs. %The EC platform, conceptualized as a third-party entity aggregating resources from multiple infrastructure providers (e.g., Google, Microsoft, Amazon, Equinix), enables SPs to select from various options. 
    Unlike previous Stackelberg models \cite{mliu18,valerio20,yswang17,wzhang20} that often rely on simplified constraints for closed-form solutions, our model integrates both computing and storage resource pricing while incorporating the practical costs and constraints related to service placement.
    
    \item \textit{Proposed Methodology and Analysis:} The inclusion of integer variables in the follower problems explicitly accounts for service placement costs; however, this renders them non-convex, making the overall bi-level problem particularly challenging to solve. We propose a decomposition-based iterative algorithm by integrating the CCG approach in RO with a reformulation procedure to achieve an exact optimal solution. 
    This reformulation process involves introducing duplicated lower-level variables, exploring potential values for lower-level integer variables to derive valid inequalities, and applying KKT conditions or LP duality, along with linearization techniques. The proposed algorithm operates in a master-subproblem framework. A detailed analysis of the feasibility and optimality of the proposed model, equivalence between the original formulation and its reformulation, as well as computational complexity, are presented.
    
    \item \textit{Numerical Simulation:} Extensive numerical results are shown to demonstrate the efficacy of the proposed scheme. It also shows that the proposed algorithm can converge efficiently after a small number of iterations.
\end{itemize}

The rest of the paper is organized as follows.
Section \ref{model} presents the system model followed by the problem formulation in  Section \ref{probfor}. 
The solution approaches and simulation
results are presented in Section   \ref{sol} and Section \ref{sim}, respectively.  Section \ref{rel} provides an overview of the related work. Finally, conclusions are presented in   Section \ref{conc}.

\section{System Model}
\label{model}
A typical edge network architecture comprises three layers: the cloud layer, an edge layer with multiple ENs, and an aggregation layer with numerous access points (APs) such as base stations, switches, and routers.  These APs collect service requests from users located in different areas.  Depending on whether
the service is hosted on the ENs or in the cloud, the requests
are routed for further processing accordingly. This approach brings the advantage of localized data processing and analytics, facilitating better scalability, %localized infrastructure, 
and edge resource sharing \cite{duong20,ella23}. In this paper, we consider a system consisting of an EC platform,  also known as an edge infrastructure provider, and a set $\mathcal{K}$ of $K$ services. 
The platform manages a set $\mathcal{J}$ of $J$ geographically distributed heterogeneous ENs to provide compute resources to the services. Each service serves users located in different areas, each of which is represented by an AP. Let  $\mathcal{I}$ denote the set of $I$ areas. The area index, EN index, and service index are %represented by 
$i$, $j$, and $k$, respectively.

\begin{figure}[hbt!]
\vspace{-0.2cm}
	\centering
       \includegraphics[width=0.45\textwidth]{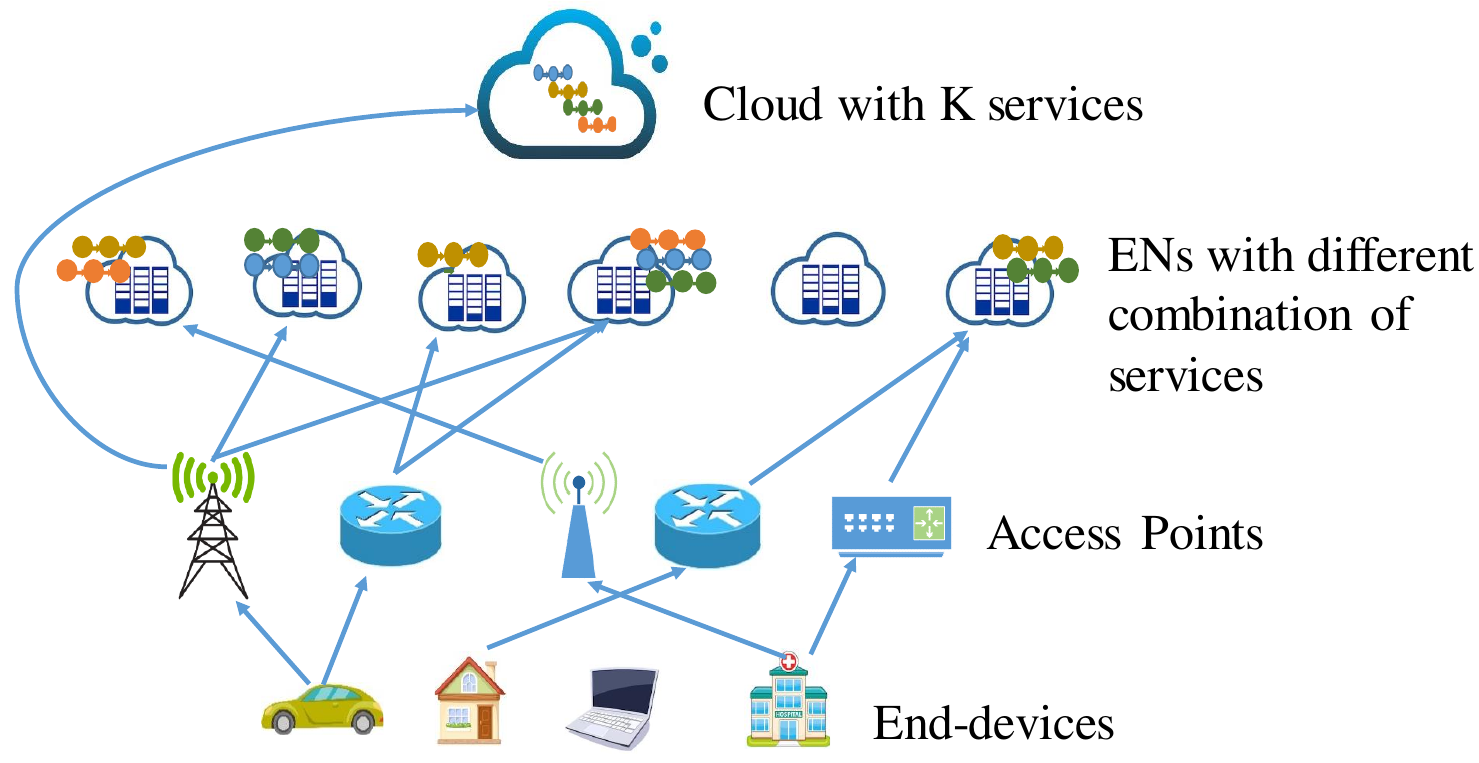}
         	\caption{System model}
	\label{fig:model}
	\vspace{-0.2cm}
\end{figure}

Requests from users are routed based on whether the service is placed on the ENs or the cloud, and further processing is carried out accordingly. The SPs have the option to enhance the Quality of Service (QoS) by delegating tasks to the ENs and directly serving user requests at the ENs with on-site data processing, which offers lower latency but incurs higher costs. Overall, the main objective of the EC platform is to maximize its profit from resource selling, while the services aim to make optimal decisions regarding service placements, resource procurement, and workload allocation. These decisions are made to minimize operating costs and improve the user experience. Fig.~\ref{fig:model} depicts the system model.

The platform maximizes its revenue by deciding the appropriate pricing for edge resources. Since the services have varying preferences for the ENs, some may become over-demanded while others remain under-demanded. Hence, efficient pricing of edge resources is essential to balance supply and demand. To this end, price vectors $\bp = (p_1, p_2, \ldots, p_J)$ and $\bp^s = (p^s_1, p^s_2, \ldots, p^s_J)$ are defined, where $p_j$ represents the unit price of computing resources and $p^s_j$ is the unit price of storage at EN $j$. The platform also needs to decide which ENs to activate to minimize total operational costs. This decision is represented by a binary vector $\bz = (z_1, z_2, \ldots, z_J)$, where $z_j$ equals $1$ if EN $j$ is active and $0$ otherwise.

Following the decisions made by the platform, each service aims to minimize the resource procurement and service placement cost while improving QoS through workload allocation decisions. To improve user experience, a delay-sensitive service prefers to have its requests processed by ENs closer to its users rather than the remote cloud. SPs can provision their services on different distributed ENs, but each service $k$ may have different size $s^k$ and specific hardware and/or delay requirements. Consequently, only a subset of ENs that meet these requirements can handle a service's requests. Each service's goal is to optimally place its services on different ENs, given the resource pricing information and activation decisions of each EN. We denote $\bt^k \!=\! (t_1^k, t_2^k, \ldots, t_J^k)$ and $\bt \!=\! (\bt^1, \bt^2, \ldots, \bt^K)$ as the installation decision of the services, where the binary variable $t_j^k$ indicates whether service $k$ chooses to install on EN $j$ or not. Each EN $j$ has a computing resource capacity of $C_j$, and each service $k$ has a budget $B^k$ for resource procurement. We use the variables$y_0^k$ and $y_j^k$ to represent the amount of computing resource that service $k$ purchases from the cloud and EN $j$, respectively, and denote $\by^k \!=\! (y_0^k,y_1^k, y_2^k, \ldots, y_J^k)$.

In addition to the procurement and placement decisions, each service aims to optimally distribute its workload for further processing. The resource demand (workload) of service $k$ at AP $i$ is denoted by $R_i^k$. The requests from users can be handled either by the cloud server or the ENs, depending on where the service is installed. Let $x_{i,0}^k$ represent the amount of service $k$'s workload at AP $i$ that is directed to the cloud, and $x_{i,j}^k$ represent the workload assigned to EN $j$. 
% \rev{
% Let $\theta^k$ be a binary indicator 
% that indicates if service k is allowed to drop user
% requests. % from its users. 
% If $\theta^k=0$, unmet demand is not allowed.
% %The services may consider dropping requests from users when $\theta^k=1$, or they may enforce that unmet demand is not allowed, i.e., when $\theta^k=0$.
% When $\theta^k=1$, service requests from each area $i$ that are not served by some ENs or the cloud will be counted as unmet demand $q_i^k$.}
Let $q_i^k$ denote the amount of unmet demand of service $k$ from area $i$.
We define $d_{i,j}$ as the network delay between AP $i$ and EN $j$, and $d_{i,0}$ as the delay between AP $i$ and the cloud. 
The goal of each service is to minimize not only the service placement and resource procurement cost but also the network delay experienced by its users. In order to capture these information, we denote $\bx_0^k = (x_{1,0}^k, x_{2,0}^k, \ldots, x_{I,0}^k)$, $\bx_i^k = (x_{i,1}^k, x_{i,2}^k, \ldots, x_{i,J}^k)$, $\bx^k = (\bx_0^k,\bx_1^k, \bx_2^k, \ldots, \bx_I^k)$, and $\bq^k = (q_1^k,q_2^k, \ldots, q_I^k)$.

\begin{table}[t!] 
\centering
\caption{NOTATIONS}
\begin{tabular}{|l|l|}
\hline
Notation   & Meaning\\	
\hline	
Indices & \\ \hline
$i$, $\mathcal{I}$, $I$ & Indices, set, and number of APs\\
$j$, $\mathcal{J}$, $J$  & Indices, set, and number of ENs\\
$k$, $\mathcal{K}$, $K$ &  Indices, set, and number of services\\
\hline
Parameters & \\ \hline
$C_j$ & Computing capacity of EN $j$\\
$S_j$ & Storage capacity of EN $j$\\
$d_{i,j}$ & Network delay between area $i$ and EN $j$\\
$d_{i,0}$ & Network delay between area $i$ and the cloud\\
$D^{k, \sf m}$ &   Delay threshold of service $k$ \\
$d_{i}^{k, \sf a}$ &   Average delay  of service $k$ in area $i$ \\
$B^{k}$ & Budget of service $k$ \\
$s^k$ & Size of  service $k$\\
$w^k$ & Delay penalty for service $k$ \\
% $\theta^k$ & Binary parameter, $1$ if service $k$ can drop requests\\
$\psi_i^k$ & Penalty for unmet demand in area $i$ for service $k$\\
$R_{i}^{k}$ & Resource demand for service $k$ in area $i$ \\
$\phi_j^k$ & Installation cost of service $k$ at EN $j$\\
$c_j$ & Variable cost at EN $j$\\
$f_j$ & Fixed cost at EN $j$\\
$p_0$ & Unit resource price at the cloud\\
\hline
Variables & \\ \hline
$p_j$ & Unit price of computing resource at EN $j$\\
$p_j^s$ & Unit price of storage resource at EN $j$\\
$z_j$ & Binary variable, $1$ if EN is active \\
$t_j^k$ & Binary variable, $1$ if service $k$ is placed at EN $j$ \\
$q_i^k$ & Unmet demand in area $i$ for service $k$\\
$x_{i,j}^k$ & Workload of service  in area $i$ assigned to EN $j$\\
$x_{i,0}^k$ & Workload of service in area $i$ assigned to the cloud \\
$y_j^k$ & Amount of resource purchased at EN $j$ for service $k$\\
$y_0^k$ & Amount of resource purchased at the cloud for service $k$\\
\hline
\end{tabular}\label{notation}
\vspace{-0.5cm}
\end{table}

Our proposed system model takes into account practical issues that arise in EC networks. One such issue is the relationship between network delay and service placement. When a user requests a service, it is advantageous to have the service available on an EN within the user's area to minimize latency. However, this may result in higher service placement costs for the SP since the service needs to be placed on distributed ENs to accommodate user requests at any time. This creates a trade-off between service placement and network delay. Additionally, there is a correlation between resource prices, delay, and service placement. ENs have different resource price configurations, assigning requests from users to cheaper ENs incurs lower resource cost and reasonable computation capacity but at the expense of greater delay. Overloading is another issue that may arise in this scenario. Therefore, selecting suitable ENs for service placement is crucial. This can be addressed by ensuring that each EN has enough computing capacity to accommodate all services and by setting an appropriate delay penalty for each service to reduce delay.

To address these issues, we propose an optimal pricing scheme for the EC platform to assist services in minimizing their costs while maximizing overall profit. The first stage of the scheme determines the optimal resource pricing and the number of active servers and EN activation to maximize profit. The second stage uses these optimal solutions to determine optimal service placement, resource procurement, and workload allocation. Table \ref{notation} summarizes the main parameters and variables used throughout the paper.  %in the profit maximization problem of the EC platform and the cost minimization problem of the services.

\begin{remark}
While our model focuses on a single EC platform, it can be conceptualized as a third-party entity aggregating resources from various infrastructure providers (e.g., Google, Microsoft, Amazon, Equinix, third-party edge data centers). %, allowing service providers (SPs) to select among multiple options. 
Our model can also be applied in scenarios where multiple SPs, especially those providing delay-sensitive services, want to procure edge resources from a specific telecom company (e.g., AT\&T, Verizon) to serve their subscribers using that carrier. Extending the proposed model to a multi-leader-multi-follower game to capture the competition among multiple EC platforms is the subject of our future work.
\end{remark}

\section{Problem Formulation}
\label{probfor}
\begin{figure}[t!]
    \centering
      \includegraphics[width=0.48\textwidth]{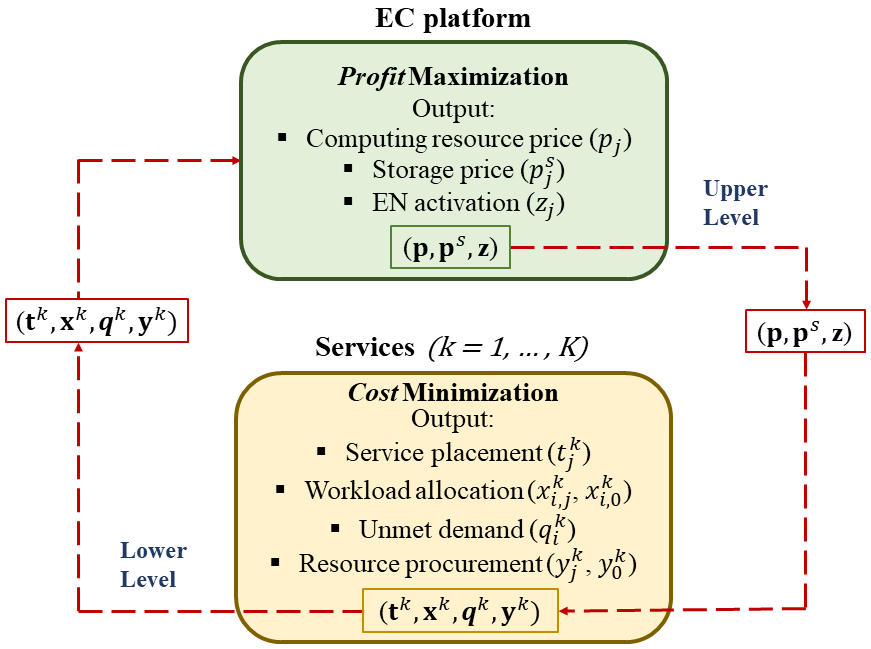}
     \caption{Interactions between the EC platform and services.}
    \label{fig: Interaction}
    \vspace{-0.5cm}
\end{figure}

This section presents the mathematical formulation of the interaction between the EC platform and the SPs as a bi-level program. The upper-level optimization problem represents the leader (the EC platform), while the $K$ lower-level problems represent the followers (the SPs), with one problem for each SP. The strategic interaction between the EC platform and the SPs is governed by Stackelberg game dynamics, as illustrated in Figure~\ref{fig: Interaction}. The platform aims to maximize its profit by solving the leader problem and subsequently informs the services of the resource price and EN activation decisions.
Each service, upon receiving this information, solves a follower problem under budget and delay constraints to minimize its costs, and then communicates the optimal solutions for service placement, resource procurement, and workload allocation back to the platform. 
 
The optimal solutions of the followers and the leader are closely intertwined.  Specifically, the followers' decisions provide input to the leader problem, while the output of the leader problem directly influences the followers’ decisions. Indeed, the follower problems act as constraints in the leader problem. This section delves into the description of the follower problem for each service, the leader problem for the platform, as well as the bi-level optimization model.

\subsection{Follower Problems}
In the following, we describe the lower-level optimization problem for each service $k$. The edge resource price and EN activation decisions computed by the leader are inputs for the follower problems.
Each service aims to optimally select  ENs for service placement and strategically allocate workload among the selected ENs and the cloud to minimize its monetary cost while enhancing the user experience.
\subsubsection{Objective Function}
The total cost of service $k$ is: 
\beqn
\label{eq:objk}
\mathcal{C}^k =  C^{{\sf c},k} + C^{{\sf e},k} + C^{{\sf p},k} + C^{{\sf d},k} + C^{{\sf u},k}
\eeqn
where $\mathcal{C}^k$ is the total cost of service $k$, and $C^{{\sf c},k}$, $C^{{\sf e},k}$, $C^{{\sf p},k}$, $C^{{\sf d},k}$, and $C^{{\sf u},k}$ denote the cloud resource procurement cost, edge resource procurement cost, service placement cost, delay penalty, and unmet demand penalty, respectively.

The total cost of service $k$ for purchasing computing resources in the cloud can be calculated as:
\beqn
\label{eq:ccost}
C^{{\sf c},k} =  p_0 y_0^k,% \quad \forall k,
\eeqn
where $p_0$ is the unit price of computing resources in the cloud. Since the capacity of a cloud data center is usually very large, we assume that every service has been installed and is available in the cloud at the beginning.

The total cost of service $k$ for acquiring computing resources at the ENs can be expressed as follows:
\beqn
C^{{\sf e},k} = \sum_j p_j y_j^k,% \quad \forall k,
\eeqn
where $p_j$ is the unit price of computing resources at EN $j$.

The service placement cost of service $k$ at EN $j$ includes a fixed cost $\phi_j^k$ % encapsulates cost of
for downloading and installing the service data and software as well as the storage cost  $s^kp_j^s$ depending on the unit price $p^s_j$ of storage at EN $j$. Thus, the cost of placing service $k$ on EN $j$ is $\phi_j^k+s^kp_j^s$. Consequently, the total placement cost of service $k$ is:
\beqn
C^{{\sf p},k} = \sum_j \left(\phi_j^k+s^kp_j^s\right) t_j^k.% \quad \forall k,
\eeqn

The delay cost between area $i$ and EN $j$ is proportional to the amount of workload allocated from area $i$ to EN $j$ as well as the network delay $d_{i,j}$ between them. Thus, the total delay cost of service $k$ can be expressed as follows: %using the following equation:
\beqn
C^{{\sf d},k} =  w^k \Big( \sum_i x_{i,0}^k d_{i,0} + \sum_{i,j} x_{i,j}^k d_{i,j} \Big). % \quad \forall k.
\eeqn

Service $k$ incurs a penalty of $\psi_i^k$ for each unit of unmet demand in area $i$, and the corresponding unmet demand penalty is expressed as:
% \beqn
% C^{{\sf u},k} =  \theta^k \sum_i \psi_i^k q_i^k . % \quad \forall k.
% \eeqn }
\beqn
C^{{\sf u},k} =   \sum_i \psi_i^k q_i^k . % \quad \forall k.
\eeqn 

Note that the delay and unmet demand penalties are virtual expenses that capture the service quality and reflect the service's sensitivity to delay and unmet demand. The actual payment incurred by the service encompasses only the initial three terms on the right-hand side (RHS) of equation \eqref{eq:objk}. The penalty parameters $w^k$, %$\theta^k$, 
and $\psi_i^k$, control the tradeoff between monetary costs and the penalty associated with service quality.
Each service can balance its payment, covering resource procurement, placement costs, and the overall penalty for delay and unmet demand, by adjusting these parameters. A higher value of $w^k$ ($\psi_i^k$) indicates increased sensitivity to delay (unmet demand) and a willingness to allocate more resources to mitigate delay (unmet demand). 

Overall, given $(\bp,\bp^s,\bz)$, the objective function $\mathcal{C}^k$ of each SP $k$ is given explicitly as follows:
% \begin{align}\label{eq-SPkObj}
% &\mathcal{C}^k(\bx^k,\bq^k,\by^k,\bt^k|\bp,\bp^s,\bz) \nonumber \\
%  =~ & p_0 y_0^k + \sum_j p_j y_j^k + \theta^k \sum_i \psi_i^k q_i^k +\sum_j \left(\phi_j^k+s^kp_j^s\right) t_j^k \nonumber \\ 
% & + w^k \Bigg( \sum_i x_{i,0}^k d_{i,0} + \sum_{i,j} x_{i,j}^k d_{i,j} \Bigg).
% \end{align}}
\begin{align}\label{eq-SPkObj}
&\mathcal{C}^k(\bx^k,\bq^k,\by^k,\bt^k|\bp,\bp^s,\bz) \nonumber \\
 =~ & p_0 y_0^k + \sum_j p_j y_j^k + \sum_i \psi_i^k q_i^k +\sum_j \left(\phi_j^k+s^kp_j^s\right) t_j^k \nonumber \\ 
& + w^k \Bigg( \sum_i x_{i,0}^k d_{i,0} + \sum_{i,j} x_{i,j}^k d_{i,j} \Bigg).
\end{align}

We are now ready to describe the constraints imposed on the optimization problems of the services. 

\subsubsection{Budget Constraints}
The total payment of  service $k$, including the resource procurement costs and the service placement cost, should not exceed its budget, i.e., we have:
\beqn
\label{spc1}
C^{{\sf c},k} + C^{{\sf e},k} + C^{{\sf p},k} \leq  B^k.% \quad \forall k.
\eeqn

\subsubsection{Service Placement Constraints}
The services should be placed on active ENs only, thus: 
\beqn
\label{spc2}
t_j^k \leq z_j, \quad \forall j.%k,
\eeqn
If EN $j$ is inactive, then $z_j = 0$, which enforces $t_j^k = 0, ~\forall k$.

\subsubsection{Capacity Constraints}
The amount of computing resource that can be purchased at each EN $j$ cannot exceed %is limited to the computing resource 
its capacity $C_j$. % available on it. 
Additionally, a service only purchases resources from the ENs that have the service installed. We have:   
\beqn
\label{spc3}
 y_j^k \leq C_j t_j^k,\quad \forall j.   
\eeqn

\subsubsection{Workload Allocation Constraints} 
The total workload of service $k$ allocated from all areas to the cloud or an EN should not exceed the amount of procured computing resources at the cloud or that node. Hence:
\beqn
\label{spc5}
y_0^k \geq \sum_i x_{i,0}^k,\quad ~~ \\
\label{spc6}
y_j^k \geq \sum_i  x_{i,j}^k, ~ \forall j. 
\eeqn
%In addition, the
The computing resource demand from each area must either be served by the cloud or the ENs or dropped, %(when $\theta^k=1$)}, 
which implies:
% \rev{\beqn \label{spc7}
% x_{i,0}^k + \sum_j x_{i,j}^k + \theta^k q_i^k = R_i^k,~ \forall i. 
% \eeqn }
\beqn \label{spc7}
x_{i,0}^k + \sum_j x_{i,j}^k +  q_i^k = R_i^k,~ \forall i. 
\eeqn

\subsubsection{Delay Constraints}
To guarantee that the average delay of a delay-sensitive service does not surpass a certain delay threshold, the following constraints can be imposed:
\beqn
\label{spc8}
d_{i}^{k, \sf a} \leq D^{k, \sf m},~\forall i, \\%,k \\
d_{i}^{k, \sf a} = \frac{x_{i,0}^k d_{i,0} +  \sum_j x_{i,j}^k d_{i,j}}{ R_i^k},~~\forall i,%k \label{spc9}
\eeqn
where $D^{k, \sf m}$ is the maximum allowable delay threshold for service $k$, and $d_{i}^{k, \sf a}$ is the average delay of service $k$ in area $i$.

\subsubsection{EN Selection Constraints}
Different services may have distinct hardware/software requirements for their operations.  %on an EN.
For instance, a latency-sensitive service may need to be placed on  ENs that are close to its service areas.
Some services can only be deployed on ENs that support specific software and operating systems (e.g., Linux). 
To ensure that these requirements are met, we incorporate a binary indicator, represented as  $a_{i,j}^k$, which determines whether an EN $j$ satisfies the prerequisites necessary for serving requests of service $k$ in area $i$. Hence, the following constraints are added to the model.
\beqn
\label{spc10}
 x_{i,j}^k \leq a_{i,j}^k R_i^k,~~ \forall i,j. %,k.
\eeqn

Overall, the follower problem for service $k$ is given as:
\begin{align}
\label{objf}
&(\Xi^k)~~~ \underset{\bx^k,\bq^k,\by^k,\bt^k}{\min} \qquad \mathcal{C}^k(\bx^k,\bq^k,\by^k,\bt^k|\bp,\bp^s,\bz)\\
&\text{subject to:} \quad~  (\ref{eq:ccost})-(\ref{spc10}) \nonumber \\
\label{spc11}
& \qquad\qquad\quad~~ x_{i,j}^k,x_{i,0}^k \geq 0,~ \forall i,j;~y_j^k,y_0^k \geq 0, ~\forall j,\\
& \qquad\qquad\quad~~ q_i^k \geq 0,~ \forall i;~ t_j^k \in \{0,~1\}, ~\forall j. \label{spc12}
\end{align} 
Note that there are $K$ follower problems, each corresponding to a specific service. These optimization problems will act as constraints for the leader (or platform) optimization problem, guiding the platform in resource provision and pricing decisions during the planning stage. Additionally, it is noteworthy that the computing resource price $\bp$, storage price $\bp^s$, and EN activation $\bz$ are decision variables in the leader problem, while they function as parameters in the follower problems.

\subsection{Leader Problem}
\subsubsection{Objective Function}
The objective of the platform is to optimize its net profit by maximizing the revenue derived from hosting services at the ENs while simultaneously minimizing the overall operational cost. The net profit is computed as the difference between revenue and cost. The platform's revenue is comprised of revenue generated from the sale of edge resources and the revenue obtained from charging services for service placement. Therefore, the platform's revenue can be expressed as:
\beqn
\label{lrev}
\textit{Revenue} = \sum_j p_j \sum_k y_j^k+\sum_j\sum_k \left(\phi_j^k+s^kp_j^s\right) t_j^k.
\eeqn

The cost of the platform is determined by the total cost associated with operating the ENs, which is dependent on the electricity price and power consumption of the nodes. To simplify the computation of the operating cost, a linear function approximation is commonly employed in the literature \cite{xsun20,tnduong21}. Specifically, the operating cost of an EN is composed of a fixed cost $f_j$ and a variable cost that is dependent on the computing resource utilization when the node is active. Conversely, if the EN is inactive (i.e., $z_j = 0$), its operating cost is considered to be zero. The operating cost of EN $j$ can be expressed as follows:
\beqn
\textit{Cost}_j^e = \left(f_j  + c_j  \frac{\sum_k y_j^k}{C_j}\right) z_j, ~~ \forall j.
\eeqn
However, (\ref{spc2}) and (\ref{spc3}) imply $y_j^k = 0$ if $z_j = 0$. Thus, $\textit{Cost}_j^e$ can be written equivalently as $f_j z_j + c_j  \frac{\sum_k y_j^k}{C_j}$. 
% Consequently, the total cost of the platform is:
% \beqn
% \label{eq-PlatformCost}
%  \textit{Cost} = \sum_j \left(f_j z_j + c_j  \frac{\sum_k y_j^k}{C_j} \right).
% \eeqn
% In summary, the platform aims to maximize the following objective function:

In summary, the platform aims to maximize the following objective function:
\begin{align}
\label{eq-PlatformObj}
&\mathcal{P}\left(\bp,\bp^s,\bz,\by,\bt\right)  =  \sum_j\sum_k \left(\phi_j^k+s^kp_j^s\right) t_j^k \nonumber\\
&\qquad\qquad+ \sum_j p_j \sum_k y_j^k - \sum_j \left(f_j z_j + c_j  \frac{\sum_k y_j^k}{C_j} \right).   
\end{align}

The leader problem's constraints are described below.

\subsubsection{Resource Capacity Constraints}
The allocation of computing resources to services can only be made from an active EN. Additionally, the allocated computing resource must not exceed the computing capacity of the EN. Thus, the following constraints are imposed:
\beqn
\label{ipc1}
\sum_k y_j^k \leq z_jC_j; ~~ z_j \in \{0,~1\}, \quad \forall j,
\eeqn
which implies if $z_j = 0$, then $y_j^k = 0$, $\forall j,k$.

\subsubsection{Storage Constraints}
The storage resources allocated to services from an EN are limited by the storage capacity of the EN. Therefore, the following constraint is applied:
\beqn
\label{ipc4}
\sum_k s^k t_j^k \leq z_jS_j,  \quad \forall j.
\eeqn

\subsubsection{Pricing Constraints}
The platform selects the unit price $p_j$ for computing resources at EN $j$ from a discrete set ${p_j^1, \ldots, p_j^V}$, where $v \in \{1, \ldots, V\}$ indexes different price levels. The unit  price $p^s_j$ for storage at each EN $j$ can also be selected from a discrete set $
\{p_j^{s1} , \ldots, p_j^{sH} \}$, where  $h \in \{1, \ldots, H\}$ indicates different storage price options. 
This assumption is reasonable as the set of price options can represent various price levels, such as low, medium, and high prices \cite{tnduong21,ChengBandit}. 
It is worth noting that, in the case of continuous prices, they can be approximated by discretizing them into numerous small intervals (e.g., through binary expansion). Discretizing the prices is essential for effectively managing the linearization procedure discussed later in the solution approach section.

To model the price selection, binary variables $r_j^v$ and $r_j^{sh}$ are introduced, where $r_j^v$ takes the value of $1$ if the unit price of computing resources at EN $j$ is $p_j^v$, and $r_j^{sh}$ is set to $1$ if the unit storage price at EN $j$ is $p_j^{sh}$. As the price can take only a single value, the following constraints are imposed:
\begin{subequations}\label{ipc56}
\begin{alignat}{3}
\label{ipc5}
\!\!&p_j = \sum_{v=1}^V p_j^v r_j^v,\forall j;\sum_{v=1}^V r_j^v = 1,\forall j; r_j^v\in\{0,1\}, \forall j, v,\\
\label{ipc6}
\!\!&p^s_j = \!\sum_{h=1}^H p_j^{sh} r_j^{sh},\forall j; \!\sum_{h=1}^H r_j^{sh} \!=\! 1,\forall j;  r_j^{sh}\!\!\in\!\!\{0,1\}, \!\forall j, h.\!\!\!
\end{alignat}
\end{subequations}

The bi-level mixed-integer program representing the leader problem is presented in a compact form below:
% \begin{align}
% \textbf{\bmip :}~~ \max_{\substack{\bp,\bp^s,\bz,\\\bx,\bq,\by,\bt}} ~ \sum_j p_j &\sum_k y_j^k + \sum_j\sum_k \left(\phi_j^k+s^kp_j^s\right) t_j^k \nonumber\\
% &- \sum_j \left(f_j z_j + c_j  \frac{\sum_k y_j^k}{C_j} \right)  \label{eq:objbmip}
% \end{align}
\begin{subequations}\label{eq:bmip}
\begin{alignat}{3}
&\textbf{\bmip :}~~~~ \max_{\substack{\bp,\bp^s,\bz,\bx,\bq,\by,\bt}} ~~~ \mathcal{P}\left(\bp,\bp^s,\bz,\by,\bt\right)  \label{eq:objbmip}\\
&\text{subject to:}~~~~ (\bp,\bp^s,\bz,\bx,\bq,\by,\bt) \in \mathcal{H}, \label{snum}\\
\label{enum}
&\left(\bx^k, \bq^k, \by^k, \bt^k\right) \!\in \!\!\!\! \argmin_{\substack{(\bx^k\!, \bq^k\!, \by^k\!, \bt^k)\\ \in \mathcal{S}^k(\bp,\bp^s,\bz)}} \mathcal{C}^k(\bx^k,\bq^k,\by^k,\bt^k|\bp,\bp^s,\bz), ~\forall k,\!
% \label{enum}
% &\left(\bx^k, \bq^k, \by^k, \bt^k\right) \!\in \!\!\!\! \argmin_{\substack{(\bx^k\!, \bq^k\!, \by^k\!, \bt^k)\\ \in \mathcal{S}^k(\bp,\bp^s,\bz)}} \!\!\Bigg\{p_0 y_0^k + \!\!\sum_j p_j y_j^k + \!\theta^k \!\sum_i \psi_i^k q_i^k\!\!\!\!\!\!\!\!\!\! \nonumber \\ 
% &\!\! +\!\!\sum_j \left(\phi_j^k+s^kp_j^s\right) t_j^k + w^k \Bigg( \!\sum_i x_{i,0}^k d_{i,0} \!+ \!\!\sum_{i,j} x_{i,j}^k d_{i,j}\! \Bigg)\!\Bigg\}, \forall k,\!
\end{alignat}
\end{subequations}
where the objective function $\mathcal{P}$ of the leader is defined in \eqref{eq-PlatformObj} while the objective function $\mathcal{C}^k(\bx^k,\bq^k,\by^k,\bt^k|\bp,\bp^s,\bz)$ of the follower $k$ is defined in \eqref{eq-SPkObj}. Here, $\mathcal{H}=\left\{(\bp,\bp^s,\bz,\bx,\by,\bt)|~(\ref{ipc1})-(\ref{ipc56})\right\}$ is the set of feasible decisions for the upper-level problem and $\mathcal{S}^k(\bp,\bp^s,\bz)=\left\{(\bx^{k},\bq^{k}, \by^{k}, \bt^{k})|~(\ref{spc1})-(\ref{spc10}),~ (\ref{spc11}),~(\ref{spc12})\right\}$ is a set of feasible decisions for the lower-level problem $k$. We define $\Omega^{IR}=\left\{(\bp,\bp^s,\bz,\bx,\by,\bt)|(\ref{snum}),(\ref{enum})\right\}$ as the feasible set of \textbf{BIMP}. This set is commonly referred to as the inducible region of the bi-level problem.

The platform aims to maximize its profit by jointly optimizing the EN activation and resource pricing decisions. The follower problems serve as constraints of the leader problem, as shown in (\ref{enum}).
In this problem, $\bz$ represents the integer variables at the upper level,  whereas $\bt^k$ represents the lower-level integer variables. Thus, the problem takes the form of a bi-level mixed-integer program. % (BMILP).

\section{Feasibility and Optimality Analysis}
In this section, we analyze the feasibility of the proposed bi-level model \textbf{BMIP} and the existence of the optimal solution to the problem. First, we consider the high point problem (HPP) \cite{jmore90} associated with \textbf{BMIP}--that is, the single-level MIP established from the bi-level problem \textbf{BMIP} when we require the lower-level variables $(\bx^{k},\bq^{k}, \by^{k}, \bt^{k})$ for all $k \in\mathcal{K}$ to be feasible, instead of optimal, given an upper-level decision $(\bp,\bp^s,\bz)$. Specifically, the HPP is defined as follows:
\begin{align}
&\textbf{HPP :} \qquad\quad \max_{\substack{(\bp,\bp^s,\bz,\bx,\bq,\by,\bt) \in \Omega }} ~~~ \mathcal{P}\left(\bp,\bp^s,\bz,\by,\bt\right)  \label{eq:objHPP}
\end{align} 
where $\Omega = \Big\{(\bp,\bp^s,\bz,\bx,\bq,\by,\bt)| (\bp,\bp^s,\bz,\bx,\bq,\by,\bt) \in \mathcal{H}$, $\left(\bx^k, \bq^k, \by^k, \bt^k\right) \in \mathcal{S}^k(\bp,\bp^s,\bz), ~ \forall k\Big\}$.

% \rev{Note that $\Omega^{IR} \subseteq \Omega$, hence, \textbf{HPP} is a relaxation of \textbf{BIMP}. Consequently, if \textbf{HPP} is infeasible then \textbf{BMIP} is infeasible.}

%We have the following result \cite[Prop. 1]{bzeng14}.
%\begin{proposition}\label{prop-HPPinfeasible}
%If \textbf{HPP} is infeasible then \textbf{BMIP} is infeasible.
%\end{proposition}}

The bilinear terms $p_j y_j^k$ and $p_j^s t_j^k$ in \textbf{HPP} can be linearized (see Appendix~\ref{app-BilinearLinearization}), resulting in a single-level MILP that can be solved using off-the-shelf solvers. %Referring to Proposition~\ref{prop-HPPinfeasible}, we
Note that $\Omega^{IR} \subseteq \Omega$, hence, \textbf{HPP} is a relaxation of \textbf{BIMP}. Consequently, if \textbf{HPP} is infeasible then \textbf{BMIP} is infeasible. Thus, we can initially verify the infeasibility of \textbf{HPP} to identify scenarios where \textbf{BMIP} is infeasible without directly solving the intricate bi-level problem. The following can be observed. %Recalling that $\theta^k \in \{0, 1\}$ for all $k \in \mathcal{K}$ is the binary parameter that controls whether service $k$ can drop requests from its users. The following are observed.}

% In the following, we analyze the feasibility of \textbf{HPP} and identify the constraints that can render the \textbf{HPP} problem infeasible.

% \rev{\begin{proposition}\label{prop-HPPanalysis}
% % Let $\theta^k \in \{0, 1\}$ for all $k \in \mathcal{K}$. 
% The feasibility is characterized as follows:
% \begin{enumerate}
%   \item If $\theta^k = 1$ for all $k \in \mathcal{K}$, then given any decision $(\bp, \bp^s, \bz)$ of the CP, the feasible set of the lower-level problem $k$, $\mathcal{S}^k(\bp,\bp^s,\bz)$, is nonempty. Furthermore, \textbf{HPP} is universally feasible.
%   \item If $\theta^k = 0$ for some $k \in \mathcal{K}$, then \textbf{HPP} may become infeasible. The budget and delay constraints in \eqref{spc1} and \eqref{spc8} are critical factors restricting feasibility.
% \end{enumerate}
% \end{proposition}
\begin{proposition}\label{prop-HPPanalysis}
% Let $\theta^k \in \{0, 1\}$ for all $k \in \mathcal{K}$. 
Given any decision $(\bp, \bp^s, \bz)$ of the leader, the feasible set of the lower-level problem $k$, $\mathcal{S}^k(\bp,\bp^s,\bz)$, is nonempty. Furthermore, \textbf{HPP} is universally feasible.
 %\end{enumerate}
\end{proposition}
\begin{proof}
%\textbf{Case 1:} When $\theta^k = 1$ for all $k \in \mathcal{K}$, 
The service requests from each area $i$ can be served either by some ENs, the cloud, or they can be dropped and will be counted as unmet demand. Let $\zero$ represent a vector with all entries equal to $0$. It can be easily verified that for any decision $(\bp, \bp^s, \bz)$ at the upper level, then $\bx_0^k = \zero$, $\bx_1^k = \zero, \ldots, \bx_I^k = \zero$, $\by_1^k = \zero, \ldots, \by_I^k = \zero$, $\by_0^k = \zero$, and $\bq^k = (R_1^k, R_2^k, \ldots, R_I^k)$, for all $k$, satisfy all constraints of the lower level problems. Thus, the set $\mathcal{S}^k(\bp, \bp^s, \bz)$ is nonempty. Furthermore, as the SPs in the lower level do not procure any resources from any EN, for any given price $(\bp, \bp^s)$ and activation decisions $\bz$, the upper-level constraints are satisfied. Consequently, \textbf{HPP} is universally feasible.
% \textbf{Case 2:} Assume that $\theta^{\tilde{k}} = 0$ for $\tilde{k} \in \tilde{\mathcal{K}}$, where $\tilde{\mathcal{K}} \subset \mathcal{K}$.  In accordance with Case 1, it is established that $\mathcal{S}^k(\bp, \bp^s, \bz)$ is nonempty for $k \in \mathcal{K}\backslash \tilde{\mathcal{K}}$ due to $\theta^k = 1$. For SPs $\tilde{k} \in \tilde{\mathcal{K}}$, considering the disallowance of unmet demand, the service requests must be served either at the cloud or at some ENs. We note that the cloud has virtually unlimited capacity with lower prices than the resources at ENs but is usually situated further from the users. The decision for SPs $\tilde{k}$ to assign requests to the cloud is constrained by the total average delay exceeding the threshold, i.e., by constraint \eqref{spc1}. Another consideration is the budget constraint in \eqref{spc8}. If the budget is insufficient, SPs $\tilde{k}$ cannot afford to acquire resources for serving requests. Since SPs $\tilde{k}$ are unable to reject requests, $\mathcal{S}^{\tilde{k}} = \varnothing$ and \textbf{HPP} is infeasible when the constraints \eqref{spc1} and \eqref{spc8} cannot be met.
\end{proof} 

% \rev{The primary goal of the proposed model is to determine the optimal pricing and activation decisions for the EC platform. The infeasibility of the \textbf{BIMP} does not yield meaningful insights, as it implies that the EC platform cannot provide edge resources to the SPs, leading to a scenario where the platform cannot generate any profit. In such instances, for meaningful insights, the EC platform can strategically set $\theta^k = 1$ for all $k \in \mathcal{K}$ to ensure that the \textbf{HPP} is feasible. Without loss of generality, we make the following assumption for the remainder of this paper.}

% \rev{\begin{assumption}
% \textbf{HPP} is feasible and $\Omega$ is nonempty.
% \end{assumption}}

\begin{remark}
It is noteworthy that, even when $\Omega$ is nonempty, the involvement of continuous variables in the upper-level decisions, coupled with binary variables in the lower level, can result in a scenario where the feasible set $\Omega$ is not closed.  Consequently, \textbf{BIMP} may not have an optimal solution. In such instances, the $\max$ in \textbf{BIMP} is replaced by $\sup$, and the solution approach proposed in the following section can be employed to obtain an $\epsilon$-optimal solution for the \textbf{BIMP}.
\end{remark}

For illustration, let's examine a simplified version of \textbf{BIMP}:
\begin{align*}
&\sup_{y,t} ~~~ t - y \\
&\text{subject to:}~~~~ 0\leq y \leq 1,\\
&\qquad \qquad \qquad t \in \argmin_t \left\{t: t\geq y, t\in \{0,1\} \right\}.
\end{align*}
This problem has the feasible set $(t,y)\in \Omega=\{(0,0)\} \cup (\{1\}\times(0,1])$, which is not closed. The supremum is $1$ and cannot be attained.

% \begin{proposition}\label{prop-HPPfeasible}
% The feasible set of the lower-level $\mathcal{S}^k(\bar{\bp},\bar{\bp}^s,\bar{\bz})$ is bounded, for all $k \in \mathcal{K}$. The feasible set of the upper-level $\Omega^{IR}$ is bounded. 
% \end{proposition}

% \begin{proof}
% For each $k\in \mathcal{K}$, the continuous variables are bounded: $0\leq x_{i,j}^k \leq \min\{R_i^k,C_j\}$, $0\leq y_j^k \leq C_j$ and $0\leq q_{i}^k \leq R_i^k$ for all $i\in \mathcal{I}$ and $j\in \mathcal{J}\cup \{0\}$. Additionally, the binary variables $\bt^k$ belong to a finite discrete set $\bT^k$ of cardinality $2^J$ with binary vectors of length $J$  which is a bounded set. Thus, $\mathcal{S}^k(\bar{\bp},\bar{\bp}^s,\bar{\bz})$ is bounded. For the second statement, we have that $p_j \in \{p_j^1 , \ldots, p_j^V \}$, $p^s_j \in 
% \{p_j^{s1} , \ldots, p_j^{sH} \}$ and $z_j \in \{0,1\}$, for all $j \in \mathcal{J} = \{1,\ldots,J\}$. $(\bar{\bp},\bar{\bp}^s,\bar{\bz})$ belong to a finite discrete set. Thus, $\Omega^{IR}$ is bounded.
% \end{proof}

To facilitate the convergence analysis, we consider the following assumption.
\begin{assumption} \label{asm-finiteOptSol}
\textbf{BIMP} has a finite optimal solution.
% $(\bp,\bp^s,\bz,\bx,\bq,\by,\bt) \in \Omega^{IR}$.
\end{assumption}

\section{Solution Approach}
\label{sol}
This section presents a solution to the challenging bi-level problem, \bmip, in which the leader aims to optimize a problem containing other optimization problems, i.e., the lower-level problems, within its constraints \cite{avra17}.
\textbf{BMIP} poses difficulties in solving due to the presence of integer variables in both levels, rendering the lower-level problems as non-convex.
In the absence of lower-level integer variables, the follower problems become convex, and one can employ KKT conditions or LP duality to transform each constraint in 
(\ref{enum}) into a set of linear integer equations. 
Subsequently, \textbf{BMIP} can be expressed as a MILP that can be solved by MILP solvers. % such as Gurobi, Cplex, or Mosek \cite{tnduong21}.

\subsection{Overview} \label{subsec:Overview}
First, we present a high-level overview of the approach used to solve the \bmip~problem. One straightforward method to handle the presence of integer variables in the follower problems is to enumerate all possible values for these variables in the lower-level problems. For each fixed set of values of the lower-level integer variables, we can transform each follower problem into an equivalent set of linear equations using KKT conditions or LP duality \cite{tnduong21}. Consequently, \textbf{BMIP} can be reformulated as a single-level MILP. However, this approach results in a large-scale MILP with an exponential number of constraints, making the single-level MILP formulation of \textbf{BMIP} intractable. To overcome this challenge, we propose an iterative decomposition algorithm that generates solutions sequentially within a finite number of iterations. This algorithm is based on the cutting plane generation method with reformulations and decomposition.

%In this work, we present an iterative algorithm for solving the proposed BiMIP model based on the cutting plane generation method using reformulations and decomposition. %\cite{bzeng14}. 
The core idea behind the proposed iterative algorithm is inspired by the  Column-and-Constraint Generation (CCG) approach proposed in \cite{zeng13} for solving two-stage RO.
The CCG approach first converts a two-stage robust problem into an equivalent MILP by enumerating all the extreme points of the uncertainty set.
%In CCG, a two-stage robust problem  can first be converted into an equivalent  MILP by enumerating all the extreme points of the uncertainty set. However, 
When the uncertainty set is extensive, this enumeration method produces a large-scale MILP that is not practically feasible to compute an optimal solution \cite{zeng13}. This situation is similar to the challenge we face in solving \textbf{BMIP}. As mentioned,  one approach to solving \textbf{BMIP} involves enumerating all possible integer values in the follower problems (\ref{enum}).
 %\textbf{BMIP} is through enumerating all possible integer values in the follower problems (\ref{enum}).
Nevertheless, this method is not scalable as we can have an enormous number of potential integer values in the follower problems. 

 %However, this enumeration method is not scalable as we can have a massive number of possible integer values in the follower problems. 
% ok;pokllfd
%To address this challenge, \cite{zeng13} proposes expanding a partial enumeration by gradually adding the most critical extreme points.
By combining the partial enumeration idea proposed for RO in \cite{zeng13}  with LP duality, as well as several reformulation and linearization techniques,
%By leveraging the partial enumeration idea proposed for robust optimization in \cite{zeng13} combined with several reformulation tricks,
we develop an efficient algorithm to solve the bi-level optimization problem \textbf{BMIP}  iteratively.
Specifically, instead of enumerating all possible integer variables in (\ref{enum}),
we use a partial enumeration formulation (i.e., a formulation defined over a subset of possible lower-level integer variables), which provides a valid relaxation and a UB to the original problem \textbf{BMIP}. 
%we consider a formulation based on a partial enumeration (i.e., a formulation defined over a subset of possible integer variables), which provides a valid relaxation, and consequently an UB, to the original maximization problem \textbf{BMIP}. 
% Then, by expanding the partial enumeration by gradually adding non-trivial integer variables in the follower problems, we obtain a stronger UB from solving the relaxed bi-level problem after each iteration. 
Gradually adding non-trivial integer variables in the follower problems to expand the partial enumeration allows us to obtain a stronger UB from solving the relaxed bi-level problem after each iteration. 
%Then, by gradually adding non-trivial integer variables in the follower problems to expand the partial enumeration, we obtain a stronger UB from solving the relaxed bi-level problem after each iteration. 
At each iteration, we also solve a subproblem to identify the most critical lower-level integer variables to include in the next iteration. The iterative algorithm operates in a master-subproblem framework until convergence. % Additionally, we show that this algorithm converges after a finite number of iterations to an exact global optimal solution.

In particular, to obtain optimal solutions ($\bx^{k*}$,~$\by^{k*}$,~$\bt^{k*}$) for the lower-level problems,
we utilize the solutions  $\bp^*$, $\bp^{s*}$ and $\bz^*$  acquired by solving the upper-level problem. 
As ($\bp^*$, $\bp^{s*}$, $\bz^*$, $\bx^{k*}$, $\by^{k*}$, $\bt^{k*}$) represents a feasible solution set of the original maximization problem, 
expanding the set $\bt^{k*}$  and evaluating the corresponding optimal value of the MP allows us to obtain stronger UBs. 
The algorithm also solves a set of subproblems at every iteration to compute an LB for \bmip. The solution to the subproblems provides a new constraint in the form of \eqref{eq:SubProbStrongDuality3} to be added to the MP after each iteration. 
By solving updated MPs and subproblems and improving the LB and UB in every iteration, this iterative algorithm converges to the optimal solution to the original bi-level problem \bmip~after a finite number of iterations.

% Solving updated MPs and subproblems iteratively leads to the convergence of this algorithm after a finite number of iterations to an exact global optimal solution to the original bi-level problem \bmip, as shown later.

\subsection{\bmip~with Duplicated Lower-Level Variables}
First, to facilitate the solution approach for solving \bmip, we propose a reformulation that results in a decomposable structure. To achieve this, we introduce a set of duplicated lower-level variables for the original problem \bmip, and obtain the following equivalent problem \bmipd.  For the detailed formulation of \bmipd, please refer to Appendix~\ref{app:bmipd}.
\begin{subequations}
\label{eq:bmipd}
\begin{alignat}{3}
% &\textbf{\bmipd:}~~ \max_{\substack{\bp,\bp^s,\bz,\\\bx^\prime,\bq^\prime,\by^\prime,\bt^\prime}} ~ \sum_j p_j \sum_k y_j^{\prime k} + \sum_j\sum_k \left(\phi_j^k+s^kp_j^s\right) t_j^{\prime k} \nonumber\\
% &\qquad\qquad\qquad- \sum_j \left(f_j z_j + c_j  \frac{\sum_k y_j^{\prime k}}{C_j} \right)   \label{objmd} \\
&\textbf{\bmipd:}~~ \max_{\substack{\bp,\bp^s,\bz,\bx^\prime,\bq^\prime,\by^\prime,\bt^\prime}} ~~~ \mathcal{P}\left(\bp,\bp^s,\bz,\by^\prime,\bt^\prime\right) \label{objmd} \\
&\text{subject to:}~~ (\bp,\bp^s,\bz,\bx',\bq',\by',\bt') \in \mathcal{H}, \label{eq-bmipdConstraints} \\
&(\bx^{\prime k}, \bq^{\prime k}, \by^{\prime k}, \bt^{\prime k}) \in \mathcal{S}^k(\bp,\bp^s,\bz),~ \forall k, \label{eq-dupConstraints}\\
% & \label{dupstart} \!\! t_j^{\prime k} \in \{0,1\}, \forall j,k;~~ x_{i,0}^{\prime k}, x_{i,j}^{\prime k}, q_{i}^{\prime k}, y_0^{\prime k}, y_{j}^{\prime k} \geq 0, \forall i,j,k,\\
% & \label{dupstart}\!\!\sum_k y_j^{\prime k} \leq z_j C_j; \sum_k s^k t_j^{\prime k} \leq z_jS_j; y_0^{\prime k} \geq \!\sum_i x_{i,0}^{\prime k}, \forall j,k,\!\\
% &p_0 y_0^{\prime k} \!+ \!\!\sum_j p_j y_j^{\prime k} \!+ \!\!\sum_j \left(\phi_j^k+s^kp_j^s\right) t_j^{\prime k}  \leq \!B^k, ~~\forall k, \\
% &t_j^{\prime k} \leq z_j; y_j^{\prime k} \geq \sum_i x_{i,j}^{\prime k};y_j^{\prime k} \leq C_j t_j^{\prime k};s^k t_j^{\prime k} \leq S_j,\forall j,k,\! \label{bilinmp1}\\
% &x_{i,j}^{\prime k} \leq a_{i,j}^k R_i^k;~\sum_j x_{i,j}^{\prime k} + x_{i,0}^{\prime k} + \theta^k q_i^{\prime k} = R_i^k, ~\forall i,k, \\
% &x_{i,0}^{\prime k} d_{i,0} +  \sum_j x_{i,j}^{\prime k} d_{i,j} \leq D^{k, \sf m} R_i^k, ~~\forall i,k, \label{dupend}\\
% \label{enumx}
% & w^k \Big( \sum_i x_{i,0}^{\prime k} d_{i,0} + \sum_{i,j} x_{i,j}^{\prime k} d_{i,j} \Big) + p_0 y_0^{\prime k} + \sum_j p_j y_j^{\prime k}\nonumber\\
% & + \theta^k \sum_i \psi_i^k q_i^{\prime k}+ \sum_j (\phi_j^k   +s^kp_j^s) t_j^{\prime k} \leq  \min_{\substack{(\bx^k, \bq^k, \by^k, \bt^k)\\ \in \mathcal{S}^k(\bp,\bp^s,\bz)}} \bigg\{p_0 y_0^k \nonumber\\
% & + \sum_j p_j y_j^k +\sum_j \left(\phi_j^k+s^kp_j^s\right) t_j^k    + \theta^k \sum_i \psi_i^k q_i^k  \nonumber\\
% & +  w^k \Big( \sum_i x_{i,0}^k d_{i,0} + \sum_{i,j} x_{i,j}^k d_{i,j} \Big) \bigg\},~ \forall k.
\label{enumx}
& \mathcal{C}^k(\bx^{\prime k},\bq^{\prime k},\by^{\prime k},\bt^{\prime k}|\bp,\bp^s,\bz) \nonumber\\
& \qquad\leq  \min_{\substack{(\bx^k, \bq^k, \by^k, \bt^k)\\ \in \mathcal{S}^k(\bp,\bp^s,\bz)}} \mathcal{C}^k(\bx^k,\bq^k,\by^k,\bt^k|\bp,\bp^s,\bz),~ \forall k.
\end{alignat}
\end{subequations}

\begin{proposition}
The formulation \bmip~ in \eqref{eq:bmip} is equivalent to the reformulation \bmipd~ in \eqref{eq:bmipd}. Furthermore, the feasible set of \bmipd~is $\Omega^{IR}$, the inducible region of \bmip.
\end{proposition}
\begin{proof}
The satisfaction of \eqref{eq-dupConstraints} and \eqref{enumx} ensures that $(\bx^{\prime k}, \bq^{\prime k}, \by^{\prime k}, \bt^{\prime k})$ constitutes an optimal solution to the lower-level problem $(\Xi^k)$ of SP $k$ for all $k \in \mathcal{K}$. Consequently, the direct equivalence between $\bmipd$ and $\bmip$ is evident. Furthermore, given $(\bp, \bp^s, \bz)$, an optimal solution to the lower-level problem can be obtained by setting $(\bx^k, \bq^k, \by^k, \bt^k) = (\bx^{\prime k}, \bq^{\prime k}, \by^{\prime k}, \bt^{\prime k})$ for all $k \in \mathcal{K}$. The second statement follows immediately.
\end{proof}

The set $(\bx^{\prime k}, \bq^{\prime k}, \by^{\prime k}, \bt^{\prime k})$ represents the duplicated lower-level variables $(\bx^{k}, \bq^{k}, \by^{k}, \bt^{k})$. At first glance, the introduction of duplicated variables and the replacement of $\bmip$ with $\bmipd$ may seem counterintuitive, as $\bmipd$ involves a larger set of variables. However, analyzing the original bi-level formulation $\bmip$ and gaining insights prove challenging due to the embedded nature of the lower-level problem within the constraints of the upper-level problem. The nested structure introduces intricate interdependencies between upper and lower decisions, significantly increasing the complexity of the optimization process. This interdependency is particularly pronounced in our paper due to the presence of integer variables at both levels.

\begin{remark}
By duplicating the lower-level variables and constraints, the reformulation $\bmipd$ provides a comprehensive variable set $(\bp, \bp^s, \bz, \bx', \bq', \by', \bt')$, encompassing both the original upper-level variables $(\bp, \bp^s, \bz)$ and the lower-level variables now represented by $(\bx', \bq', \by', \bt')$, under the control of the leader. This extension empowers the leader to utilize the duplicated variables   $(\bx', \bq', \by', \bt')$ for simulating follower responses across various upper-level decisions.
\end{remark}

We are now ready to present the  iterative decomposition algorithm for solving \bmipd ~for optimal
pricing design, resource management, and service placement in a master-subproblem framework. First, for each value of $\bt^k$, we can convert each follower problem constraint in \eqref{enumx} into an equivalent set of linear equations.  Then, \bmipd ~can be expressed as a large-scale MILP. The number of constraints in the resulting MILP is proportional to the number of possible values of $\bt^k$, which can be exponential. Therefore, we further develop an efficient algorithm to solve this MILP iteratively. In each iteration, we only need to solve a relaxed MILP with a subset of constraints in the original MILP and gradually add critical constraints to the relaxed MILP until convergence.

\subsection{Reformulation of the Follower Problems}
Solving the lower-level problems is challenging due to the presence of discrete variables. To this end, we propose a series of transformations combined with LP duality to convert each follower problem in \eqref{enumx} into a set of equivalent linear equations. 
Specifically, for each follower problem, we can separate discrete and continuous variables  as follows:
% \begin{align}\label{sepenum}
% &w^k \Big( \sum_i x_{i,0}^{\prime k} d_{i,0} + \sum_{i,j} x_{i,j}^{\prime k} d_{i,j} \Big) + \!\theta^k \!\sum_i \psi_i^k p_i^{\prime k}  +  \sum_j p_j y_j^{\prime k} \nonumber\\
% &+ p_0 y_0^{\prime k}+ \sum_j \left(\phi_j^k+s^kp_j^s\right) t_j^{\prime k} \leq \underset{\bt^k \in  \bT^k}{\text{min}}   \Bigg\{ \sum_{j} \left(\phi_j^k+s^kp_j^s\right) t_j^k  \nonumber\\
% &+\underset{(\bx^k\!,\by^k\!,\bq^k) \in \mathcal{S}^k(\bp,\bp^s,\bz,\bt^k)}{\text{min}} \bigg[ p_0 y_0^k+ \!\!\sum_j p_j y_j^k \!  + \theta^k \!\sum_i \!\psi_i^k q_i^k\nonumber \\ 
%  &+ w^k \Big(\!\sum_i x_{i,0}^k d_{i,0} \!+ \!\!\sum_{i,j} x_{i,j}^k d_{i,j} \Big)\! \bigg] \Bigg\}, \forall k,
% \end{align}
% \rev{\begin{align}\label{sepenum}
% &\mathcal{C}^k(\bx^{\prime k},\bq^{\prime k},\by^{\prime k},\bt^{\prime k}|\bp,\bp^s,\bz) \leq \underset{\bt^k \in  \bT^k}{\text{min}}   \Bigg\{ \sum_{j} \left(\phi_j^k+s^kp_j^s\right) t_j^k  \nonumber\\
% &+\underset{(\bx^k\!,\by^k\!,\bq^k) \in \mathcal{S}^k(\bp,\bp^s,\bz,\bt^k)}{\text{min}} \bigg[ p_0 y_0^k+ \!\!\sum_j p_j y_j^k \!  + \theta^k \!\sum_i \!\psi_i^k q_i^k\nonumber \\ 
%  &+ w^k \Big(\!\sum_i x_{i,0}^k d_{i,0} \!+ \!\!\sum_{i,j} x_{i,j}^k d_{i,j} \Big)\! \bigg] \Bigg\}, \forall k,
% \end{align}}
\begin{align}\label{sepenum}
&\mathcal{C}^k(\bx^{\prime k},\bq^{\prime k},\by^{\prime k},\bt^{\prime k}|\bp,\bp^s,\bz) \leq \underset{\bt^k \in  \bT^k}{\text{min}}   \Bigg\{ \sum_{j} \left(\phi_j^k+s^kp_j^s\right) t_j^k  \nonumber\\
&+\underset{(\bx^k\!,\by^k\!,\bq^k) \in \mathcal{S}^k(\bp,\bp^s,\bz,\bt^k)}{\text{min}} \bigg[ p_0 y_0^k+ \!\!\sum_j p_j y_j^k \!  +  \!\sum_i \!\psi_i^k q_i^k\nonumber \\ 
 &+ w^k \Big(\!\sum_i x_{i,0}^k d_{i,0} \!+ \!\!\sum_{i,j} x_{i,j}^k d_{i,j} \Big)\! \bigg] \Bigg\}, \forall k,
\end{align}
%\noindent
where $\bT^k$ represents the collection of all possible $\bt^k$ for the lower-levels of the problem, and 
% $\mathcal{S}^k(\bp,\bp^s,\bz,\bt^k)$
$\mathcal{S}^k(\bp,\bp^s,\bz,\bt^k)=\{(\bx^{ k}, \bq^{ k}, \by^{ k})|~(\ref{spc1})-(\ref{spc10}),~ (\ref{spc11}),~(\ref{spc12})\}$ 
is the set $\mathcal{S}^k(\bp,\bp^s,\bz)$ for a given $\bt^k$. 
Let $\mu_1^k \ge 0$, $\nu_j^k\ge 0$, $\Gamma_j^k\ge 0$, $\mu_2^k\ge 0$, $\sigma_j^k\ge 0$, $\xi_i^k\ge 0$, $\eta_i^k \in \mathbb{R}$, and $\tau_{i,j}^k\ge 0$ be continuous dual variables associated with the constraints in the feasible set  $\mathcal{S}^k(\bp,\bp^s,\bz,\bt^k)$ (see Appendix~\ref{app-FeasibleSetSk}). By employing LP duality in the inner minimization problem in \eqref{sepenum}, the lower-level problems become:
\begin{align} 
% &w^k \Big(\sum_i x_{i,0}^{\prime k} d_{i,0} + \sum_{i,j} x_{i,j}^{\prime k} d_{i,j}\Big) + \!\theta^k \!\sum_i \psi_i^k p_i^{\prime k} + \sum_j p_j y_j^{\prime k} \nonumber\\ 
% &+p_0 y_0^{\prime k}+\sum_j \left(\phi_j^k+s^kp_j^s\right) t_j^{\prime k}
% \leq \min_{\bt^k \in  \bT^k} \sum_j \left(\phi_j^k+s^kp_j^s\right) t_j^k\nonumber\\
&\mathcal{C}^k(\bx^{\prime k},\bq^{\prime k},\by^{\prime k},\bt^{\prime k}|\bp,\bp^s,\bz) \leq \min_{\bt^k \in  \bT^k} \sum_j \left(\phi_j^k+s^kp_j^s\right) t_j^k\nonumber\\
&+ \!\bigg\{  \max_{\substack{\mu_1^k,\nu^k_j,\Gamma^k_j,\mu_2^k,\\\sigma^k_j,\xi^k_i,\eta^k_i,\tau^k_{i,j}}}  -\mu_1^k \Big(B^k \!- \!\!\sum_j\! \left(\phi_j^k+s^kp_j^s\right)t_j^k \Big)+\!\sum_i R_i^k \eta_i^k  \nonumber\\
& + \nu_j^k(t_j^k\!-\!z_j) - \!\sum_i\! R_i^k D^{k,m} \xi_i^k  - \!\sum_j\! C_j t_j^k \sigma_j^k - \!\sum_{i,j}\! a_{i,j}^k \tau_{i,j}^k R_i^k \nonumber \\
&  \text{subject to:} \quad  p_0 \big( 1 + \mu_1^k \big) - \mu_2^k \geq 0,   \nonumber \\ 
&  \quad \quad   \quad  \! p_j \mu_1^k -\Gamma_j^k + \sigma_j^k + p_j \geq 0,~\forall j, \nonumber \\
% &   \quad  \quad    \quad    \!  \!
% \mu_2^k + d_{i,0} \xi_i^k - \eta_i^k \geq -w^k d_{i,0},~~ \theta^k\eta_i^k \leq \theta^k\psi_i^k, ~~ \forall i, \nonumber \\
&   \quad  \quad    \quad    \!  \!
\mu_2^k + d_{i,0} \xi_i^k - \eta_i^k \geq -w^k d_{i,0},~~ \eta_i^k \leq \psi_i^k, ~~ \forall i, \nonumber \\
&   \quad  \quad  \quad  \!
\Gamma_j^k + d_{i,j} \xi_i^k + \tau_{i,j}^k - \eta_i^k \geq -w^k d_{i,j}~,\forall i,j, \nonumber \\
&   \quad  \quad    \quad  \!
\mu_1^k, ~\nu_j^k, ~\Gamma_j^k, ~\mu_2^k, ~\sigma_j^k, ~\xi_i^k,~\tau_{i,j}^k \geq 0,~\forall i,j   \!\bigg\}, \forall k. \!\!\!\! \label{eq:SubProbStrongDuality1}
\end{align}

Recall that $\bt^k = (t_1^k,\ldots,t_J^k), \forall k$, where each $t_j^k$ is a binary variable. Thus, $\bT^k = (\bt^{k,1}, \ldots, \bt^{k,l},\ldots, \bt^{k,Q})$ where $Q = | \bT^k | =  2^J$ and  $\bt^{k,l}$ is a $J-$dimensional vector with binary entries. Then, \eqref{eq:SubProbStrongDuality1} is equivalent to:
\begin{subequations}
\label{eq:SubProbStrongDuality2}
\begin{alignat}{3} 
% &\!\!\!w^k \Big(\sum_i x_{i,0}^{\prime k} d_{i,0} + \sum_{i,j} x_{i,j}^{\prime k} d_{i,j}\Big) + \!\theta^k \!\sum_i \psi_i^k p_i^{\prime k} +p_0 y_0^{\prime k}  \nonumber\\ 
% &\!\!\!+ \sum_j p_j y_j^{\prime k}+\sum_j \left(\phi_j^k+s^kp_j^s\right) t_j^{\prime k}
% \leq  \sum_j \left(\phi_j^k+s^kp_j^s\right) t_j^{k,l} \nonumber\\
&~~\mathcal{C}^k(\bx^{\prime k},\bq^{\prime k},\by^{\prime k},\bt^{\prime k}|\bp,\bp^s,\bz)
\leq  \sum_j \left(\phi_j^k+s^kp_j^s\right) t_j^{k,l} \nonumber\\
&\!+ \!\Bigg\{  \max_{\substack{\mu_1^{k,l},\nu_j^{k,l},\Gamma_j^{k,l},\mu_2^{k,l},\\\sigma_j^{k,l},\xi_i^{k,l},\eta_i^{k,l},\tau_{i,j}^{k,l}}}  -\mu_1^{k,l} \Big(B^k - \sum_j \left(\phi_j^k+s^kp_j^s\right)t_j^{k,l} \Big)  \nonumber\\
&  \qquad + \nu_j^{k,l}(t_j^{k,l}-z_j) + \sum_i R_i^k \eta_i^{k,l} - \sum_i R_i^k D^{k,m} \xi_i^{k,l}  \nonumber \\
& \qquad - \sum_j C_j t_j^{k,l} \sigma_j^{k,l} - \sum_{i,j} a_{i,j}^k \tau_{i,j}^{k,l} R_i^k  \label{eq:SubProbStrongDualityConstObj1} \\
&   \text{subject to:} \quad p_0 \big( 1 + \mu_1^{k,l} \big) - \mu_2^{k,l} \geq 0,  \label{eq:SubProbStrongDualityConst1}  \\ 
&       p_j \mu_1^{k,l} -\Gamma_j^{k,l} + \sigma_j^{k,l} + p_j \geq 0,~\forall j, \\
&    
\mu_2^{k,l} + d_{i,0} \xi_i^{k,l} - \eta_i^{k,l} \geq -w^k d_{i,0},~\forall i\\
% &\theta^k\eta_i^{k,l} \leq \theta^k\psi_i^k,~\forall i, \\
% &    
&\eta_i^{k,l} \leq \psi_i^k,~\forall i, \\
&    
\Gamma_j^{k,l} + d_{i,j} \xi_i^{k,l} + \tau_{i,j}^{k,l} - \eta_i^{k,l} \geq -w^k d_{i,j},~\forall i,j, \\
&      
\mu_1^{k,l}\!, \Gamma_j^{k,l}\!, \mu_2^{k,l}\!, \sigma_j^{k,l}\!, \xi_i^{k,l}\!, \tau_{i,j}^{k,l} \geq 0,\forall i,j,  \label{eq:SubProbStrongDualityConste} \!\!\Bigg\}, \forall k,  1 \!\leq\! l \leq Q.\!\!\!
\end{alignat}
\end{subequations}

We can remove its max operator to obtain the following:
%\begin{align} \label{eq:SubProbStrongDuality3}
\begin{subequations}
 \label{eq:SubProbStrongDuality3}
\begin{alignat}{3}
% &\bigg\{ w^k \Big( \sum_i x_{i,0}^{\prime k} d_{i,0} + \sum_{i,j} x_{i,j}^{\prime k} d_{i,j} \!\Big) + p_0 y_0^{\prime k} +\sum_j p_j y_j^{\prime k}\nonumber\\
% & + \sum_j  \left(\phi_j^k+s^kp_j^s\right) t_j^{\prime k}\leq   \sum_j  \left(\phi_j^k+s^kp_j^s\right) t_j^{k,l} \!+\!\sum_i\! R_i^k \eta_i^{k,l} \nonumber\\
&~~\Big\{\mathcal{C}^k(\bx^{\prime k},\bq^{\prime k},\by^{\prime k},\bt^{\prime k}|\bp,\bp^s,\bz) \leq   \sum_j  \left(\phi_j^k+s^kp_j^s\right) t_j^{k,l} \!\nonumber\\
&+\!\sum_i\! R_i^k \eta_i^{k,l} \!\!-\!\mu_1^{k,l} \!\Big(\!B^k\!\! -\!\! \sum_j  \!\left(\phi_j^k+s^kp_j^s\right) t_j^{k,l} \Big) \!+\! \nu_j^{k,l}(t_j^{k,l}\!-\!z_j)  
\nonumber\\
&- \!\sum_j\! C_j t_j^{k,l} \sigma_j^{k,l}-\sum_i R_i^k D^{k,m} \xi_i^{k,l}  \!- \!\sum_{i,j} a_{i,j}^k \tau_{i,j}^{k,l} R_i^k,\! \label{eq:SubProbStrongDualityObj} \\
& \quad \eqref{eq:SubProbStrongDualityConst1}- \eqref{eq:SubProbStrongDualityConste}
 ~~\Big\},~\forall k,  1\leq l \leq Q.\!\!\!
%\end{align}
\end{alignat}
\end{subequations}

\begin{proposition} \label{prop-BMIPeEquivalence}
Let Assumption~\ref{asm-finiteOptSol} hold. The formulation \bmipd~ in \eqref{eq:bmipd} is equivalent to the following expanded single-level formulation:
\begin{align*}
&\textbf{\ebmipd :} ~~\max_{\substack{\bp,\bp^s,\bz,\bx^\prime,\bq^\prime,\by^\prime,\bt^\prime}} ~~~ \mathcal{P}\left(\bp,\bp^s,\bz,\by^\prime,\bt^\prime\right) \\
& \textnormal{subject to:} ~~~ \eqref{eq-bmipdConstraints}, ~\eqref{eq-dupConstraints}, ~\eqref{eq:SubProbStrongDuality3}.
\end{align*}
\end{proposition}

\begin{proof}
We first note that \ebmipd~ is derived from \bmipd~ by reformulating the constraint \eqref{enumx}. We will examine each reformulation step to establish its equivalence.\\
\indent \textbf{Step 1:} \eqref{enumx}$\Leftrightarrow$\eqref{sepenum} This is straightforward since we only restructure the RHS of the constraint by separating the discrete and continuous variables.\\
\indent \textbf{Step 2:} \eqref{sepenum}$\Leftrightarrow$\eqref{eq:SubProbStrongDuality1} due to LP duality and Assumption~\ref{asm-finiteOptSol}.\\
\indent \textbf{Step 3:} \eqref{eq:SubProbStrongDuality1} $\Leftrightarrow$ \eqref{eq:SubProbStrongDuality2} since $\bT^k = (\bt^{k,1}, \ldots, \bt^{k,l},\ldots, \bt^{k,Q})$.\\
\indent \textbf{Step 4:} We remove the $\max$ operator in the constraint \eqref{eq:SubProbStrongDuality2} to obtain \eqref{eq:SubProbStrongDuality3}. The resulting maximization problem $\{(\mathcal{M}_2): \eqref{objmd} - \eqref{eq-dupConstraints}, \eqref{eq:SubProbStrongDuality3}\}$  is equivalent to $\{(\mathcal{M}_1): \eqref{objmd} - \eqref{eq-dupConstraints}, \eqref{eq:SubProbStrongDuality2}\}$. First, the optimal solution of $(\mathcal{M}_2)$ is also the optimal solution of $(\mathcal{M}_1)$. This is because the set of variables satisfying \eqref{eq:SubProbStrongDuality3} will also satisfy \eqref{eq:SubProbStrongDuality2} since RHS\eqref{eq:SubProbStrongDualityObj} $\leq$ RHS\eqref{eq:SubProbStrongDualityConstObj1}  . Alternatively, the optimal solution $(\bp^{*1}$, $\bp^{s*1}$, $\bz^{*1}$, $\bx^{\prime *1}$, $\bq^{\prime *1}$, $\by^{\prime *1}$, $\bt^{\prime *1})$ of $(\mathcal{M}_1)$ is also the optimal solution of $(\mathcal{M}_2)$. This is demonstrated by considering two cases:\\
\textbf{Case 1.} The constraint \eqref{eq:SubProbStrongDualityObj} is binding, and the optimal solutions of the maximization problem on RHS\eqref{eq:SubProbStrongDualityObj} is $(\mu^{*1}$, $\nu^{*1}$, $\Gamma^{*1}$, $\sigma^{*1}$, $\xi^{*1}$, $\eta^{*1}$, $\tau^{*1})$. Since RHS\eqref{eq:SubProbStrongDualityObj} $\leq$ RHS\eqref{eq:SubProbStrongDualityConstObj1}, the optimal value of $(\mathcal{M}_2)$ is less than or equal to that of $(\mathcal{M}_1)$. Equality can be attained with the optimal solution to $(\mathcal{M}_2)$ chosen as $(\bp^{*1}$, $\bp^{s*1}$, $\bz^{*1}$, $\bx^{\prime *1}$, $\bq^{\prime *1}$, $\by^{\prime *1}$, $\bt^{\prime *1})$ and $(\mu^{*1}$, $\nu^{*1}$, $\Gamma^{*1}$, $\sigma^{*1}$, $\xi^{*1}$, $\eta^{*1}$, $\tau^{*1})$ which is feasible for \eqref{eq:SubProbStrongDuality2}.\\
\textbf{Case 2.} The constraint \eqref{eq:SubProbStrongDualityObj} is not binding. The optimal solutions to $(\mathcal{M}_2)$ are $(\bp^{*1}$, $\bp^{s*1}$, $\bz^{*1}$, $\bx^{\prime *1}$, $\bq^{\prime *1}$, $\by^{\prime *1}$, $\bt^{\prime *1})$ and $(\mu^{*2}, \nu^{*2}, \Gamma^{*2}, \sigma^{*2}, \xi^{*2}, \eta^{*2}, \tau^{*2})$, and the optimal values of $(\mathcal{M}_1)$ and $(\mathcal{M}_2)$ equal. However, the RHS of \eqref{eq:SubProbStrongDualityConstObj1} at $(\mu^{*2}, \nu^{*2}, \Gamma^{*2}, \sigma^{*2}, \xi^{*2}, \eta^{*2}, \tau^{*2})$ is strictly less than the RHS of \eqref{eq:SubProbStrongDualityObj} at $(\mu^{*1}$, $\nu^{*1}$, $\Gamma^{*1}$, $\sigma^{*1}$, $\xi^{*1}$, $\eta^{*1}$, $\tau^{*1})$.\\
\indent In either cases, the optimal solutions $(\bp^{*}$, $\bp^{s*}$, $\bz^{*}$, $\bx^{\prime *}$, $\bq^{\prime *}$, $\by^{\prime *}$, $\bt^{\prime *})$ are identical for both $(\mathcal{M}_1)$ and $(\mathcal{M}_2)$ with the same optimal values, demonstrating the equivalence of the two reformulations. Proposition~\ref{prop-BMIPeEquivalence} follows accordingly.
\end{proof}

\begin{remark}
Recalling that for a given $(\bp,\bp^s,\bz,\bt^k)$, the remaining follower problem, identified as the inner minimization problem in \eqref{sepenum}, is an LP referred to in this proof as the primal problem $(P^k)$. The corresponding dual problem, denoted as $(D^k)$, is the inner maximization problem in \eqref{eq:SubProbStrongDuality1}.\\
\indent Assumption~\ref{asm-finiteOptSol} ensures the existence of a specific optimal pair of primal variables $(\bx, \by, \bq)$ associated with dual variables $(\mu, \nu, \Gamma, \sigma, \xi, \eta, \tau)$. Consequently, the feasible set of the dual problems $(D^k)$ is nonempty for all $k$. In other words, the dual problem $(D^k)$ is either finitely optimal or unbounded.\\
\indent By the Duality Theorem, two possibilities emerge:
\begin{enumerate}
    \item Both $(P^k)$ and $(D^k)$ are finitely optimal. In this case, a non-trivial upper bound parameterized by $(\bp, \bp^s, \bz)$ exists for the RHS of \eqref{eq:SubProbStrongDualityObj}.
    \item $(P^k)$ is infeasible, and $(D^k)$ is unbounded (as $(D^k)$ is always feasible, as shown above). In this case, the upper bound on the RHS of \eqref{eq:SubProbStrongDualityObj} is trivial, being $+\infty$.
\end{enumerate}
\end{remark}

\subsection{Single-Level MILP Reformulation of \bmip}
\label{singlemilp}
By replacing constraints \eqref{enum} in the original bi-level problem \bmip ~by \eqref{eq:SubProbStrongDuality3}, we obtain an equivalent single-level optimization problem \ebmipd~(refer to Proposition~\ref{prop-BMIPeEquivalence}).  
However, the resulting problem contains bilinear terms $p_j^s \mu_1^{k,l}$, $\nu_j^{k,l}z_j$, $p_j \mu_1^{k,l}$ and $p_j y_j^{\prime k}$ in \eqref{eq:SubProbStrongDuality3}, as well as the bilinear terms $p_j^s t_j^{\prime k}$ in \eqref{eq:SubProbStrongDuality3} and the objective function \eqref{objmd}. Note that $t_j^{k,l}$ is a known parameter in \eqref{eq:SubProbStrongDuality3}. 
By utilizing  (\ref{ipc5}), we can express the bilinear terms $p_j^s \mu_1^{k,l}$ in the following  form:
\beqn
\label{psmu1linear}
p_j^s \mu_1^{k,l} = \sum_h p_j^{sh} r_j^{sh} \mu_1^{k,l} = \sum_h p_j^{sh}  \kappa_j^{h,k,l},
\eeqn
% \label{ipc5}
% &p_j = \sum_{v=1}^V p_j^v r_j^v,\forall j; 
% &&\sum_{v=1}^V r_j^v = 1,\forall j; &&r_j^v\in\{0,1\}, \forall j, v,\\
% \beqn
% \label{psmu1linearccg}
% p^s_j \mu_1^{k,l} = \sum_h p_j^{sh} r_j^{sh} \mu_1^{k,l} = \sum_h p_j^{sh} \kappa_j^{h,k,l},
% \eeqn
where $\kappa_j^{h,k,l} = r_j^{sh} \mu_1^{k,l}$, which is a product of a binary variable and a non-negative continuous variable. 
Let $b$ and $u$ be a binary variable and a non-negative continuous variable, respectively. Consider the bilinear term $U = ub$. %between the binary variable $b$ and the continuous variable $u$
We have $U = u$ if $b = 1$ and $U = 0$ if $b = 0$. Thus, the bilinear term $ub$ can be implemented through the following linear equations \cite{linear2}:
\begin{subequations}
\label{bilinearBC}
\begin{align}
\label{bilinearBC1}
U \leq M b;~~ U \leq u + M - M  b, \\
\label{bilinearBC2}
U \geq 0;~~ U \geq u + M b - M,
\end{align}
\end{subequations}
where $M$ is a sufficiently large number. By applying the linearization procedures in \eqref{bilinearBC} to $\kappa_j^{h,k,l} = r_j^{sh} \mu_1^{k,l}$ in \eqref{psmu1linear}, the bilinear term $p_j^s \mu_1^{k,l}$ can be rewritten as a set of linear equations, as follows: 
\begin{subequations}\label{psmulinearccg}
\begin{align}
\label{psmulinear1ccg}
\kappa_j^{h,k,l} \leq M r_j^{sh};~
\kappa_j^{h,k,l} \leq \mu_1^{k,l},~\forall j, k, h,l,\\
\label{psmulinear2ccg}
\kappa_j^{h,k,l} \geq 0;~ \kappa_j^{h,k,l} \geq \mu_1^{k,l} + M r_j^{sh} - M, ~ \forall j, k, h,l.
\end{align}
\end{subequations}
Similarly, we can linearize bilinear terms $p_j^s \mu_1^{k,l}$, $p_j y_j^{\prime k}$ and $p_j^s t_j^{\prime k}$ using (\ref{ipc56}), and \eqref{bilinearBC}. Finally, bilinear terms $\nu_j^{k,l}z_j$ are products of a discrete and binary variable, which can be linearized using \eqref{bilinearBC}. Refer to Appendix~\ref{app-BilinearLinearization} for further details.

\begin{remark}
% The linearization process utilizing the Big-M method is, by nature, an approximation and does not inherently ensure a strict equivalence between the linearized model \bmipl~and the original problem \ebmipd. 
Choosing an appropriate value for the Big-M parameter, which ensures that the linearized constraints are as tight as possible without making them infeasible, is essential for achieving an accurate linearized model.
\end{remark}

By employing the linearization procedures described above for the bilinear terms in \ebmipd, \bmip~ is transformed into the equivalent problem \bmipl, which is indeed an MILP. However, \bmipl ~contains $K Q = K 2^J $ constraints in the form of \eqref{eq:SubProbStrongDuality3}. Consequently, the problem size grows exponentially with the number of ENs. Inspired by the CCG algorithm in RO, we propose to address the resulting large-scale MILP \bmipl~ iteratively within a master-subproblem framework. For a detailed description of the proposed master-subproblem framework, please refer to Section~\ref{subsec:Overview}.

% The MP represents a reduced form of \bmipl ~with only a subset of constraints in \eqref{eq:SubProbStrongDuality3}. Thus, solving the MP provides an UB for the original bi-level maximization problem \bmip.  The algorithm also solves a subproblem at every iteration to compute a LB for \bmip. The solution to the subproblem provides a new constraint in form of \eqref{eq:SubProbStrongDuality3} to be added to the MP after each iteration. \rev{By improving the LB and UB every iteration, we anticipate that this iterative algorithm converges to the optimal solution to the original bi-level problem \bmip.}

% Solving the updated MP helps improve the UB. We shows that this iterative algorithm converges to an exact optimal solution to the original bi-level problem \bmip. In the following, we  present the MP, the subproblem, and the iterative algorithm to solve \bmipl.

% We note that this reformulation, as compared to the reformulation based on KKT conditions, has simpler structure, fewer variables and constraints and is usually far less computationally expensive. Indeed, it can be estimated by professional MIP solvers which are conveniently accessible in practice. 

% In the following subsection, we will employ this reformulation approach to develop a master-subproblem framework for our algorithm.

\subsection{Master Problem}
\label{masterccg}
The \textbf{MP} is a relaxed problem of \bmipl ~that is equivalent to \bmip. At iteration $L$, the \textbf{MP} consists of $L$ constraints in \eqref{eq:SubProbStrongDuality3}. 
The \textbf{MP} at iteration $L$  is given by:
\begin{subequations}
\label{MPfinal}
\begin{alignat}{3}
% & \textbf{MP}:~~ % \underline{\Theta}
% \Theta ~ =\max_{\substack{\bp,\bp^s,\bz,\\\bx^\prime,\bq^\prime,\by^\prime,\bt^\prime}} ~ \sum_{j,k}  \Bigg[p_j y_j^{\prime k} +  \left(\phi_j^k+s^kp_j^s\right) t_j^{\prime k}\Bigg] \nonumber\\
% & \quad \quad \quad \quad \quad \quad \quad \quad - \sum_j \left(f_j z_j + c_j  \frac{\sum_k y_j^{\prime k}}{C_j} \right)  \label{objmpd} \\
& \textbf{MP}:~~ \Theta ~ =\max_{\substack{\bp,\bp^s,\bz,\bx^\prime,\bq^\prime,\by^\prime,\bt^\prime}} ~~~ \mathcal{P}\left(\bp,\bp^s,\bz,\by^\prime,\bt^\prime\right) \label{objmpd} \\
% & \text{subject to:} ~  (\bp,\bp^s,\bz,\bx',\bq',\by',\bt') \in \mathcal{H}, ~ \eqref{dupstart}-\eqref{dupend}, \!\!\\
& \text{subject to:} ~~~ \eqref{eq-bmipdConstraints}, ~\eqref{eq-dupConstraints}, \!\!\\
&\quad \quad \quad \quad \quad \Big\{   \eqref{eq:SubProbStrongDualityConst1}- \eqref{eq:SubProbStrongDualityConste},~ \eqref{eq:SubProbStrongDualityObj} \Big\}, ~~\forall k, 1 \!\leq\! l \!\leq\! L.
%&z_j\in \{0,1\};~~ t_j^{\prime k} \in \{0,1\}; ~~ \eta_i^{k,l} \in\mathbb{R}
% &p_j = \sum_v p_j^v r_j^v,\forall j;
% \sum_v r_j^v = 1,\forall j; r_j^v\in\{0,1\}, ~~\forall j, v\\
% &p^s_j = \sum_h p_j^{sh} r_j^{sh},\forall j;
% \sum_h r_j^{sh} = 1,\forall j; r_j^{sh}\in\{0,1\}, \forall j, h.
\end{alignat}
\end{subequations}
 
% For the bilinear terms $p_j^s t_j^k$ in the objective function of \textbf{MP}, discretizing $p_j^s$ using (\ref{ipc6}), we have,
% \beqn
% \label{pstlinearccg}
% p^s_j t_j^k = \sum_l p_j^{sl} r_j^{sl} t_j^k = \sum_l p_j^{sl} \zeta_j^{l,k},
% \eeqn
% where $\zeta_j^{l,k}$ are products of binary variables $r_j^{sl}$ and $t_j^k$. Hence, the constraints $\zeta_j^{l,k}=r_j^{sl} t_j^k$ can be implemented through the following linear inequalities, for all $l,j,k$:
% \beqn
% \label{pmulinear1ccgz}
% \zeta_j^{l,k} \leq r_j^{sl};~~
% \zeta_j^{l,k} \leq t_j^k; ~~
% \zeta_j^{l,k} \geq r_j^{sl}+t_j^k-1.
% \eeqn
After the linearization steps in Section \ref{singlemilp}, the \textbf{MP} is a MILP that can be solved by MILP solvers. Refer to Appendix~\ref{app-singlemilp-dual} for the MILP \bmipl ~of the \textbf{MP}. Since \textbf{MP} is a relaxed version of the original \bmip, the optimal value of \textbf{MP} provides a valid UB for \bmip. 
% In \textbf{MP}, duplicated lower-level variables (and constraints) are denoted in the upper-level by $x^{\prime k}$, $y^{\prime k}$, and $t^{\prime k}$ to provide a complete variable set which is under the control of the leader. Since each new iteration adds more constraints to the master problem, we obtain a tighter LB. Thus:
Therefore:
\beqn \label{eq-UB}
UB = \Theta^*.   
\eeqn
Since each new iteration adds more constraints to the MP, we obtain a tighter UB after each iteration. 

\subsection{Subproblems}
We formulate and evaluate the following subproblems, $\textbf{SP1}^k$, for a given upper-level resource selling price and EN activation decision $(\bp^*\!,\bp^{s*}\!,\bz^*)$. Each subproblem $\textbf{SP1}^k$ corresponds to the lower-level problem for service $k$.
% \begin{subequations}
% \label{SP1}
% \begin{alignat}{3}
% &\!\textbf{SP1}^k \!\!: \varphi^k(\bp^*\!,\!\bp^{s*}\!\!,\!\bz^*)\!= \!\! \underset{\bx^k\!, \bq^k\!, \by^k\!, \bt^k}{\text{min}}  w^k\! \Big(\!\sum_i x_{i,0}^k d_{i,0}\!+\!\!\sum_{i,j} x_{i,j}^k d_{i,j}\! \Big) \!\!\!\nonumber\\
% &\!\!\! + p_0 y_0^k +\! \sum_j p_j^* y_j^k +\theta^k\! \sum_i \psi_i^k q_i^k + \!\sum_j \left(\phi_j^k+s^kp_j^{s*}\right) t_j^k\!\!\!\!\! \\
% &\text{subject to:}~~~~ (\bx^k, \bq^k, \by^k, \bt^k) \in \mathcal{S}^k(\bz^*).
% \end{alignat}
% \end{subequations}
\begin{align*}
\!\textbf{SP1}^k \!\!: \varphi^k(\bp^*\!,\!\bp^{s*}\!\!,\!\bz^*)\!=\!\!\!\!\min_{\substack{(\bx^k\!, \bq^k\!, \by^k\!, \bt^k)\\ \in \mathcal{S}^k(\bp^*\!,\bp^{s*}\!,\bz^*)}} \!\!\!\! \mathcal{C}^k(\bx^k,\bq^k,\by^k,\bt^k|\bp^*\!,\!\bp^{s*}).
\end{align*}
The solution $(\bx^{k*},\bq^{k*},\by^{k*},\bt^{k*})$ to $\textbf{SP1}^k$ provides the optimal workload allocation, unmet demand, resource procurement, and service placement of service $k$ under the resource price and EN activation decision $(\bp^*\!,\bp^{s*}\!\!,\bz^*)$ of the platform. 
Note that each follower problem $\textbf{SP1}^k$ is an LP that may have multiple optimal solutions. Furthermore, the service placement and resource procurement decisions of the services directly affect the revenue of the platform. To derive a valid LB for the original bi-level problem \bmip, we consider the following subproblem, \textbf{SP2}, to select the optimal solutions to problems $\textbf{SP1}^k$, which support the objective of the leader. By using valid inequalities, we define \textbf{SP2} as follows:
% The second subproblems, \textbf{SP2}, are generally necessary because the lower-level problems may have multiple optimal solutions for \textbf{SP1} and doing so derives the one in support of the upper-level problem (i.e., valid inequalities). Thus we denote the two subproblems as follows:
\begin{subequations}
\label{SP2}
\begin{alignat}{3}
&\!\textbf{SP2}\!:~ \Theta_o(\bp^*, \bp^{s*}, \bz^*) = \underset{\bx, \bq, \by, \bt}{\max} ~\sum_{j,k}p^*_j y_j^k +\nonumber \\  &\!\quad\sum_{j,k}\left(\phi_j^k+s^kp^{s*}_j\right) t_j^k- \sum_j \Bigg(f_j z^*_j + c_j  \frac{\sum_k y_j^k}{C_j} \Bigg)\\
\label{sp2opt1}
&\!\text{subject to:} ~~~ (\bx^k, \bq^k, \by^k, \bt^k) \in \mathcal{S}^k(\bp^*,\bp^{s*},\bz^*),~\forall k,  \\
\label{sp2opt2}
& \mathcal{C}^k(\bx^k,\bq^k,\by^k,\bt^k|\bp^*\!,\!\bp^{s*})   \leq \varphi^k(\bp^*\!, \bp^{s*}\!, \bz^*),~\forall k,
% \label{sp2opt}
% & w^k \Big(\! \sum_i x_{i,0}^k d_{i,0} + \! \sum_{i,j} x_{i,j}^k d_{i,j} \Big) +\theta^k\! \sum_i \psi_i^k q_i^k +\!\sum_j p_j^* y_j^k \! \nonumber\\
% &+ p_0 y_0^k + \!\sum_j \left(\phi_j^k+s^kp_j^{s*}\right) t_j^k   \leq \varphi^k(\bp^*\!, \!\bp^{s*}\!\!,\! \bz^*),~\forall k,
\end{alignat}
\end{subequations}
where 
 $\varphi^k(\bp^*, \bp^{s*}, \bz^*)$ is the optimal value of $\textbf{SP1}^k$. 
%As integer variables are involved at both levels, they can be easily computed using existing MIP solvers such as Gurobi and Mosek. 
Note that \textbf{SP2} is a MILP that can be solved by existing MIP solvers.
%Also we can express the UB as follows:
Furthermore, it is easy to see that  ($\bp^*, \bp^{s*},  \bz^*, \bx^{*},\bq^{*},\by^{*},\bt^{*}$) is a feasible solution set of the original bi-level maximization problem \bmip. Hence, the optimal value $\Theta_o(\bp^*, \bp^{s*}, \bz^*)$ of \textbf{SP2} provides a LB to \bmip. To ensure convergence of the proposed algorithm, 
the LB to \bmip~  should be weakly improved after each iteration. Thus, we define:
\beqn
\label{LBdef}
LB = \text {max} \Big\{LB, ~ \Theta_o(\bp^*,  \bp^{s*}, \bz^*) \Big\}.
\eeqn
The solution to the subproblem \textbf{SP2} provides a new constraint in the form of \eqref{eq:SubProbStrongDuality3} to be added to the MP after each iteration.

\subsection{Proposed Algorithm}
\begin{figure}[t!]
\vspace{-0.4cm}
    \centering
    \includegraphics[width=0.45\textwidth]{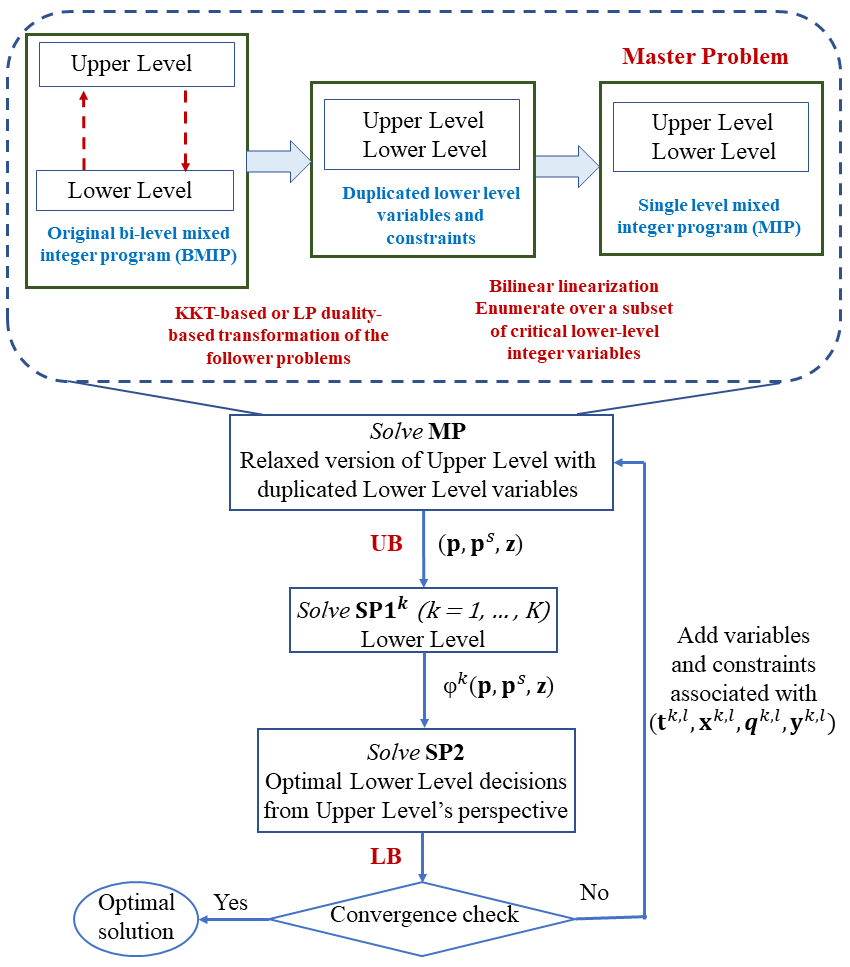}
    \caption{Illustration of the proposed solution algorithm.}
    \label{fig:algdiagram}
    \vspace{-0.5cm}
\end{figure}

Building upon the master and subproblems outlined earlier, we are ready to present \textbf{Algorithm 1}, an iterative decomposition approach for solving \bmip, or equivalently, \bmipl. %Based on the illustration of the master and subproblems above, we can solve \bmip~ to 
% Through the implementation of \textbf{Algorithm 1}, we can effectively optimize both the EN activation and edge resources pricing decisions for the platform, as well as the service placement, workload allocation, and resource procurement decisions for the services. 
Fig.~\ref{fig:algdiagram} further summarizes the proposed solution technique. 

\begin{algorithm}[h!]
Set $LB = -\infty$, $UB = +\infty$, and $L = 0$. \\
% Check, \linebreak
%  (I) If  $k \geq n$, then the  $m_{2}$. \\
Solve the \textbf{MP} in (\ref{MPfinal}) given $(\bt^{1*},\ldots,\bt^{L*})$ and derive the optimal solution $\big( \bp^*,\bp^{s*},\bz^*,\bx^{\prime *},\bq^{\prime *},\by^{\prime *},\bt^{\prime*}\!$, $\Gamma^{1*},\ldots,\Gamma^{L*}\!, \mu^{1*}\!,\ldots,\mu^{L*}\!, \nu^{1*},\ldots,\nu^{L*},\sigma^{1*},\ldots,\sigma^{L*}$, $ \xi^{1*},\ldots,\xi^{L*}, \eta^{1*},\ldots,\eta^{L*},\tau^{1*},\ldots,\tau^{L*} \big)$ at iteration $L$. Then, update $UB = \Theta^*$.\\
If $\frac{UB - LB}{UB}$ $\leq$ $\epsilon$, return the associated $LB$ and terminate. Otherwise, go to step 4.\\
Solve $\textbf{SP1}^k,~\forall k,$ for a given optimal solution of the \textbf{MP} in step 2 and obtain $\varphi^k(\bp^*, \bp^{s*},  \bz^*)$.\\
Solve \textbf{SP2}.  Report the optimal $\bt^{L+1,*}=\{\bt^{k*}, \forall k\},$ and $\Theta_o(\bp^*, \bp^{s*},  \bz^*)$. Update $LB$ using (\ref{LBdef}).\\
Add cuts \eqref{eq:SubProbStrongDualityConst1}--\eqref{eq:SubProbStrongDualityConste}, and  \eqref{eq:SubProbStrongDualityObj} %(\ref{cut1})-(\ref{cut5})
 with $l = L+1$ to the \textbf{MP}. Set $L = L+1$ and go to Step 2.
\caption{Iterative Algorithm for Solving \bmip}
\label{algo:CCG}
\end{algorithm}

% \vspace{0.2cm}
% \hline
% \vspace{0.1cm}
% \textbf{Algorithm 1}: \textit{CCG algorithm for BiMIP}\\
% \vspace{-0.25cm}
% \hline
% \vspace{0.1cm}
% \begin{enumerate}
%   \item Set $LB = -\infty$, $UB = +\infty$, and $L = 0$. 
%   \item Solve the master problem in (\ref{masterccg}) and derive optimal solutions ($p^*$, $z^*$, $x^{\prime *}$, $y^{\prime *}$, $t^{\prime*}$, $x^{1*}$,\ldots,$x^{L*}$, $y^{1*}$,\ldots,$y^{L*}$, $t^{1*}$,\ldots,$t^{L*}$, $\Gamma^{1*}$,\ldots,$\Gamma^{L*}$, $\mu^{1*}$,\ldots,$\mu^{L*}$, $\sigma^{1*}$,\ldots,$\sigma^{L*}$, $\xi^{1*}$,\ldots,$\xi^{L*}$, $\eta^{1*}$,\ldots,$\eta^{L*}$, $\tau^{1*}$,\ldots,$\tau^{L*}$) and update $UB = \underline{\Theta}^*$. 
%   \item If $\frac{UB - LB}{UB}$ $\leq$ $\epsilon$, return the associated $UB$ and terminate. Otherwise, go to step 4.
%   \item Solve \textbf{SP1} for given optimal solutions of step 2 and obtain $\varphi^k(\bp^*, \bp^{s*},  \bz^*)$.
%   \item Solve \textbf{SP2}.  Report optimal ($x^{k*}$,$y^{k*}$,$t^{k*}$) and $\Theta_o(\bp^*, \bp^{s*},  \bz^*)$. Update $LB$ using (\ref{ubccg}).
%   \item Add cuts (\ref{cut1})-(\ref{cut5}) to the \textbf{MP}. Set L = L+1 and go to Step 2.
% \end{enumerate}
% \hline

To implement our proposed solution, certain information from SPs is required. However, this information is typically not private or sensitive and is commonly disclosed in practice. %Specifically:
\begin{remark}
    Our model addresses the long-term relationship between the platform and SPs. In practical settings, SPs typically declare their resource demand and QoS requirements (e.g., maximum average delay) to the platform, similar to the practices observed in cloud computing systems offered by providers such as Amazon and Google Cloud. Additionally, SPs specify their requirements for ENs, such as proximity and support for specific operating systems (e.g., Ubuntu), which are essential for their operations. These types of requirements are generally not considered highly sensitive. 
%Furthermore, both the platform and SPs can leverage historical data to forecast future demand during the scheduling period, thereby mitigating the need for complete real-time information. 
This approach is analogous to the utilization of Service Level Agreements (SLAs) in cloud computing, where declared requirements form the basis for resource allocation and quality assurance.
\end{remark}

\section{Computational Complexity}
\begin{table*}[t!]
\centering
\begin{tabular}{|l|c|c|c|}
\hline
& \textbf{MP}                  & $\textbf{SP1}^k (k\in\mathcal{K}$) & \textbf{SP2} \\ \hline
\!\!\! Number of constraints          & \!\!\!$6J+K\big[J+L+(L+1)\big(1+4J+I(J+2)+3J(V+H)\big)\big]$\!\!\! & \!\!\!$2+2I+3J+IJ$\!\!\!                       & \!\!\!$K(3+2I+3J+IJ)$\!\!\!               \\ \hline
\!\!\! Number of continuous variables\!\!\! & $K\big[1+I+J+2L(I+2J+1)+J(I+H+V)(L+1)\big]$           & \!\!\!$1+2I+J+IJ$\!\!\!                        & \!\!\!$K(1+2I+J+IJ)$\!\!\!                \\ \hline
\!\!\! Number of binary variables     & $J(1+V+H+K)$                                  & $J$                                & $KJ$                          \\ \hline
\end{tabular}
\caption{Problem size of the master problem (MP) and subproblems in \textbf{Algorithm 1} at iteration $L$.}
\label{tab:NoConstVars}
\end{table*}

Bi-level programs are inherently difficult to solve due to their generic non-convex and non-differentiable nature. Even in the simplest case of linear bi-level programming, where all involved functions are linear programs, the complexity is substantial and known to be strongly NP-hard \cite[Section 1.2.2]{Marcotte2005}. Moreover, verifying strict local optimality in linear bi-level programming is also NP-hard \cite[Theorem 2.2]{Marcotte2005}.

Our proposed approach involves transforming the formulated BMIP into an MILP. Note that MILP, a challenging class of optimization problems in itself, is NP-hard, and when treated as a decision problem, it is NP-complete \cite{Conforti2014}. In this study, our focus is not on proposing efficient approximation algorithms for solving \bmip. Instead, we aim to devise an iterative exact algorithm for computing an optimal solution by converting the original challenging bi-level mixed-integer nonlinear program \bmip~ into an MILP (\bmipl), which can be solved by off-the-shelf solvers. This exact algorithm can serve as a benchmark for future research to evaluate the performance of new approximation algorithms.

Most MILP solvers commonly integrate heuristic algorithms during a pre-solve phase to effectively reduce problem size before utilizing techniques such as branch-and-bound (BnB) and branch-and-cut (BnC) \cite{Conforti2014}. The computational time of these methods is primarily influenced by the problem size, particularly the number of integer variables \cite{Conforti2014}. The specifics of the problem, including the number of constraints and binary/continuous variables for each of the master and subproblems, are detailed in Table~\ref{tab:NoConstVars}.

The following result demonstrates the finite convergence of the proposed method in \textbf{Algorithm 1}.
\begin{proposition}
\textbf{Algorithm 1} converges within a finite number of iterations to the optimal value of \bmip~in $\mathcal{O}(Q)$ iterations, where $Q = |\bT^k| = 2^J$ represents the number of possible values for the placement decision $\bt^k$ in the lower-level problems, for all $k \in \mathcal{K}$.
\end{proposition}

\begin{proof}
By solving the \textbf{MP}, we obtain an optimal solution  ($\bp^*, \bp^{s*},  \bz^*$), which serves as input to the subproblems. We then solve subproblems and obtain the optimal solution ($\bx^{k*},\bq^{k*},\by^{k*},\bt^{k*}$). As ($\bp^*, \bp^{s*},  \bz^*, \bx^{k*},\bq^{k*},\by^{k*},\bt^{k*}$) is a feasible solution set of \bmip,  $\Theta_o(\bp^*, \bp^{s*}, \bz^*)$ provides a LB to the problem. 
By definition \eqref{LBdef}, the LB is non-decreasing. 

After solving the subproblems, we expand the \textbf{MP} by adding new cuts related to the new optimal solution $\bt^{k*}$, which also expands the partial enumeration over the sets $\bT^k$. Solving the newly updated \textbf{MP} with the added constraints provides a stronger  UB to \bmip. By iteratively solving the MP and subproblems, we anticipate the LBs and UBs to converge to the optimal value. 

%the  set of possible lower-level integer values  by adding the new optimal solution $\bt^{k*}$ to each set $\bT^k$. The \textbf{MP} is also expanded by adding new cuts. 
%Next, we update the $\bt^k$ set by including $\bt^{k*}$ and expand the master problem $\underline{\Theta}$.
% Solving the augmented master problem with the new set $\bt^k$, that involves the associated constraints up to iteration $k$, leads to new $\bp^*$ and $\bz^*$ optimal values, as well as a stronger UB (weakly-decreasing) for the new iteration. By iteratively solving the master problem and subproblems, we anticipate the lower and UBs to converge to the optimal value. 

Note that a new $\bt^{k*}$ is generated after each iteration. Thus, \textbf{to show the convergence within finite $\mathcal{O}(Q)$ iterations}, we can instead prove that \textit{any repeated $\bt^{k*}$ implies $LB = UB$}.
This can be shown by contradiction as follows. 
%Indeed, the convergence can  be shown by contradiction. 
%that any repeated $\bt^{k*}$ implies $LB = UB$.
Let $L_1$ denote the current iteration index, ($\bp^{*}, \bp^{s*}, \bz^*$, $\bx^{\prime k*}$, $\by^{\prime k*}$, $\bt^{\prime k*}$) is obtained in Step 2 with $\frac{UB - LB}{UB} > \epsilon$, and ($\bx^{k*}$, $\by^{k*}$, $\bt^{k*}$) is obtained from Step 4 of \textbf{Algorithm 1}. From Step 5, we have: 
\begin{align*}
&LB \geq \Theta_o(\bp^*, \bp^{s*},  \bz^*) \\
= & \sum_{j,k} \left[ p^*y_j^{k*} \!+\!\left(\phi_j^k+s^kp^{s*}_j\right) t_j^{k*} \right] \!-\! \sum_j \!\bigg(\!f_j z_j^* + c_j  \frac{\sum_k y_j^{k*}}{C_j} \bigg).
\end{align*} 
%and optimal ($\bx^{k*}$, $\by^{k*}$, $\bt^{k*}$) that is in favor of the upper-level decision maker.
% Now assume that $\bt^{k*}$ has appeared in some previous iteration $L_0 < L_1$. At the current iteration $L_1$, since $\frac{UB -LB}{UB} > \epsilon$, we augment the \textbf{MP} with a set of new variables and constraints associated with $\bt^{k*}$ ($=\bt^{k,L_1 + 1}$), for all $k$. However, since those variables and constraints are the same as those created and included in iteration $L_0$. Thus, the \textbf{MP} remains the same after iteration the current iteration $L_1$.
% Consequently, in iteration $L_1 + 1$, the \textbf{MP} yields the same optimal value as the value in iteration $L_1$. Hence, when the algorithm goes from iteration $L_1$ to iteration $L_1 + 1$, the UB does not change and the set of optimal solution obtained at iteration $L_1 + 1$ remains the same as ($\bp^*, \bp^{s*},  \bz^*, \bx^{k*},\bq^{k*},\by^{k*},\bt^{k*}$).
Suppose that $\bt^{k*}$ has appeared in some previous iteration, say $L_0 < L_1$. In the current iteration $L_1$, we have $\frac{UB - LB}{UB} > \epsilon$, and augment \textbf{MP} with a set of new variables and constraints associated with $\bt^{k*}$ ($=\bt^{k,L_1 + 1}$), for all $k$. However, since these variables and constraints are identical to those created and included in iteration $L_0$. Therefore, the \textbf{MP} remains unchanged after the current iteration $L_1$. As a consequence, in the next iteration $L_1 + 1$, the \textbf{MP} yields the same optimal value as in iteration $L_1$. This means that the UB remains the same, and the optimal solution obtained in iteration $L_1 + 1$ is identical to that obtained in iteration $L_1$, i.e., ($\bp^*, \bp^{s*}, \bz^*, \bx^{k*},\bq^{k*},\by^{k*},\bt^{k*}$).

We will show that at $L_1 + 1$, UB $\leq$ LB. According to the update for UB in \eqref{eq-UB}, we have: 
% \begin{align*}
%     &~~UB ~~=~~\mathcal{P}\left(\bp^*,\bp^{s*},\bz^*,\by^*,\bt^*\right)\\
%     \le&~~\underset{\bx^k, \bq^k, \by^k, \bt^k \in \mathcal{S}^k(\bz^*)}{\max} ~{P}\left(\bp^*,\bp^{s*},\bz^*,\by,\bt\right) \\
%     =&~~\Theta_o(\bp^*, \bp^{s*}, \bz^*) ~~\le~~ LB,
% \end{align*}
\begin{align*}
    UB \!=&\! \! \sum_{j,k} \!\left[ p^*y_j^{\prime k*} \!\!+\!\left(\phi_j^k\!+\!s^kp^{s*}_j\right)\! t_j^{\prime k*} \!\right] \! -\!\! \sum_j\!\! \left(\!f_j z_j^* \!+\! c_j  \frac{\sum_k y_j^{\prime k*}}{C_j}\!\right)\!\nonumber \\ 
    \le&~~\underset{\bx, \by, \bt}{\max} ~~\Bigg\{\sum_{j,k}\left[p^*_j y_j^k +\left(\phi_j^k+s^kp^{s*}_j\right) t_j^k\right] \nonumber \\ 
    &- \sum_j \left(f_j z^*_j + c_j  \frac{\sum_k y_j^k}{C_j} \right):\eqref{sp2opt1},\eqref{sp2opt2}\Bigg\}\\
    =&\Theta_o(\bp^*, \bp^{s*}, \bz^*) \le LB,
\end{align*}
where the last equality follows from solving \textbf{SP2} in the previous iteration $L_1$ and observing that the set of solutions ($\bp^*, \bp^{s*}, \bz^*$) in iteration $L_1+1$ is the same as that in iteration $L_1$. Therefore, the optimal value of \textbf{SP2}, $\Theta_o(\bp^*, \bp^{s*}, \bz^*)$, remains unchanged, which is less than or equal to the current LB. 
Since UB $\leq$ LB, Step 3  terminates the algorithm, and we obtain the optimal solution. %This typically occurs in a finite number of iterations. %, indicating the convergence of the algorithm.
%The worst-case number of iterations is 
% As $\bt^{k,L_1 + 1} = \bt^{k*}$ is optimal to %$\varphi^k(\bp^*, \bp^{s*},  \bz^*)$, $\forall k$
% $\textbf{SP1}^k$ and constraints from strong duality ensures that at iteration $L_1 + 1$:
% \begin{align*}
% \varphi^k(\bp^*,\bp^{s*},\bz^*) = \sum_j (\phi_j^k+s^kp^{s}_j) t_j^{k,L_1 + 1} - \mu_1^{k,L_1 + 1} \Big(\!B^k - \\
% \sum_j (\phi_j^k+s^kp^{s}_j) t_j^{k,L_1 + 1}\! \Big) \!\!-\!\! \sum_i \!R_i^kD^{k,m} \xi_i^{k,L_1 + 1} \!\!+\!\! \sum_i \!R_i^k \eta_i^{k,L_1 + 1}\\ - \sum_j C_j t_j^{k,L_1 + 1} \sigma_j^{k,L_1 + 1}-\sum_{i,j} a_{i,j}^k \tau_{i,j}^{k,L_1 + 1} R_i^k ,~\forall k.
% \end{align*}
\end{proof}

\begin{remark}
Similar to well-established methodologies like Benders decomposition, column generation, and Dantzig-Wolfe decomposition, the time complexity of our proposed decomposition algorithm is contingent on the intricacy involved in solving the associated master and subproblems. Decomposition methods are widely used because, despite the lack of a polynomial guarantee, they often lead to more efficient solutions in practice.
\end{remark}

\section{Numerical Results}
\label{sim}
\subsection{Simulation Setting}
Similar to the prior work \cite{tnduong21,duong20}, we employ the widely-used Barabasi-Albert model with an attachment rate of $2$ to generate a topology for an edge network that comprises $100$ nodes. The link delay between adjacent nodes is randomly generated within the range of $[2, 5]~ms$, and the network delay between any two nodes is determined as the shortest path between them \cite{tnduong21}. 
The delay between each area and the remote cloud is set to  $60~ms$.  The maximum delay thresholds $D^{k,m}$ for the services are randomly chosen from the range of $30~ms$ to $100~ms$. In the  \textbf{base case}, we consider a small-scale system comprising $K=4$ services, $I=12$ areas, and $J=8$ ENs, which are randomly selected from the set of $100$ nodes. We also explore larger system sizes during sensitivity analyses. %to examine the sensitivity of the proposed algorithms.

%During the scheduling horizon, t
For each area, the workload (i.e., resource demand) is randomly drawn from the range $[20, 35]$ \textit{vCPUs} \cite{tnduong21}. 
% During the scheduling horizon, the resource demand (i.e., workload) in each area is randomly drawn in the range of $[20, 35]$ vCPUs.
Each EN is selected at random from the set of Amazon EC2 M5 instances.
%Each EN is randomly chosen from the set of Amazon EC2 M5 instances. 
Using the hourly price of an \textit{m5d.xlarge} Amazon EC2 instance \footnote{https://aws.amazon.com/ec2/pricing/} as reference, the unit price of computing resources at the cloud is set to  $0.02~ \$ $\textit{/vCPU}, while the set of possible unit prices of computing resources at the ENs is $[0.01, 0.02, 0.03, 0.04, 0.05]~ \$ $/\textit{vCPU}. The set of possible unit storage prices is $[0.005, 0.01, 0.015] \$ $/\textit{TB}. 
To account for the operational costs of each EN, we set the fixed and variable operational costs 
($f_j$ and $c_j$) %of each EN $j$
within the ranges of  $[\$0.1, \$3.6]$ and $[\$0.04, \$1.44]$, respectively. 
%, depending on the size of the EN. 
The delay penalty parameters, $w^k$, are randomly generated within the range of $[10^{-5}, 10^{-3}] \$ $/$(vCPU.ms)$. The installation cost of each service at each EN is set to $\$0.2  $ (i.e., $\phi_j^k = 0.2$, $\forall j, k$). The size of each service is randomly generated between $20$GB and $1$TB, and the budget of each service ranges from $\$20 $ to $\$50 $.

% In our experiments, unless stated otherwise, we use the default setting. We conduct our computational study using Matlab 2022a software and solve the optimization problem using the Mosek solver on an Apple M1 CPU with 8 GB RAM.

In all simulations, we have employed the default settings unless explicitly mentioned otherwise.  Our computational experiments were conducted using \textit{Matlab 2022a} and %, with the optimization problem solved via 
Mosek
% \footnote{https://www.mosek.com/} 
on an Apple M1 CPU with 8 GB of RAM.% Apple M1 CPU equipped with 8 GB of RAM.

\subsection{Performance Evaluation}
\subsubsection{Base Case Solution}
% First, we solve for the base case scenario.
%where the cloud price $p_0$ is set to $0.02$. 
% The convergence behavior of the proposed algorithm for this case is depicted in Fig.\ref{fig:BaseConv}. 
The proposed model incorporates both storage and computing resource pricing design. First, we %evaluate and 
compare the performance of the proposed scheme with two other schemes: (i) \textit{No storage cost}: the services do not have to pay for the storage costs. The storage prices are implicitly set to zero; and (ii) \textit{No placement cost} \cite{tnduong21}: the services do not have to pay for the service placement cost, including both service installation and storage costs. %solve the proposed scheme with dynamic pricing for both computing resources and storage within the base case configuration. To gain a more comprehensive understanding of workload allocation, we conduct a comparative analysis between this solution and a scenario where the platform disregards the costs associated with storage and placement. 
The optimal solutions on EN activation decisions as well as unit prices of computing and storage resources are presented in Table~\ref{tab:optSol}. For instance, in the proposed scheme, ENs $1, 2, 3$ and $5$ are active with the corresponding unit computing resource prices being $0.03$, $0.02$, $0.04$, and $0.04$, respectively. The unit storage price is the same at $0.005$ across all ENs. Notably, the active ENs differ among the schemes.  Furthermore, the computing resource prices tend to decrease when considering the service placement cost.
This is due to services having a reduced budget for purchasing edge resources, prompting the platform to lower resource prices to encourage resource consumption.
%It is because the services are left with less budget for purchasing edge resources. Therefore, the platform may reduce the resource prices to encourage more resource consumption. 

\begin{table}[t!]
\centering
\begin{tabular}{|c|ccc|cc|cc|}
\hline
 & \multicolumn{3}{c|}{\!\!\!\!\!Proposed Scheme\!\!\!\!\!}  & \multicolumn{2}{c|}{\!\!\!\!\!No Storage Cost\!\!\!\!\!}& \multicolumn{2}{c|}{\!\!\!\!No Placement Cost\!\!\!\!}  \\ \hline
\!\!\!\!\!\! Active ENs $j$ \!\!\!\!\!\!      & \multicolumn{1}{c|}{$1$}    & \multicolumn{1}{c|}{$2$}    &\multicolumn{1}{c|}{$\!\!\!\!\{3,5\}\!\!\!\!$}    & \multicolumn{1}{c|}{$2$}    & $\!\!\!\!\{3,4,5,6,7,8\}\!\!\!\!$ & \multicolumn{1}{c|}{$\!\!\!\!\{1,5\}\!\!\!\!$}    & $\!\!\!\!\{3,7\}\!\!\!\!$ \\ \hline
\!\!\!\!\!$p_j$ \!\!\!\!\! & \multicolumn{1}{c|}{\!\!\!$0.03$\!\!\!} & \multicolumn{1}{c|}{\!\!\!$0.02$\!\!\!} & \multicolumn{1}{c|}{\!\!\!$0.04$\!\!\!} & \multicolumn{1}{c|}{\!\!\!$0.02$\!\!\!} & \!\!\!$0.03$\!\!\! & \multicolumn{1}{c|}{\!\!\!$0.04$\!\!\!} & \!\!\!$0.05$\!\!\!            \\ \hline
$p^s_j$  & \multicolumn{3}{c|}{$0.005$} & \multicolumn{1}{c|}{\!\!\!NA\!\!\!}     & NA   & \multicolumn{1}{c|}{NA}     & NA                \\ \hline
\end{tabular}
\caption{Optimal decisions.}
\vspace{-0.5cm}
\label{tab:optSol}
\end{table}

\begin{figure}[t!]
    \centering
	\vspace{-0.2cm}
	\subfigure[With storage pricing]{
 \includegraphics[width=0.228\textwidth,height=0.11\textheight]{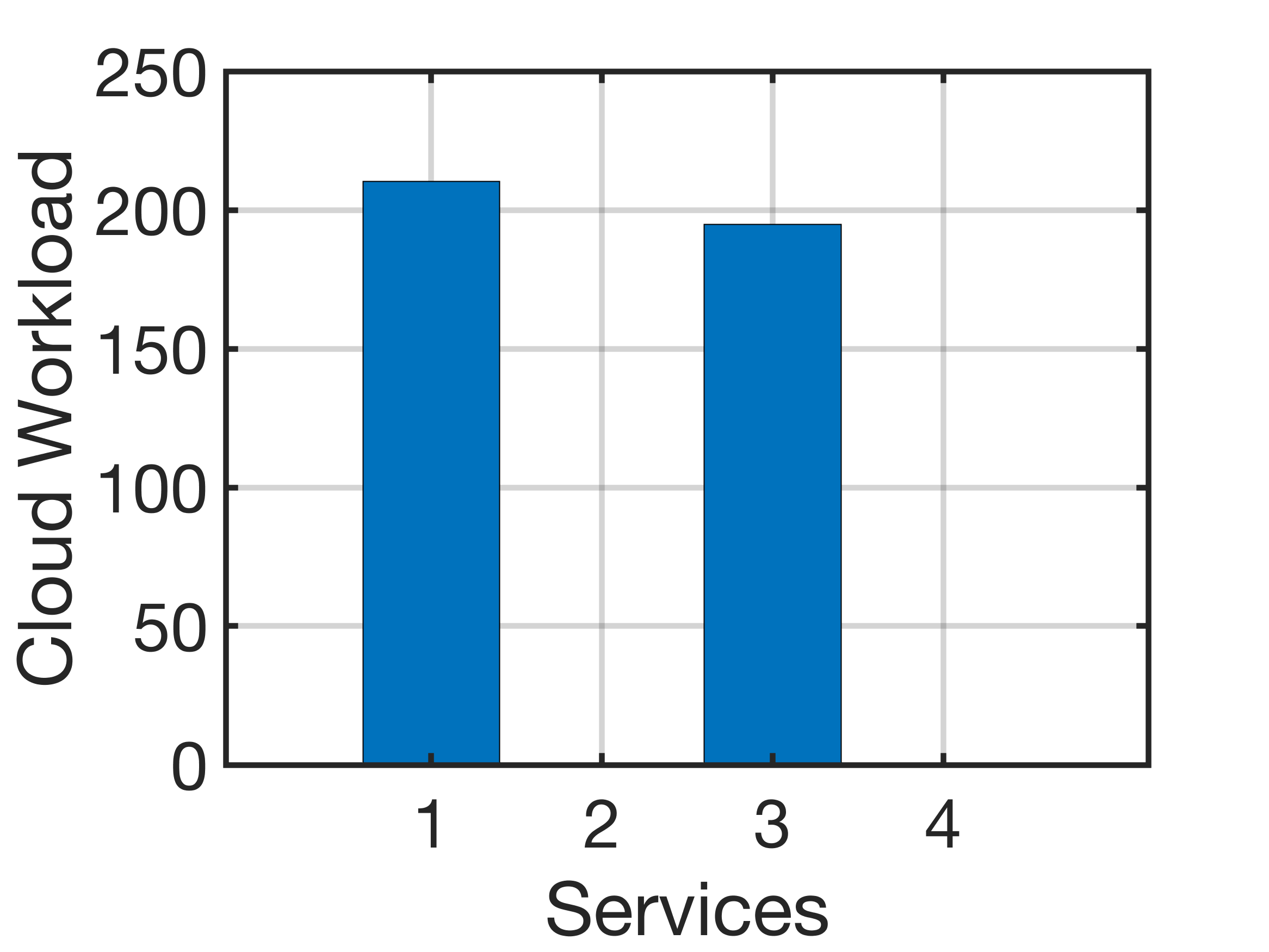}
	\label{fig:BaseCloudWLS}
	}   
	\hspace*{-1em}
	\subfigure[No storage cost]{
	 \includegraphics[width=0.232\textwidth,height=0.11\textheight]{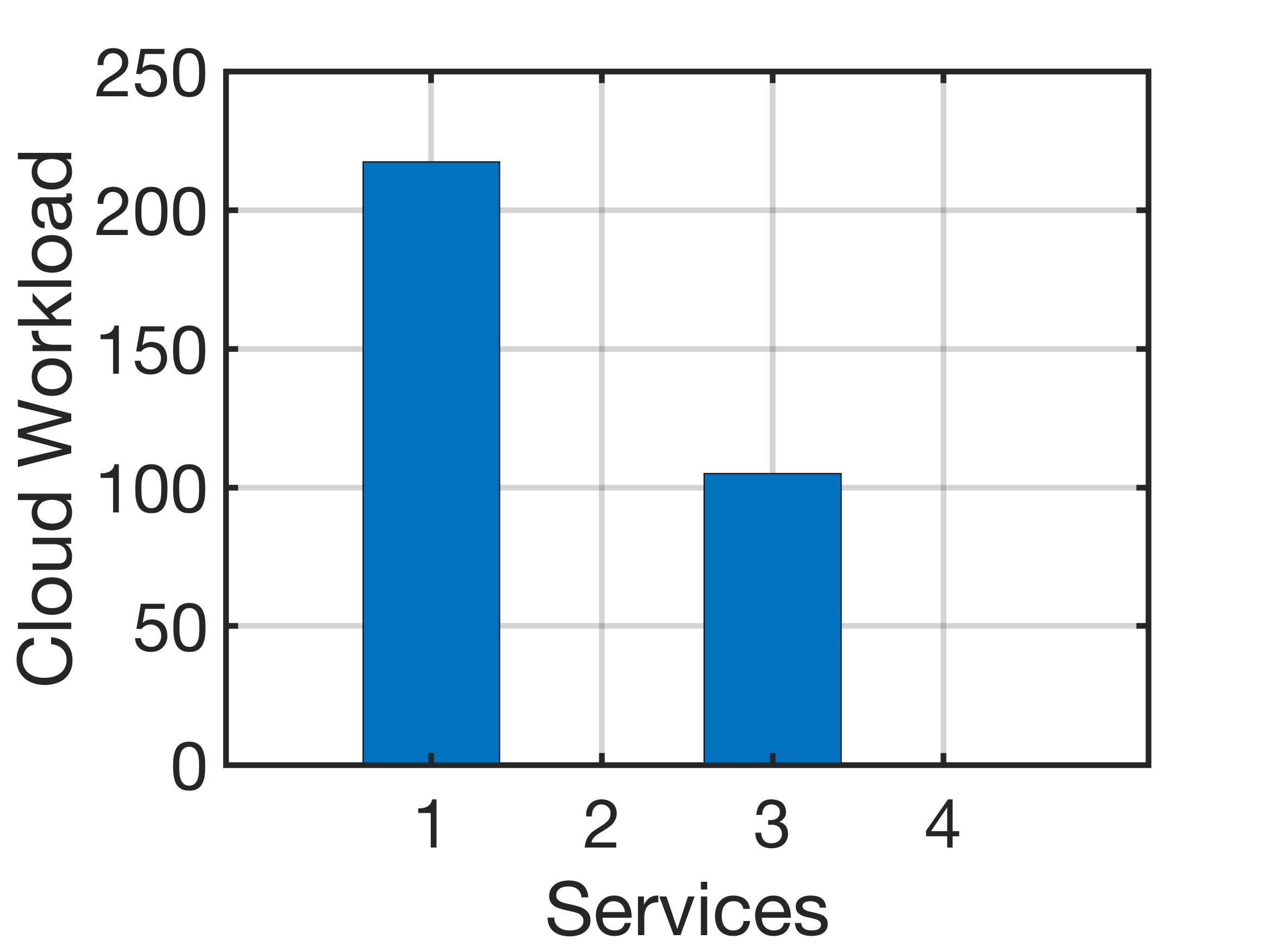}
	\label{fig:BaseCloudWL}
	}
	\subfigure[With storage pricing]{
	 \includegraphics[width=0.232\textwidth,height=0.11\textheight]{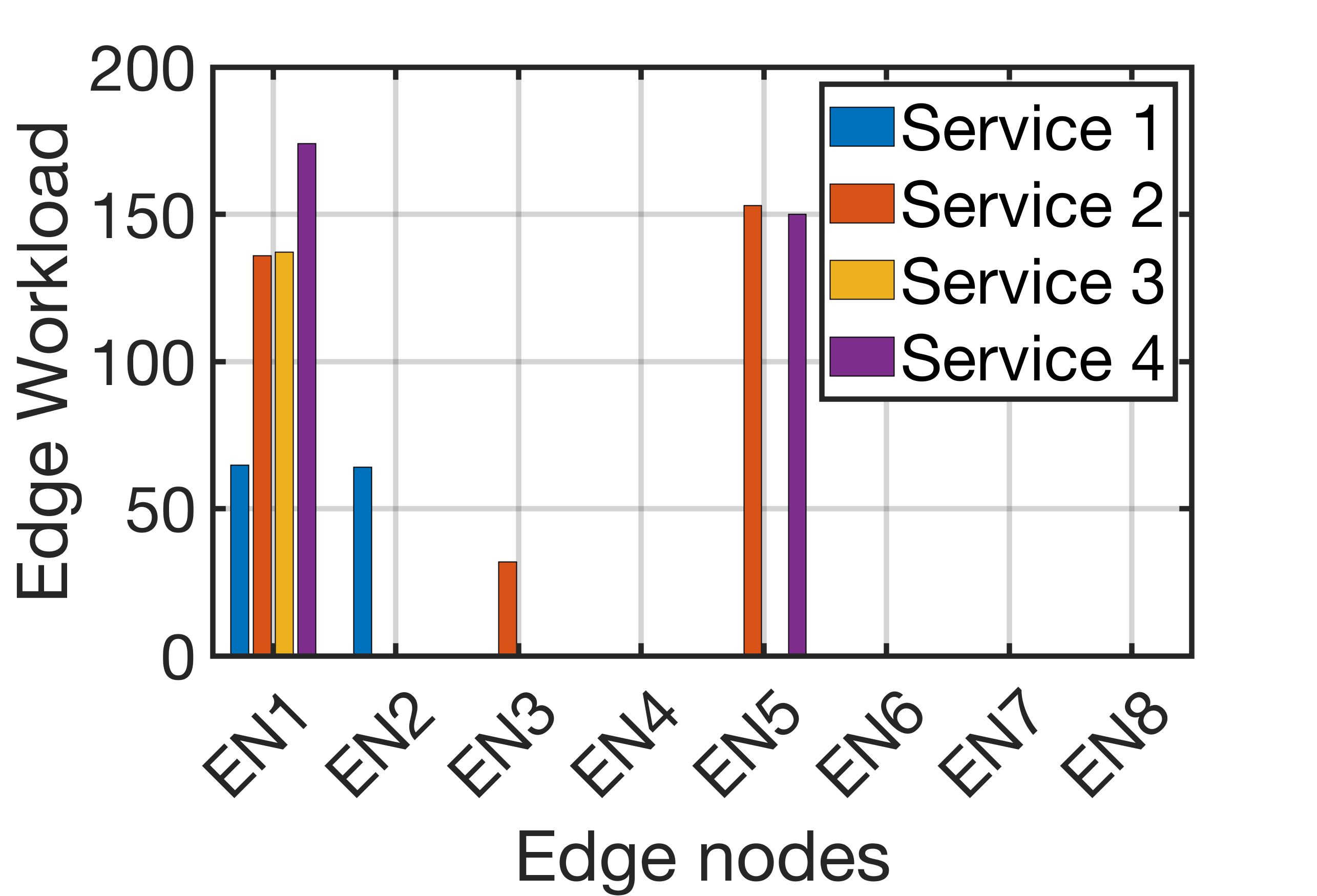}
	\label{fig:BaseEdgeWLS}
	} 
	\hspace*{-1.7em}
	\subfigure[No storage cost]{
	 \includegraphics[width=0.232\textwidth,height=0.11\textheight]{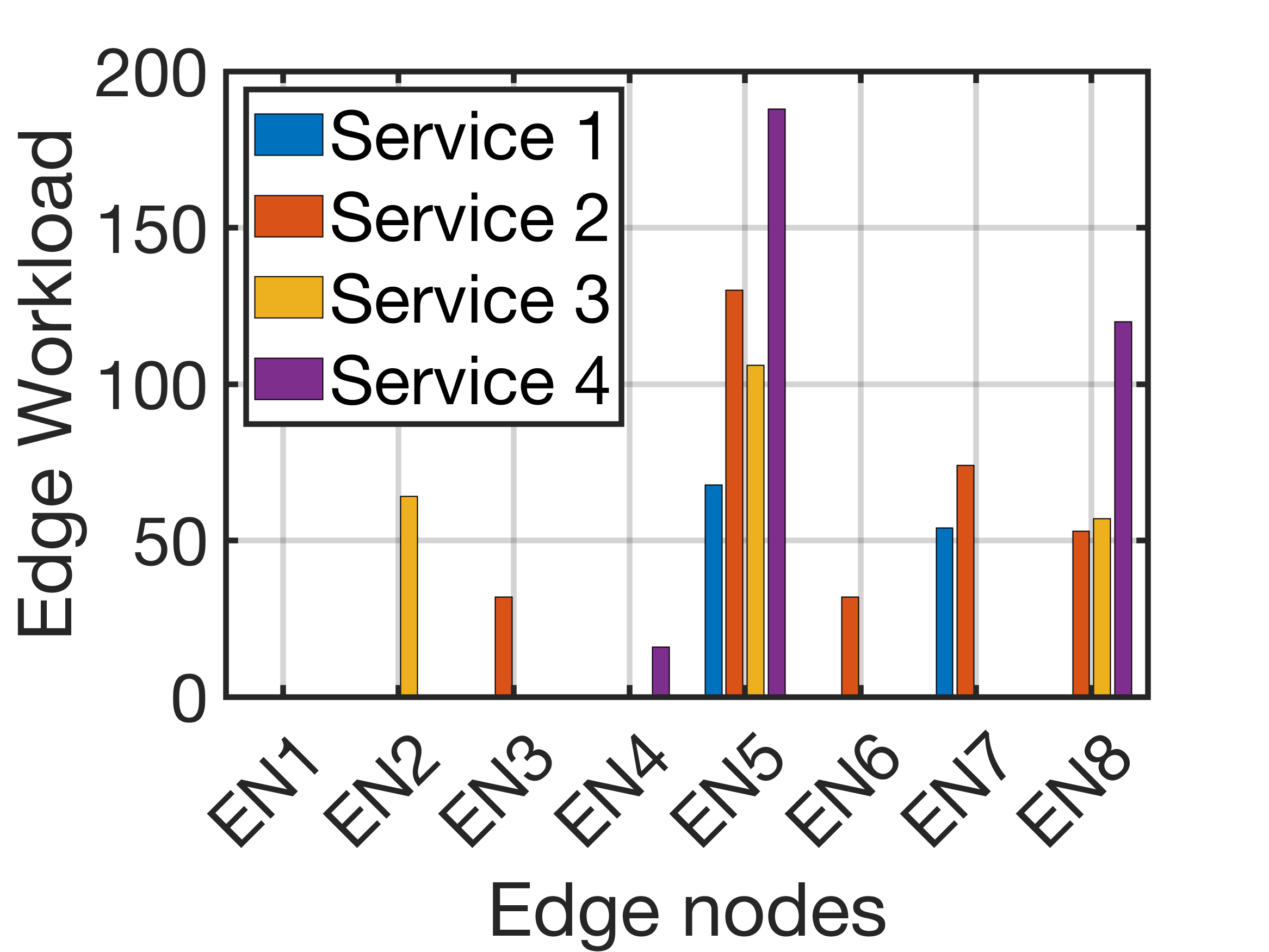}
	\label{fig:BaseEdgeWL}
	}
	\vspace{-0.2cm}
% 	\caption{Total workload allocated at the clouds and ENs for the proposed model with dynamic storage pricing (left) versus no storage cost (right).}
	\caption{Workload allocation.}
	\label{fig:BaseCase}
	\vspace{-0.7cm}
\end{figure}

Resource procurement of each service from both the cloud and ENs is visually depicted in Fig.\ref{fig:BaseCase} for the proposed model with storage pricing versus no storage costs. The figure highlights that when storage costs are not considered, services tend to acquire more computing resources from multiple ENs to minimize service delay.  
 In practical scenarios where storage costs are significant, as modeled in this paper, services need to carefully consider the trade-off between storage costs and resource acquisition. The cost of storage may outweigh the benefits of reducing delay. 
 % However, in practice, storage costs are usually a significant factor, as modeled in this paper. Therefore, services must carefully consider storage costs while acquiring resources from multiple ENs since the cost of storage may outweigh the benefits of minimizing delay.
 Consequently, services tend to obtain computing resources from fewer ENs to avoid paying for storage costs (i.e., less placement costs) for a small amount of computing resources from that EN, and instead, they offload the remaining extra tasks to the cloud. 
 This behavior is particularly observed in services 1 and 3, which have lower delay penalties and are less sensitive to delay compared to services 2 and 4, which prioritize ENs and are more sensitive to delay.
 %Specifically, this applies to services $1$ and $3$, which have a lower delay penalty and are less sensitive to delay compared to services $2$ and $4$, which prioritize ENs and are more sensitive to delay.
 It is worth noting that in the base case, although a smaller amount of resources are procured from ENs, the resulting profit is comparable to scenarios where storage cost is not considered.
 This is because the profit generated from charging storage costs and higher computing costs from services compensates for the reduced amount of rented computing resources.
 %This is because the profit generated from charging storage costs and higher computing costs from services compensates for the reduced amount of computing resources rented out. 

\begin{figure}[h!]
    \centering
    \subfigure[No placement cost]{
	\includegraphics[width=0.45\textwidth,height=0.12\textheight]{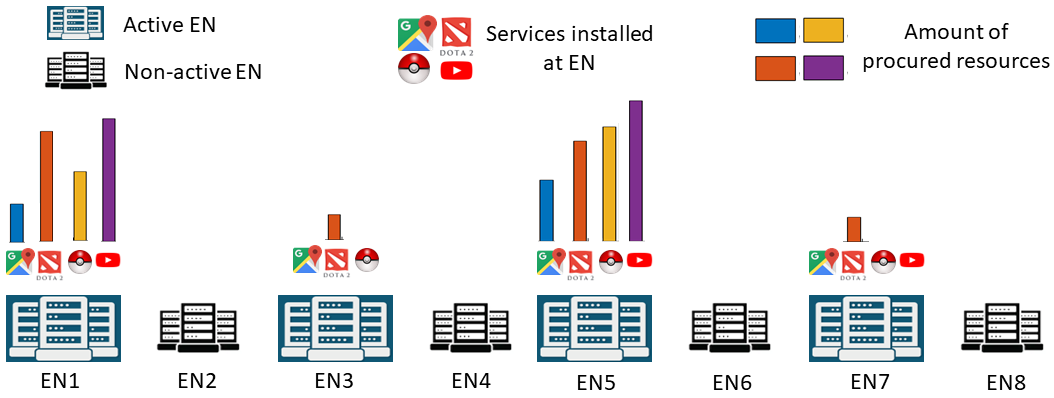}
	\label{fig:noplacement}}
 \vspace{-0.2cm}
    \subfigure[Proposed scheme]{
	\includegraphics[width=0.45\textwidth,height=0.09\textheight]{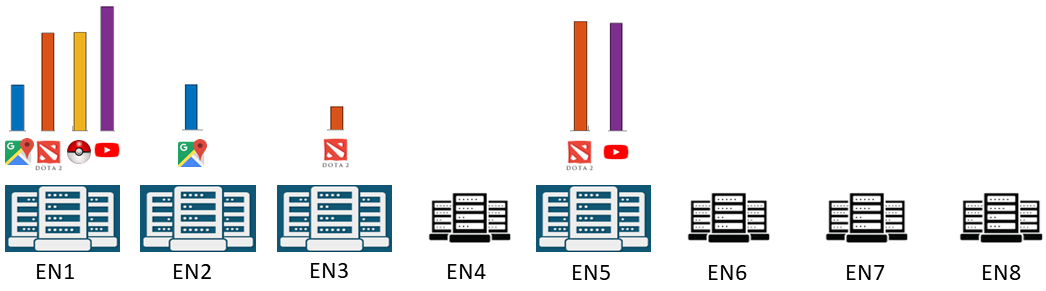}
	\label{fig:proposedscheme}}
	\caption{Services placement and procured computing resources.}
	\label{fig:noplacementandproposed}
 \vspace{-0.2cm}
\end{figure}

\begin{figure}[t!]
    \centering
    \vspace{-0.2cm}
	\includegraphics[width=0.3\textwidth]{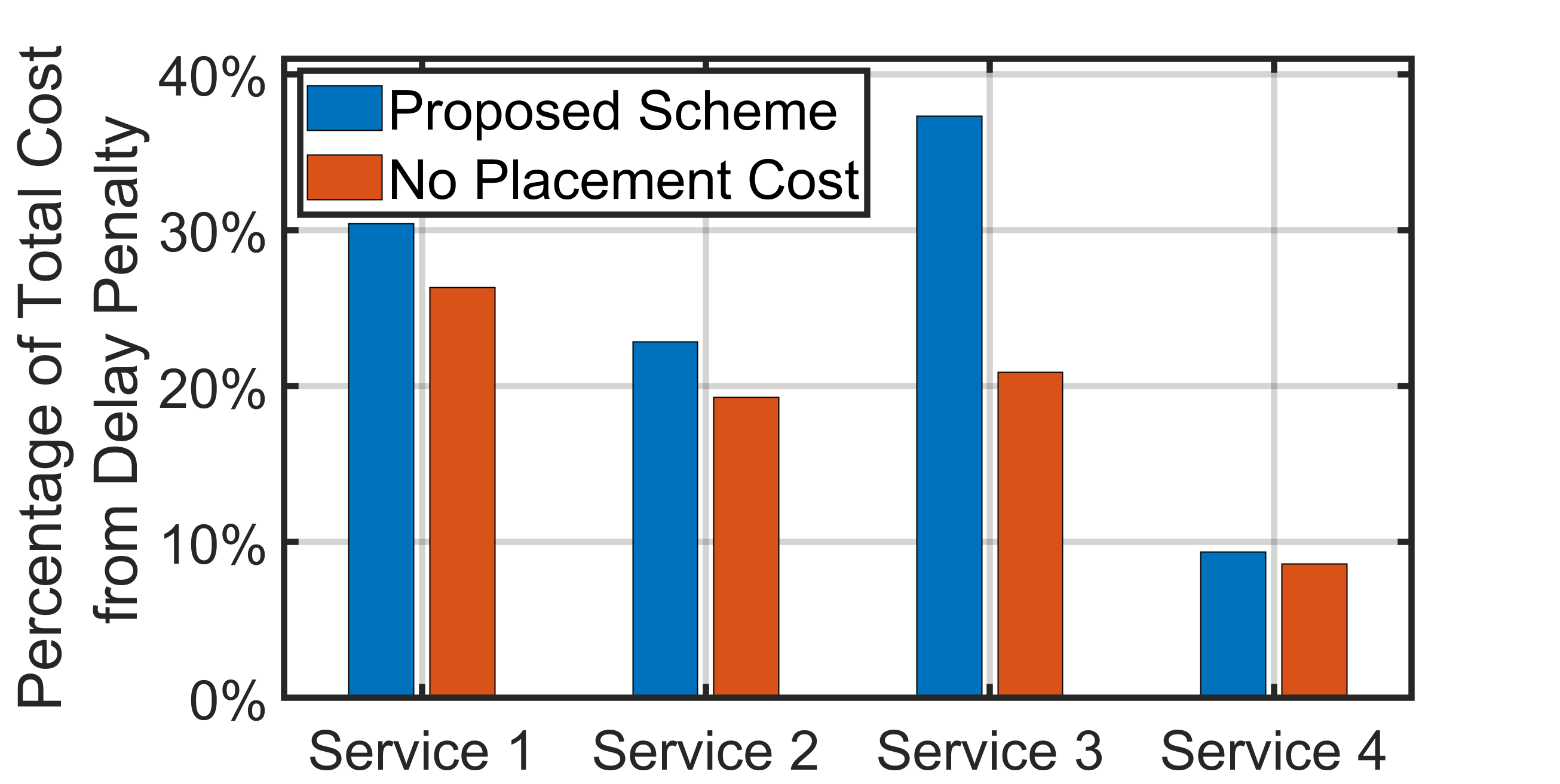}
	%\label{fig:CompareDelayPenalty}
	%\caption{The percentage of total cost resulting from  delay penalty.}
	\caption{Delay penalty comparison.}
	\label{fig:CompareDelayPenalty}
 \vspace{-0.2cm}
\end{figure}

Figure \ref{fig:noplacementandproposed} provides a further comparison of the placement and resource procurement decisions in the proposed model, which incorporates placement costs, and the model without placement costs as discussed in \cite{tnduong21}. %The figure clearly illustrates the distinction between these two approaches. 
In the absence of placement costs (Fig. \ref{fig:noplacement}), the services are installed at all active ENs, regardless of whether they acquire computing resources or serve any demand at those ENs. The placement decision is subject only to the storage capacity of each EN. %Indeed, services are installed 
%even when they do not procure any computing resources and do not serve any demand at those ENs.
This is as expected since there are no costs %or negative consequences 
associated with placing services at ENs. However, this indiscriminate installation of services at ENs without utilization poses an economic burden on the platform.  %as it ties up storage resources that remain idle. 
In contrast, our proposed scheme accounts for the placement cost, thus preventing over-installation (Fig. \ref{fig:proposedscheme}) and improving resource utilization. %, that is, the allocation of storage resources without actually serving any demand. 
% Our observation emphasizes the significance of considering both the service installation and storage costs. 
%and placement cost 
%as critical factors in the design of resource allocation schemes for EC systems. 
%By incorporating these costs into the decision-making process, more practical and cost-effective resource allocation strategies can be developed, resulting in better system performance and economic benefits.
Additionally, as expected, Fig. \ref{fig:CompareDelayPenalty} demonstrates that services experience higher delays when they have to bear service placement costs 
%shows that the services suffer higher delays when they have  to pay for service placement costs 
since they have less budget and flexibility for edge resource procurement.
The impact of incorporating service placement costs is more pronounced on less-delay-sensitive services 
%Incorporating service placement costs has more impacts on less-delay sensitive services 
(e.g., services 1 and 3), as they offload a larger workload to the cloud.
%since they have to offload more workload to the cloud.

\subsubsection{Computational Time}
Table~\ref{table:execution_times} provides a comparison of average and worst-case computational times between the proposed algorithm (\textbf{Algorithm 1}) and a Brute Force method, conducted across different network settings with each configuration involving $100$ simulations. It is noteworthy that the Brute Force method employs our proposed reformulation techniques to convert the formulated bi-level problem (\bmip) into an MILP problem (\bmipl). This reformulation itself is nontrivial. However, the Brute Force method opts out of leveraging the proposed decomposition and iteration framework present in \textbf{Algorithm 1}. Instead, it tackles the \bmipl~ problem by considering all $Q=2^J$ possible combinations of $\bt^k$ for the lower-level problems, in the form of constraint \eqref{eq:SubProbStrongDuality3}, for each $k$. The resulting MILP becomes considerably large and computationally expensive, particularly with the increase in $J$, as evident in Table~\ref{table:execution_times}. The Brute-Force method proves ineffective in solving the problem within the specified time limit of $18,000$ seconds, beginning at $J=10$. In contrast, \textbf{Algorithm 1} only engages constraint \eqref{eq:SubProbStrongDuality3} for crucial values of $\bt^k$ identified through solving the subproblems. While we lack theoretical results, our numerical simulations consistently demonstrate that \textbf{Algorithm 1} typically converges within a few iterations. Subsequent sections present additional results to further illustrate this convergence pattern (refer to Fig.~\ref{fig:ConvPlots}).

\begin{table}[t!]
\centering
\begin{tabular}{|c|cc|cc|}
\hline
\multirow{2}{*}{} & \multicolumn{2}{c|}{\textbf{Algorithm 1}}                                   & \multicolumn{2}{c|}{Brute Force}          \\ \cline{2-5} 
& \multicolumn{1}{c|}{\!\!\! Average (s)\!\!\!}   & \multicolumn{1}{c|}{\!\!\!Worst-case (s)\!\!\!} & \multicolumn{1}{c|}{\!\!\!Average (s)\!\!\!}    & \!\!\!Worst-case (s)\!\!\!  \\ \hline
$I = 6, J = 4$    & \multicolumn{1}{c|}{$4.145$}   & $14.151$                       & \multicolumn{1}{c|}{$3.933$}    & $7.172$    \\ \hline
$I = 9, J = 6$    & \multicolumn{1}{c|}{$27.328$}  & $116.670$                        & \multicolumn{1}{c|}{$111.949$}  & $211.110$  \\ \hline
\!\!\!$I = 12, J = 8$ \!\!\!  & \multicolumn{1}{c|}{$121.105$} & $422.4164$                      & \multicolumn{1}{c|}{$1451.033$} & $1654.544$ \\ \hline
% \!\!\!$I = 15, J = 10$ \!\!\!  & \multicolumn{1}{c|}{$1153.125$} & $3777.055$                      & \multicolumn{1}{c|}{\!\!\!\!$>\!4090.528$\!\!\!\!} & $>\!18000$ \\ \hline
\!\!\!$I = 15, J = 10$ \!\!\!  & \multicolumn{1}{c|}{$1153.125$} & $3777.055$                      & \multicolumn{1}{c|}{\!\!\!\!NA\!\!\!} & NA \\ \hline
\end{tabular}
\caption{Comparison of average and worst-case computational times between \textbf{Algorithm 1} and a Brute Force method over $100$ simulations.} 
\label{table:execution_times}
\vspace{-0.5cm}
\end{table}

\subsubsection{Sensitivity Analysis}
%For sensitivity analyses, we will examine the impact of various design parameters on the system performance. 
To perform sensitivity analyses, we will assess the impact of various design parameters on the system performance. Specifically, we compare the outcomes under three different pricing strategies: (i)  Dynamic pricing (\textit{Dyn}), where the platform solves the proposed model using \textbf{Algorithm 1} to determine the resource prices; (ii) Flat pricing  (\textit{Flat}), where the platform solves the proposed model under the condition that the resource prices at all the ENs are the same;  and (iii) Flat average pricing (\textit{Avg}), where the computing resource prices at all ENs are set to the average value within the price range (i.e., $p_j = 0.03, \forall j$). %In the \textit{Dyn} approach, the resource prices at the ENs are allowed to vary in order to balance the supply and demand. In the \textit{Flat} price scheme, the leader offers the same computing resource price at all ENs while in the \textit{Avg} approach, the price is set to  the average   price. %range. %We solve \textbf{SP1} and \textbf{SP2} for each price decision.
To quantify the effects of different system parameters, we introduce scaling factors $\delta$, $\gamma_0$, $\Lambda$, and $\beta_0$ for the resource demand $R$, EN capacities $C$, delay threshold $D^{k,m}$, and delay penalty $w$, respectively. For instance, when $\delta = 0.5$, the resource demand used in the experiment is half of that used in the \textit{base case}.

\begin{figure}[h!]
\vspace{-0.4cm}
    \centering
    \hspace*{-0.5em}
    \subfigure[Varying $p_0$]{
 \includegraphics[width=0.242\textwidth,height=0.11\textheight]{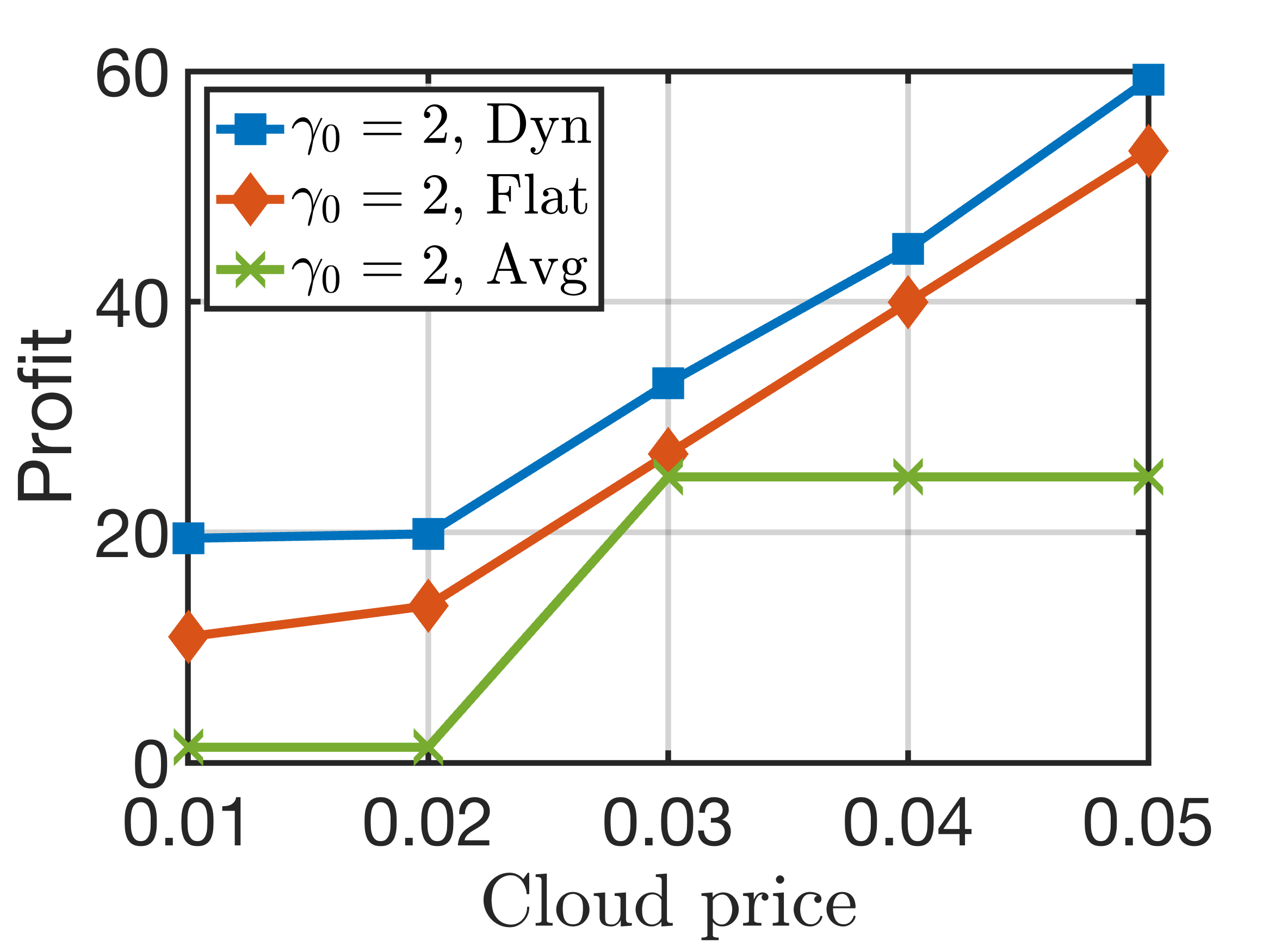}
	\label{fig:CompareCPAvg}
	}
	\hspace*{-1.8em} %\hspace\fill 
	\subfigure[Varying $M$]{
 \includegraphics[width=0.242\textwidth,height=0.11\textheight]{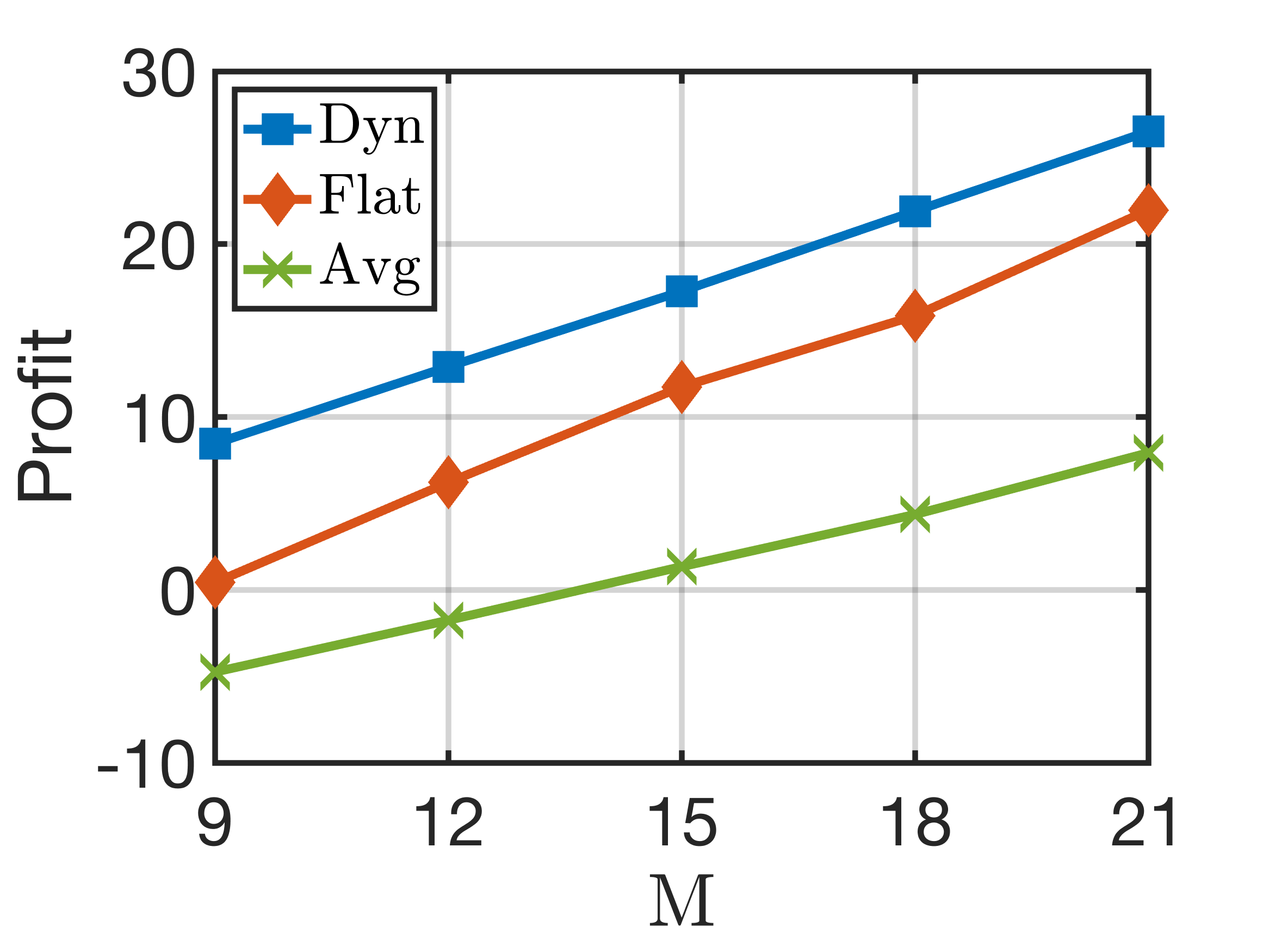}
	\label{fig:CompareM}
	}
	\vspace{-0.2cm}
	\caption{Performance comparison between Dyn, Flat, and Avg.}
	\label{fig:Compareschemes}
	\vspace{-0.2cm}
\end{figure}

%Based on our analysis, the proposed dynamic pricing scheme, as shown in 

Fig.\ref{fig:Compareschemes} illustrates that the proposed \textit{Dyn} scheme consistently outperforms the other pricing schemes. 
This can be attributed to the platform's ability to dynamically adjust prices for under-demanded ENs, thereby encouraging services to redistribute their workload to these nodes. As a result, resource utilization and platform revenue from these nodes improve. 
%This is primarily due to the ability of the platform to adjust prices for under-demanded ENs, which encourages services to reallocate their workload to these nodes. This results in improved resource utilization and revenue for the platform from these nodes. 
Furthermore, highly-demanded ENs are priced at maximum values, leading to higher profits. Increased resource procurement from ENs tends to decrease the average network delay for each service.
%Moreover, the pricing of highly-demanded ENs is often set to the maximum value, leading to %higher demand for edge resources and, in turn,
%higher profits. With increased resource procurement from ENs, the average network delay for each service tends to decrease.
In contrast, the \textit{Flat} pricing scheme, which sets the same price for all ENs, fails to incentivize services to shift their workload to under-demanded nodes. The Avg pricing scheme performs even worse since it employs a fixed price without considering the cost of cloud pricing. 
%since the price is fixed, without considering the cost of cloud pricing. We can also obverse from 
Fig.\ref{fig:Compareschemes} demonstrates that the platform's profit increases as the cloud resource price increases. This result is further supported in Fig.~\ref{fig:CompareCPAvg}. Additionally, as the number of service areas $M$ increases, resulting in a higher workload, the platform's profits also increase.

\begin{figure}[h!]
\vspace{-0.4cm}
    \centering
    \hspace*{-0.5em}
    \subfigure[Varying $M$ and demand $R$]{
	 \includegraphics[width=0.242\textwidth,height=0.11\textheight]{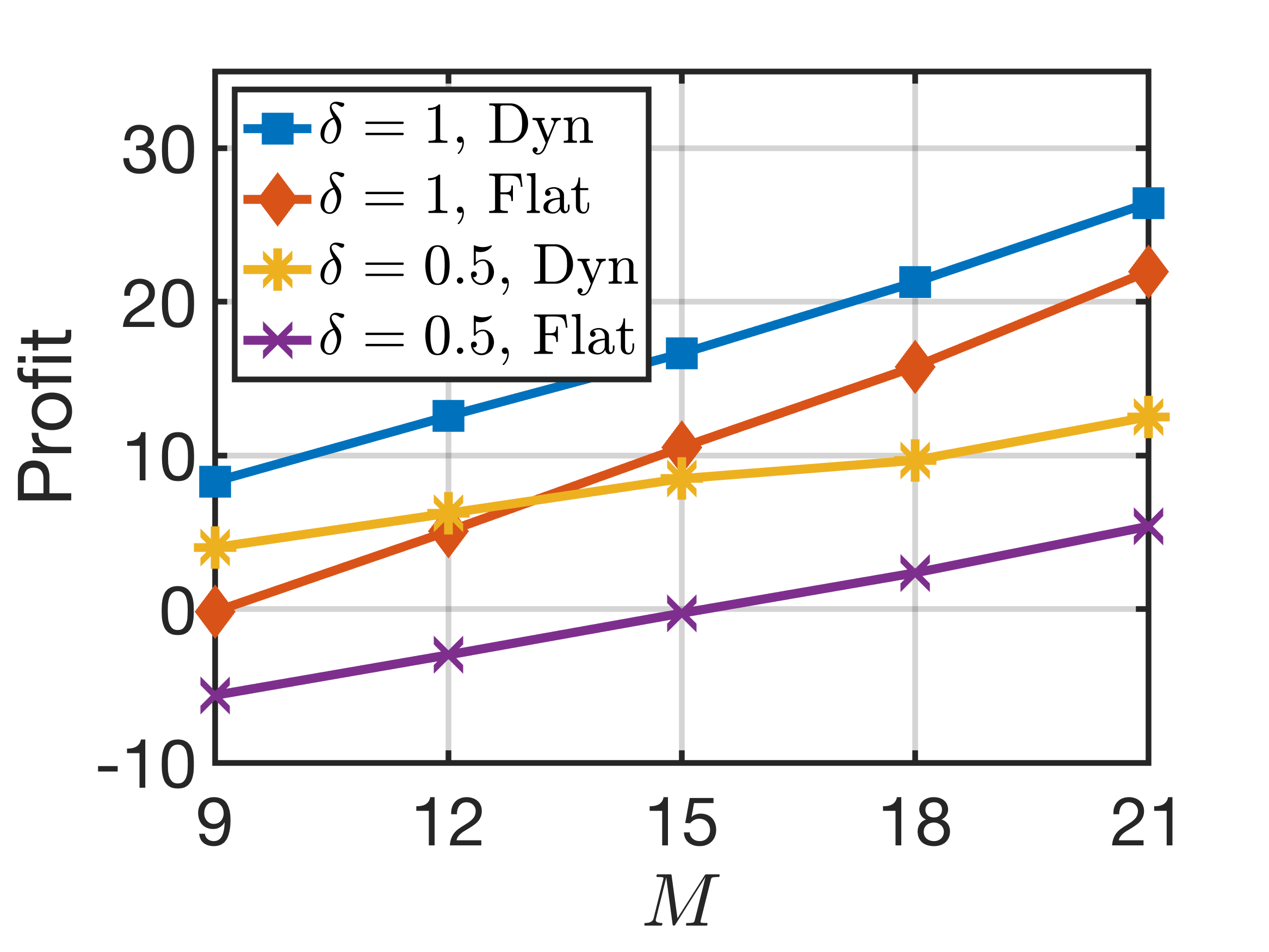}
	\label{fig:ENR}
	}
	\hspace*{-1.8em} %\hspace\fill 
	\subfigure[Varying $M$ and EN capacities]{
	 \includegraphics[width=0.242\textwidth,height=0.11\textheight]{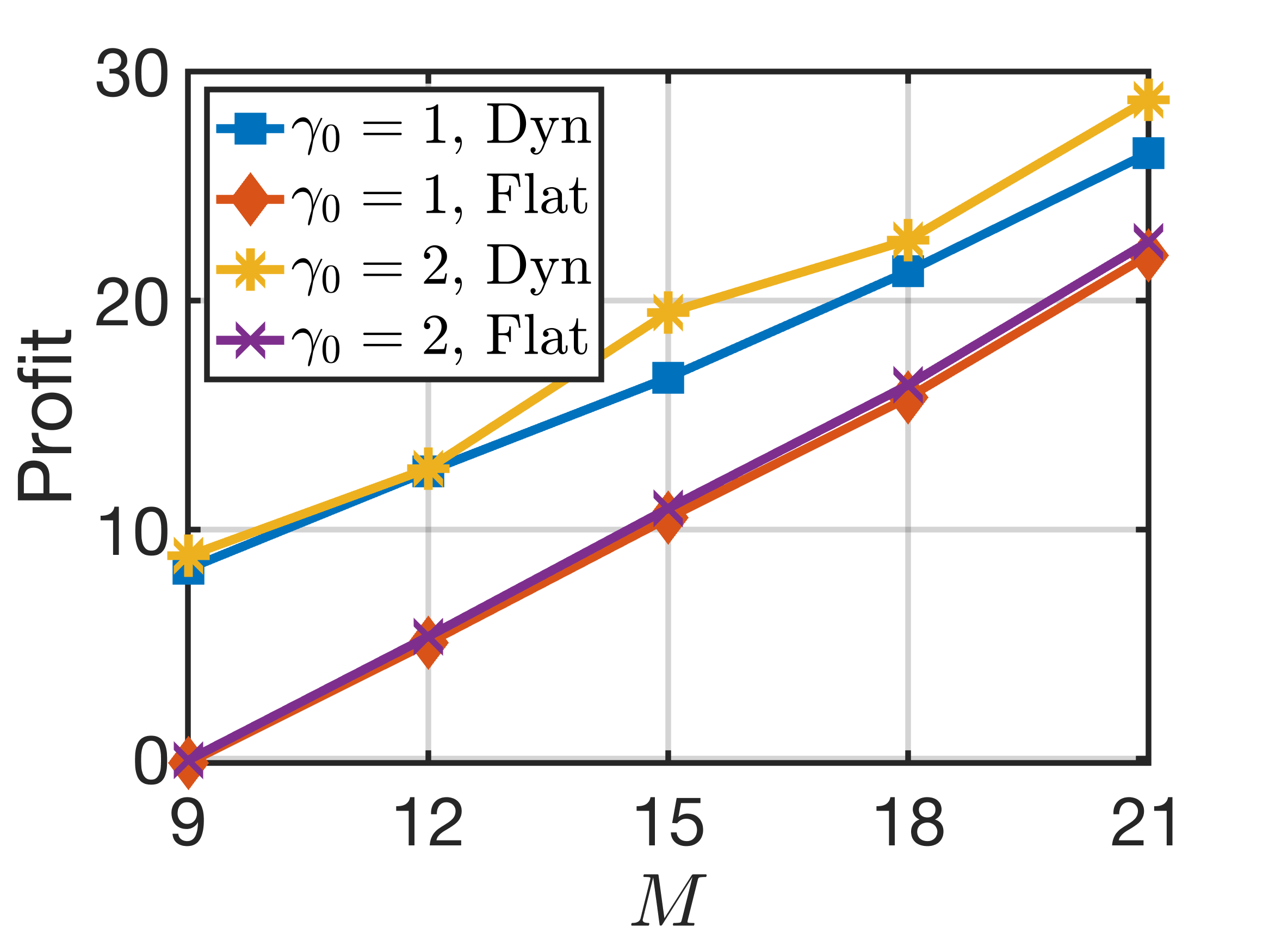}
	\label{fig:ENCap}
	}
	\vspace{-0.2cm}
	\caption{Impacts of number of areas.} % when $p_0=0.01$.}
	\label{fig:CompareEN}
	\vspace{-0.2cm}
\end{figure}

\begin{figure}[h!]
\vspace{-0.4cm}
    \centering
    \hspace*{-0.5em}
    \subfigure[Varying $\delta$, $p_0$ and $\gamma_0$]{
	 \includegraphics[width=0.242\textwidth,height=0.11\textheight]{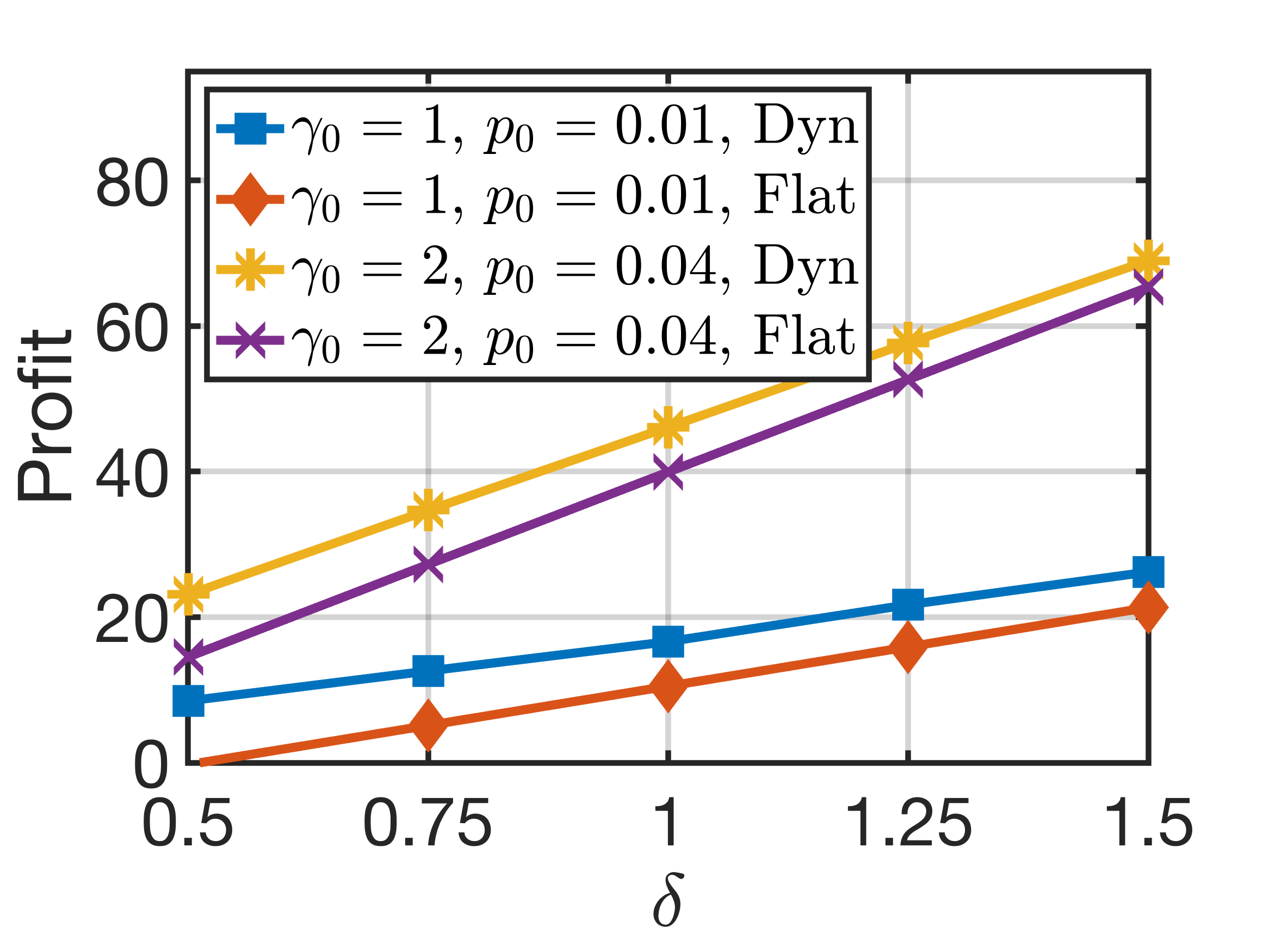}
	\label{fig:DemandCPCap}
	}
	\hspace*{-1.8em} %\hspace\fill 
	\subfigure[Varying $\delta$ and $\Lambda$]{
	 \includegraphics[width=0.242\textwidth,height=0.11\textheight]{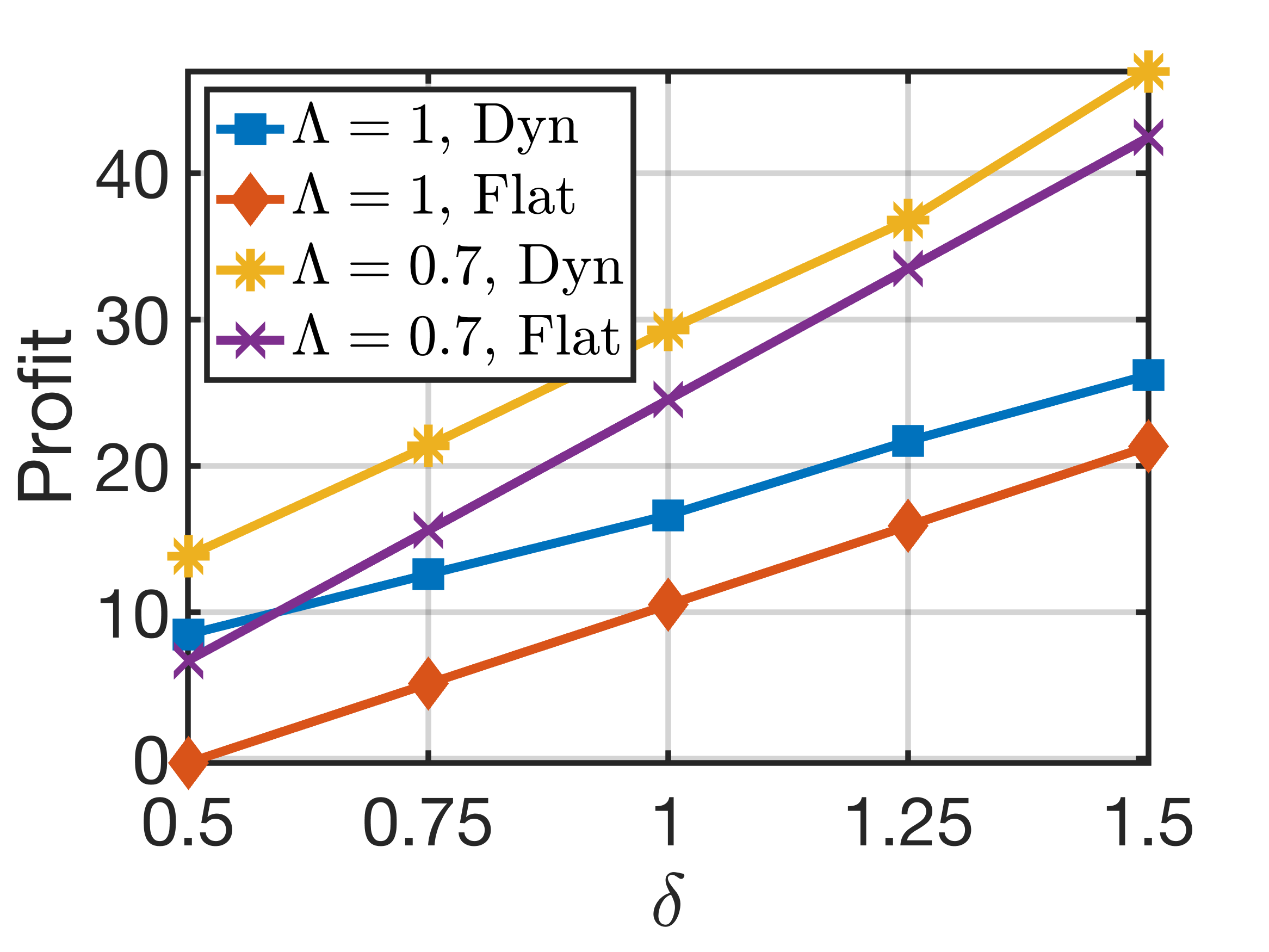}
	\label{fig:DemandDelay}
	}
	\vspace{-0.2cm}
	\caption{Impacts of the resource demand.}
	\label{fig:CompareDemand}
	\vspace{-0.5cm}
\end{figure}

Fig.\ref{fig:CompareEN} summarizes the impact of the number of areas on the system performance with varying demands and ENs' capacities by factors of $\delta$ and $\gamma_0$, respectively. These findings reinforce the superior performance of the proposed scheme. Moreover, an increase in the number of areas leads to an increased workload, resulting in increased profits.  Similarly, an increase in the demand or capacities of ENs further increases the profit, as services can acquire more resources from beneficial ENs. The impact of resource demands of services on system performance is illustrated in Fig.\ref{fig:CompareDemand}, where an increase in resource demands and capacities of ENs results in increased profits, consistent with the results observed in Fig.~\ref{fig:CompareEN}.

\begin{figure}[h!]
    \vspace{-0.4cm}
    \centering
    \hspace*{-0.5em}
    \subfigure[Varying $p_0$ and $\delta$]{
	 \includegraphics[width=0.242\textwidth,height=0.11\textheight]{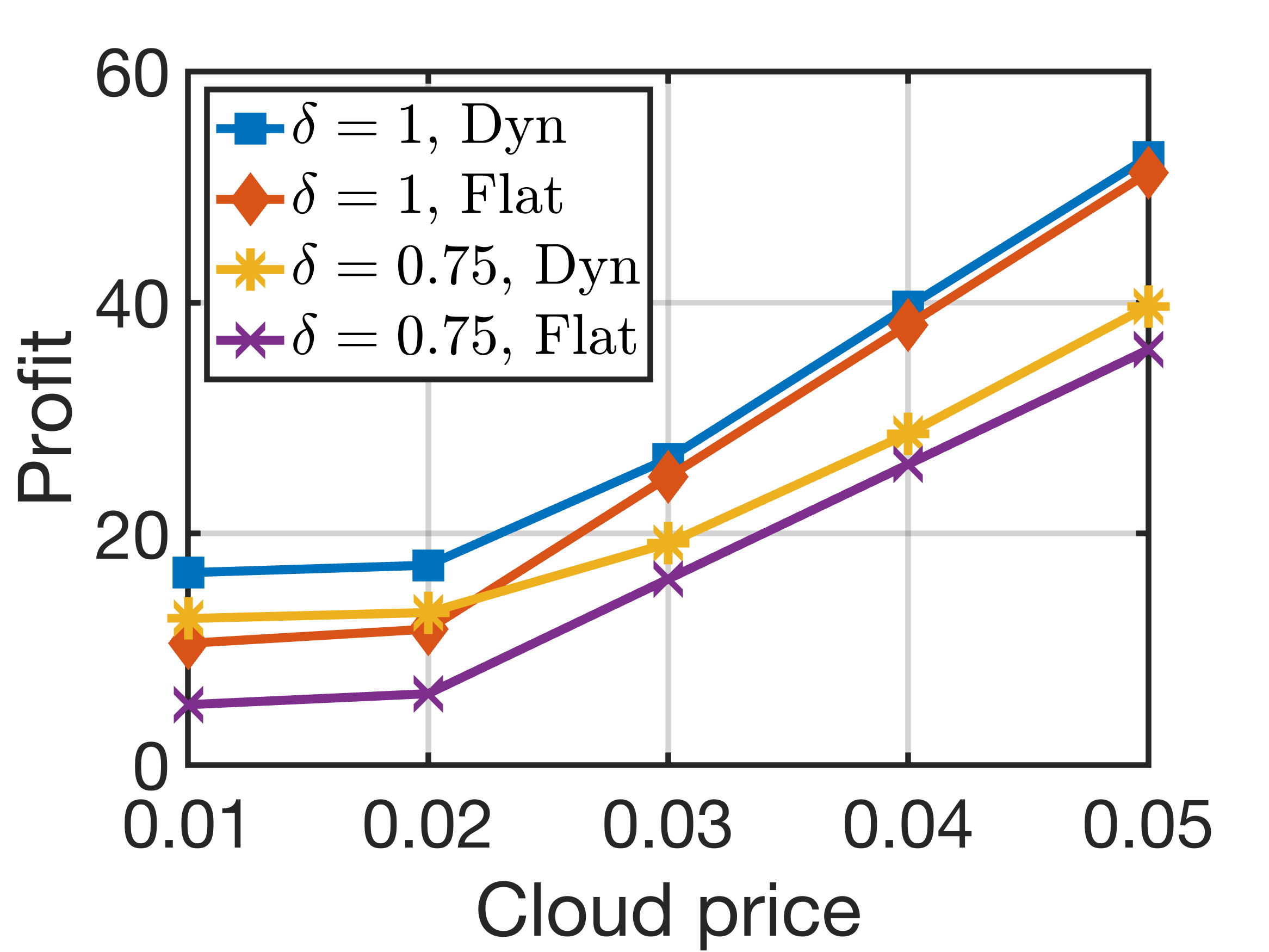}
	\label{fig:CPR}
	}
	\hspace*{-1.8em} %\hspace\fill 
	\subfigure[Varying $p_0$ and  $\gamma_0$]{
	 \includegraphics[width=0.242\textwidth,height=0.11\textheight]{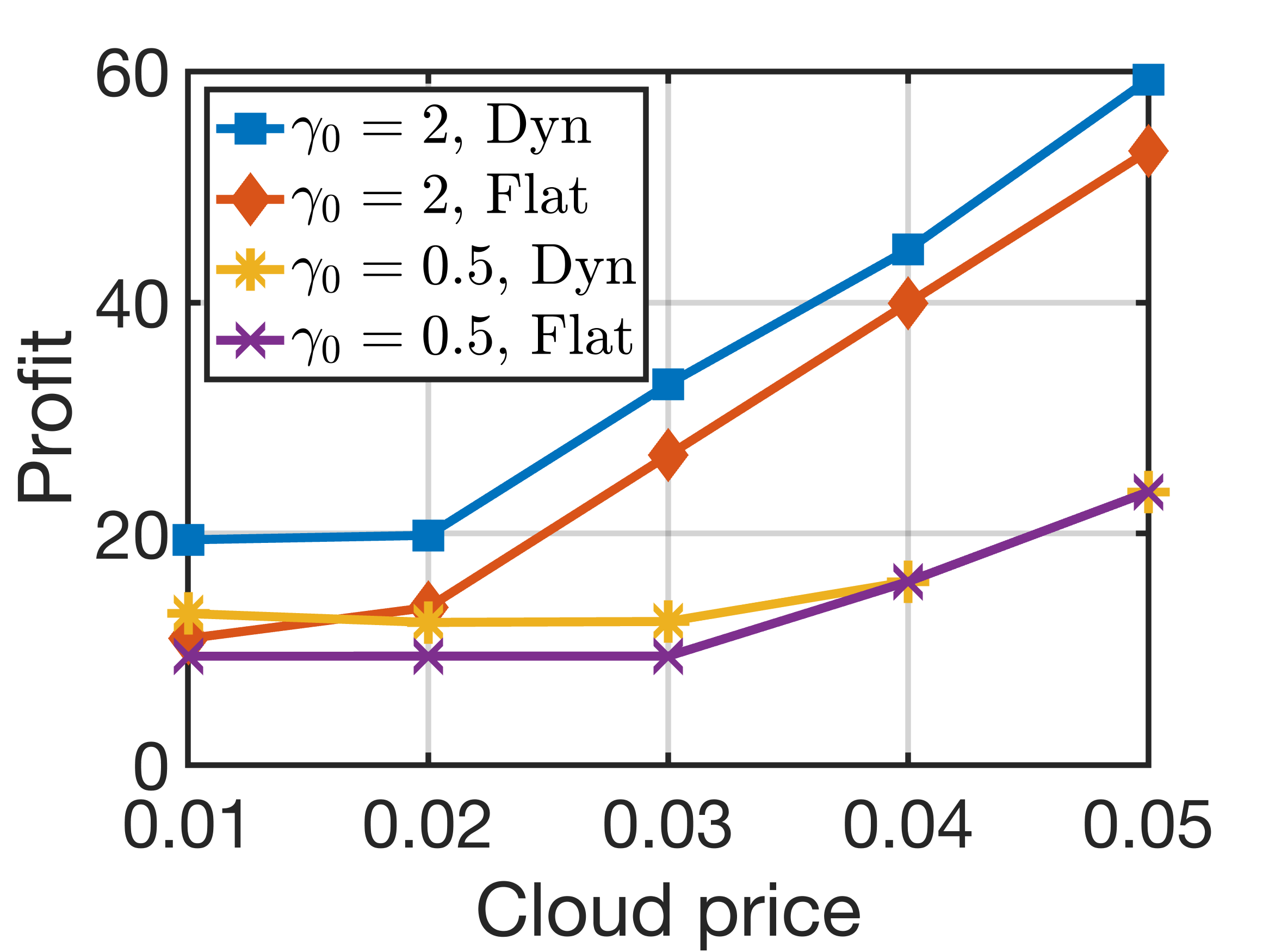}
	\label{fig:CPCap}
	}
	\vspace{-0.2cm}
	\caption{Impacts of cloud price $p_0$ on the system performance.}
	\label{fig:CompareCP}
	\vspace{-0.2cm}
\end{figure}

The impact of the cloud resource price on system performance is illustrated in Fig.\ref{fig:CompareCP}. %, where the profits are plotted against varying cloud resource prices. 
As shown in the figure, the platform's profit increases with higher cloud prices. %an increase in the cloud resource price leads to an increase in profit. 
This is because when the cloud resource price rises, services are motivated to offload more workload to ENs, leading to higher demand for edge resources and increased revenue and profit for the platform. It is worth noting that the difference in profit between the two pricing schemes becomes smaller as the cloud resource price increases. For instance, in Fig.\ref{fig:CPCap} where the cloud resource price is $p_0=0.05$, both \textit{Dyn} and \textit{Flat} pricing schemes result in the same profit. This outcome  is expected because when the cloud price is high, it is most profitable for the platform to set the price at all ENs equal to the cloud price, resulting in the dynamic pricing scheme behaving similarly to the flat pricing scheme.

\begin{figure}[h!]
\vspace{-0.4cm}
    \centering
    \hspace*{-.5em}
    \subfigure[Varying $\Lambda$ and $p_0$]{
	 \includegraphics[width=0.242\textwidth,height=0.11\textheight]{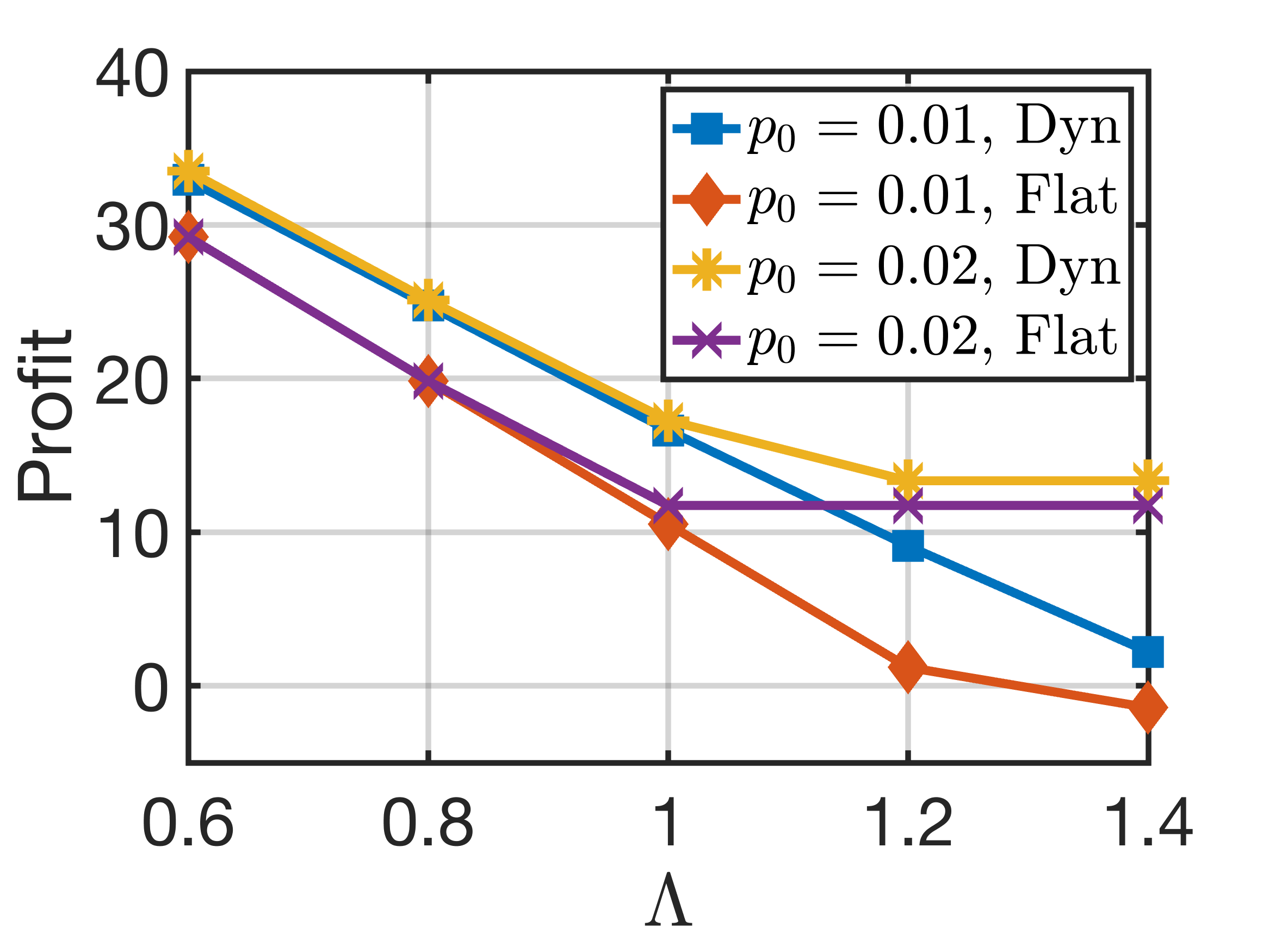}
	\label{fig:DelayCP}
	}
	\hspace*{-1.8em}
	\subfigure[Varying $\Lambda$ and $\delta$]{
	 \includegraphics[width=0.242\textwidth,height=0.11\textheight]{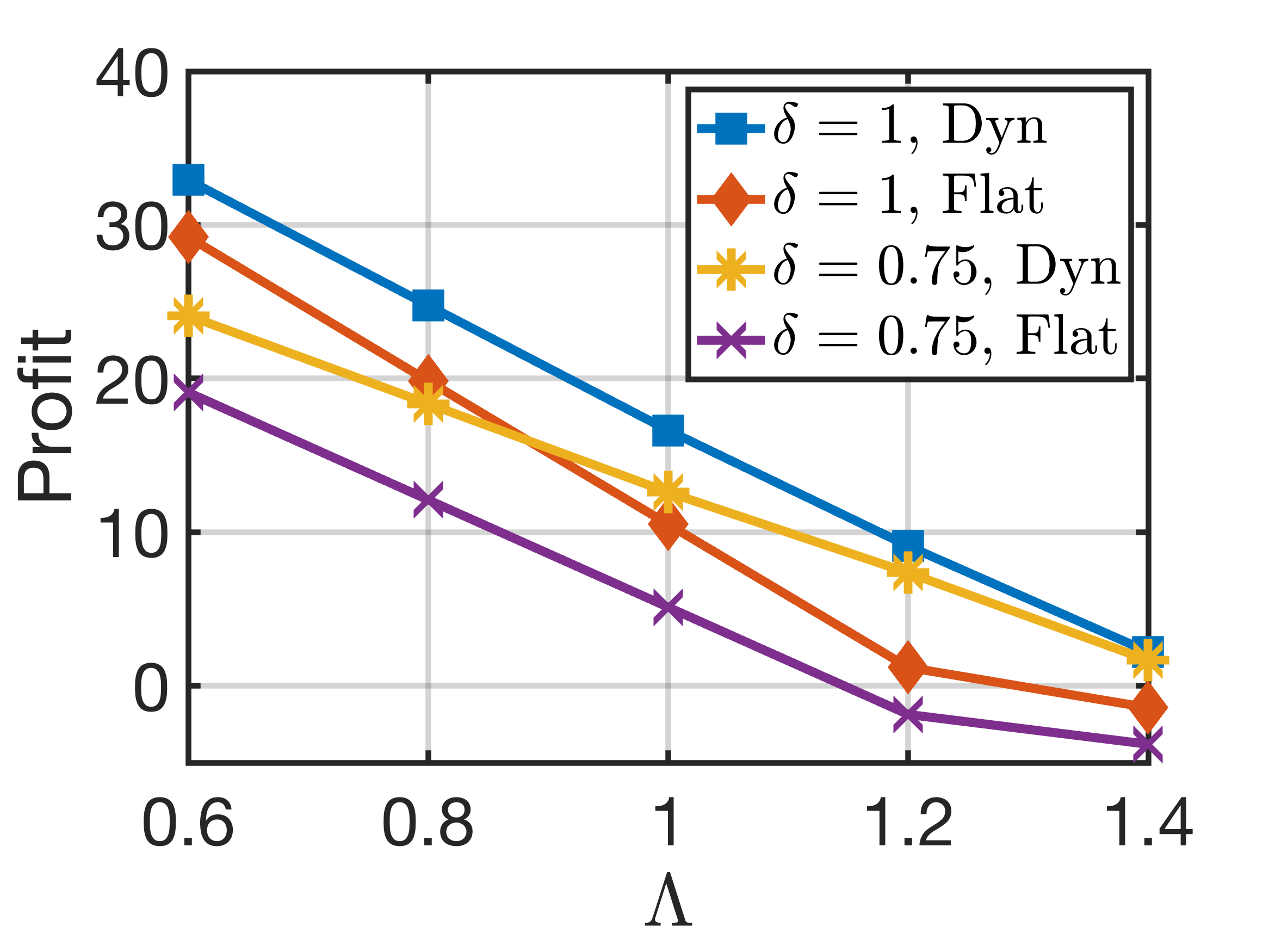}
	\label{fig:DelayR}
	}
	\vspace{-0.2cm}
	\caption{Impacts of the delay threshold.}
	\label{fig:CompareDelay}
	\vspace{-0.2cm}
\end{figure}

Fig.\ref{fig:CompareDelay} presents the impact of the delay threshold $D^{k,m}$ on the platform's profit.   Recall that $\Lambda$ is the scaling factor of the delay threshold. The figure demonstrates that the proposed scheme outperforms the flat pricing scheme. % the superior performance of the dynamic pricing scheme compared to the flat schemes.
Also, the profit increases as the delay threshold parameters decrease. This is because as services become more delay-sensitive, they prioritize placing services on multiple ENs that are close to users to minimize network delay, 
even if it results in higher placement costs. 
%even if it means incurring higher placement costs. %However, selecting appropriate ENs for service placement can be challenging due to ENs' heterogeneity and service-specific specifications.
Consequently, the platform can increase the edge resource prices to increase its profit.

\begin{figure}[h!]
\vspace{-0.4cm}
    \centering
    \hspace*{-0.5em}
    \subfigure[$p_0=0.02$, $\beta_0=0.5$]{
	\ \includegraphics[width=0.242\textwidth,height=0.11\textheight]{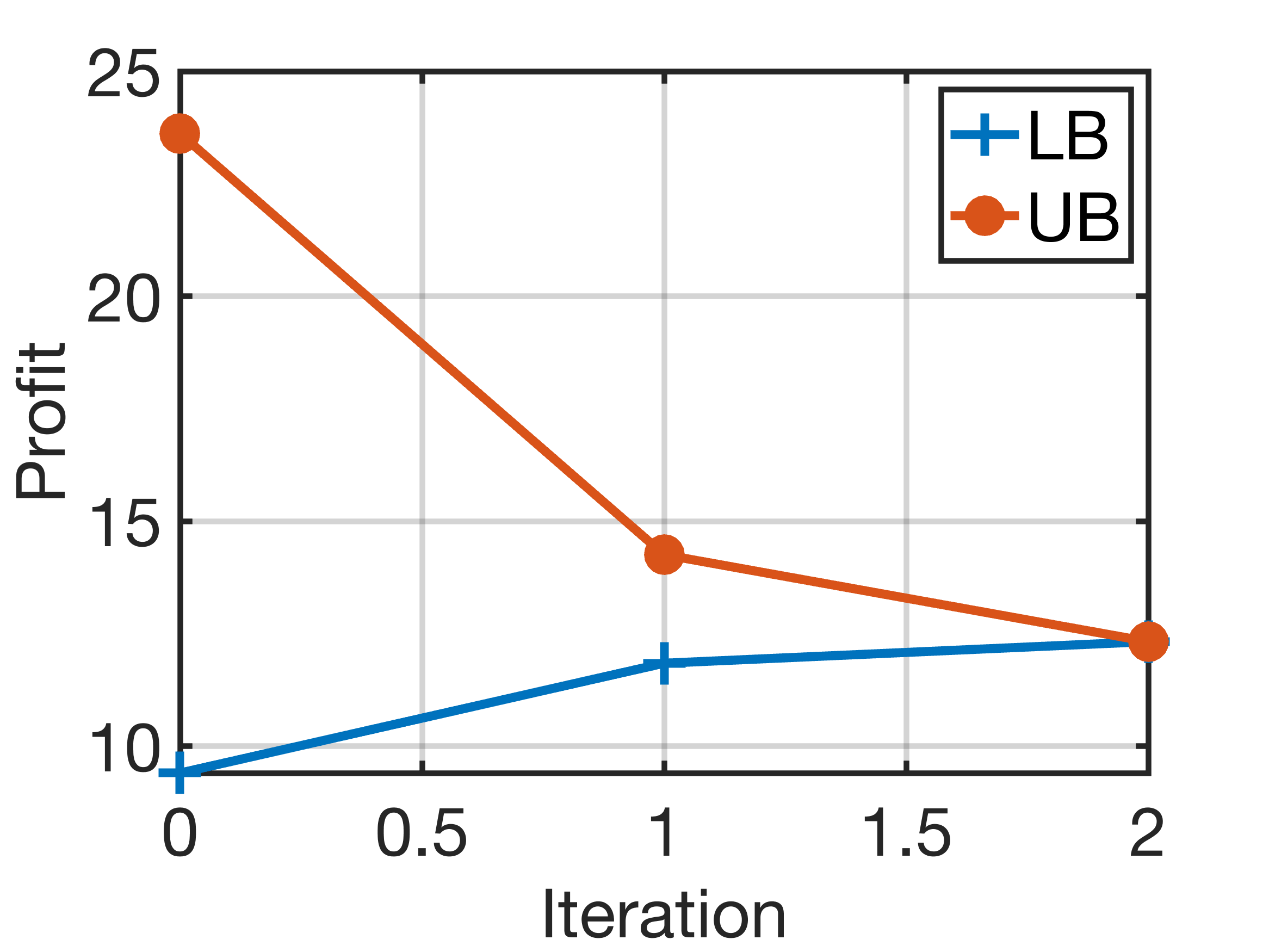}
	\label{fig:Conv1}
	}
	\hspace*{-2.1em} %\hspace\fill 
	\subfigure[$p_0=0.04$, $\beta_0=0.5$]{
	 \includegraphics[width=0.242\textwidth,height=0.11\textheight]{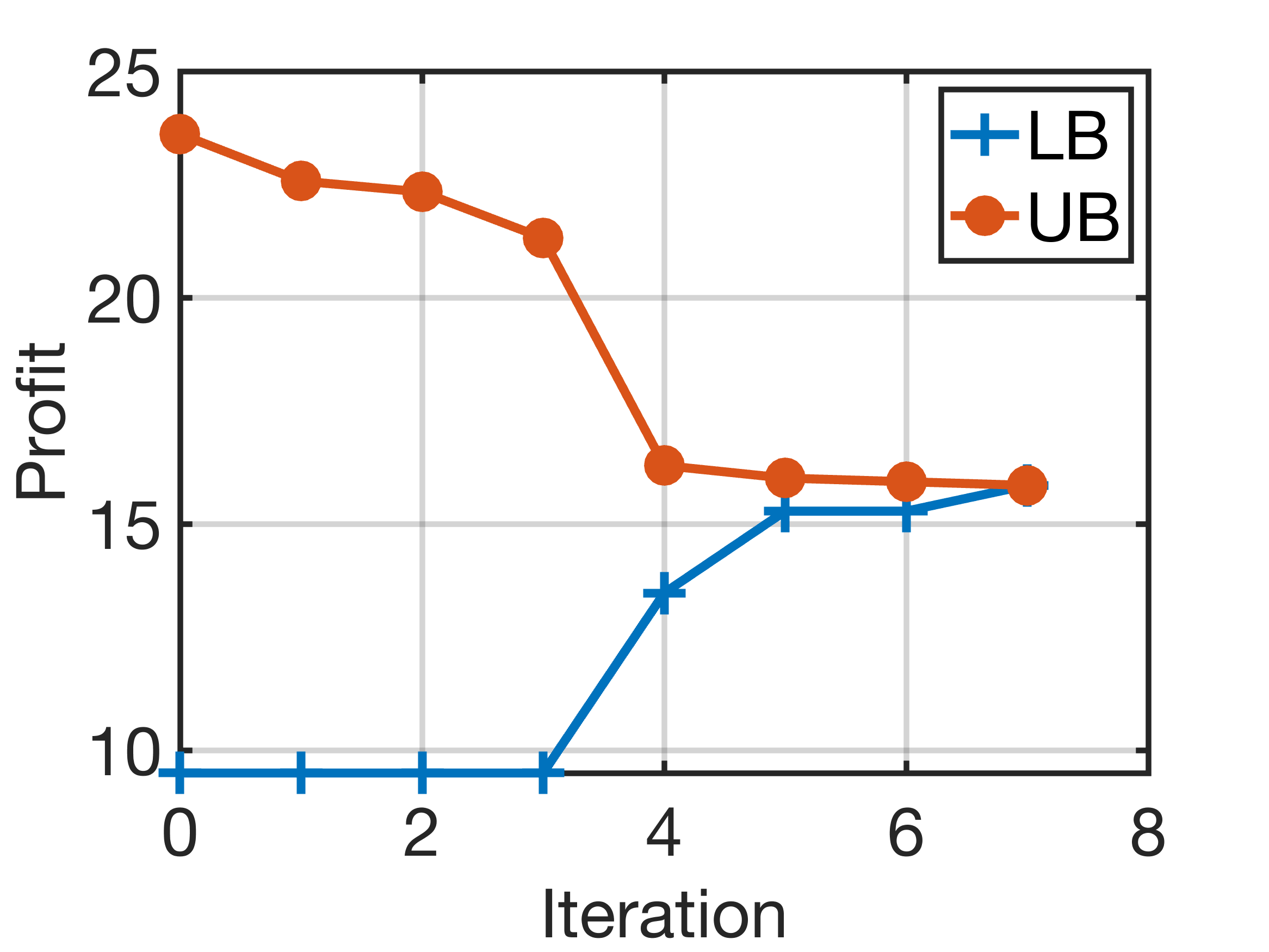}
	\label{fig:Conv2}
	}
	\hspace*{-0.5em}
    \subfigure[$I=18$, $\delta_0=0.5$]{
	 \includegraphics[width=0.242\textwidth,height=0.11\textheight]{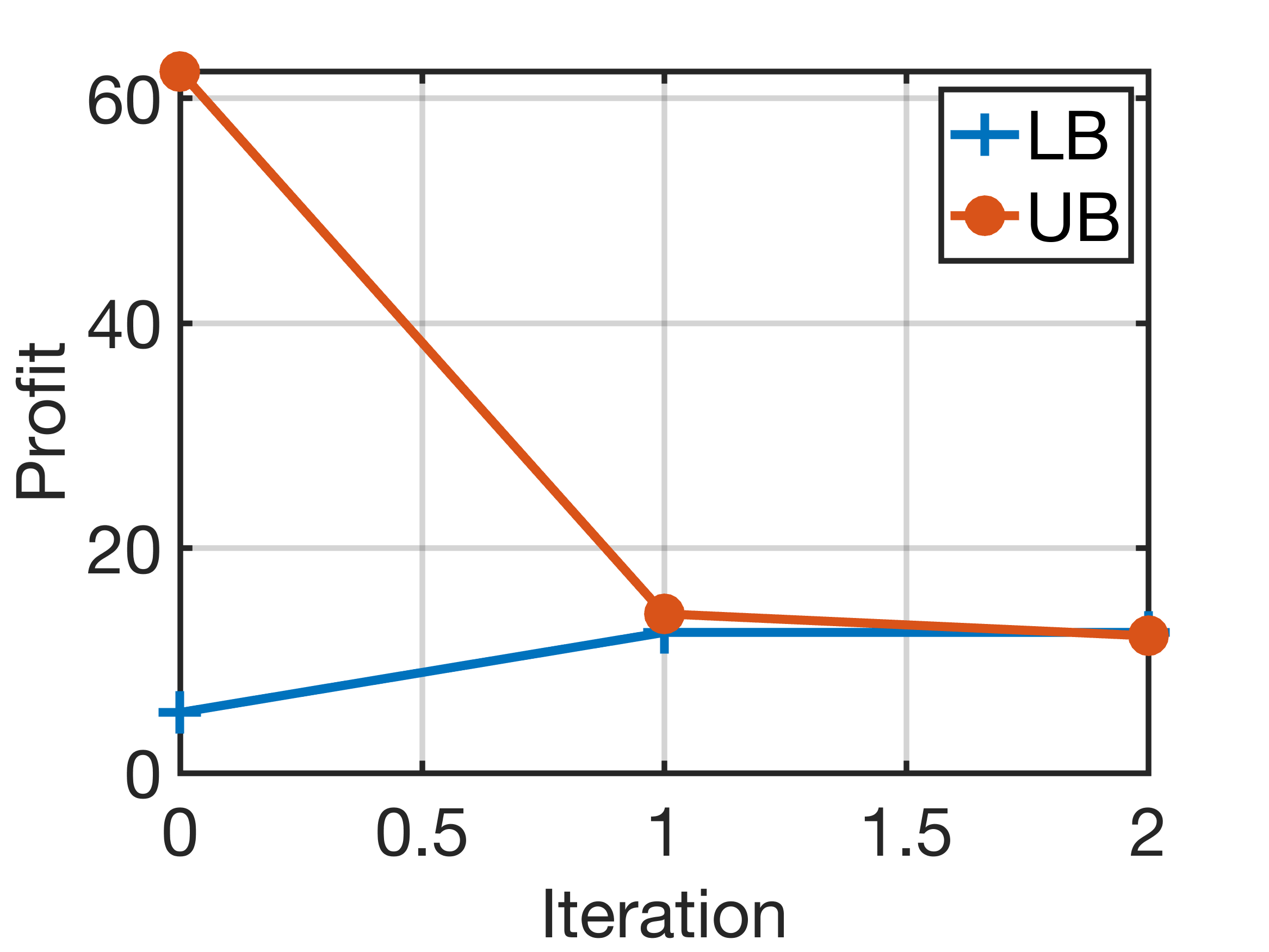}
	\label{fig:Conv3}
	}
	\hspace*{-2.1em} %\hspace\fill 
	\subfigure[$I=21$, $\gamma_0=2$]{
	 \includegraphics[width=0.242\textwidth,height=0.11\textheight]{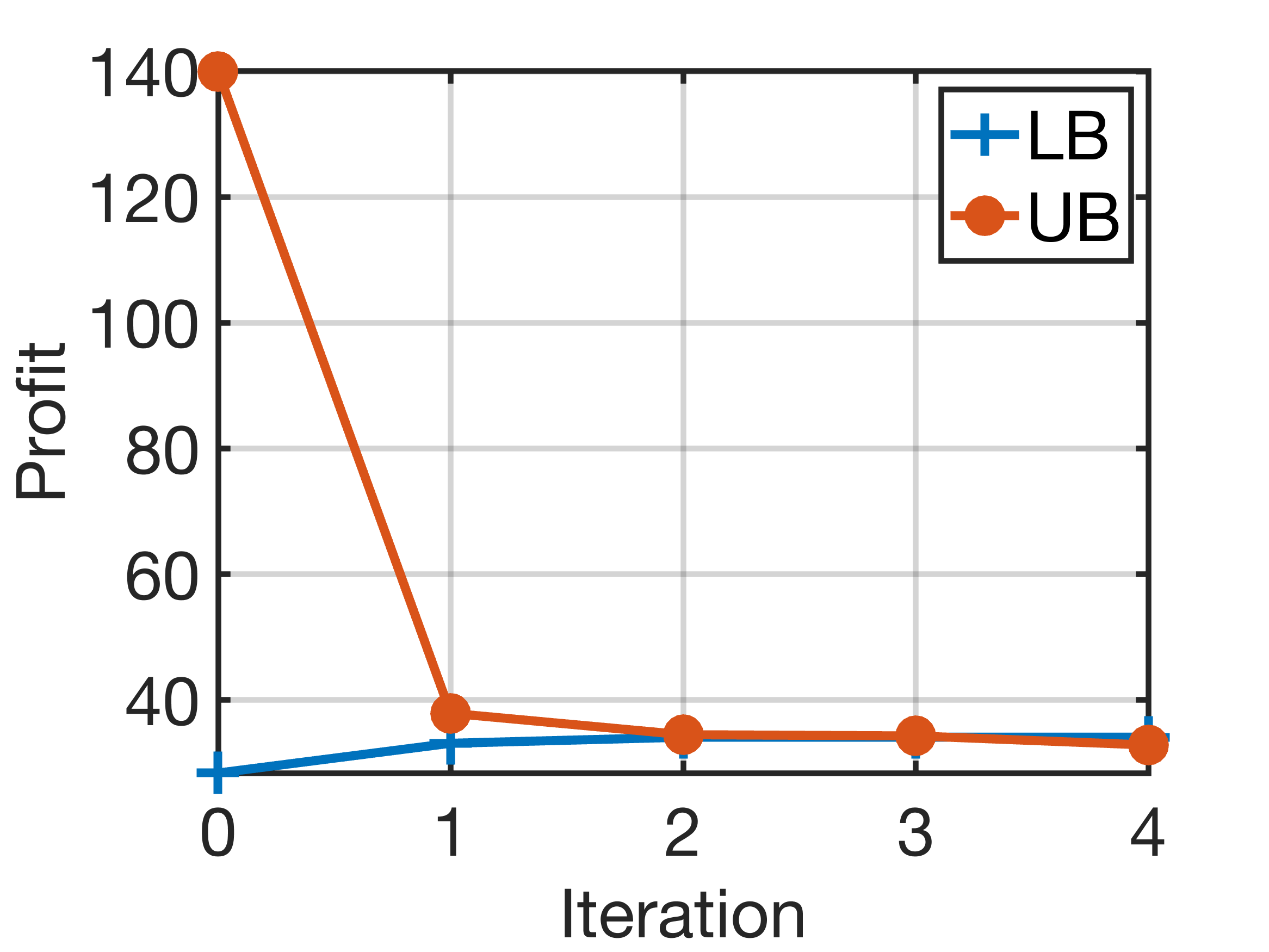}
	\label{fig:Conv4}
	}
	\vspace{-0.2cm}
	\caption{Convergence property.}
	\label{fig:ConvPlots}
	\vspace{-0.4cm}
\end{figure}

Finally, Fig.\ref{fig:ConvPlots}
illustrates the convergence property of the proposed algorithm for different problem instances.
%we demonstrate the convergence property of the proposed algorithm in Fig.\ref{fig:ConvPlots} for different problem instances. 
Our numerical experiments demonstrate that the algorithm converges rapidly towards optimal solutions in just a few iterations for various system sizes.

\section{Related Work}
\label{rel} 
%Edge and cloud computing are important network architectures for their potential to deliver scalable and efficient resource provisioning across various applications.

The existing literature on service placement and workload scheduling in EC is extensive. In \cite{duong20}, a two-stage optimization framework is proposed to optimize service placement and sizing decisions for an SP, considering uncertain demand. Reference \cite{mjia17}  utilizes queuing models to jointly optimize cloudlet placement and workload allocation decisions, aiming to minimize system response time with a fixed number of cloudlets. 
%In [25], the authors propose a ranking-based heuristic algorithm for efficient cloudlet deployment in an IoT network. 
Reference \cite{syan19} focuses on optimizing cloudlet placement and task allocation to minimize energy consumption under %while adhering to
delay constraints.
%In [9], the authors present a two-timescale optimization framework to optimize service placement and request scheduling under the budget and multi-dimensional resource constraints. 
%In [10], A. Yousefpour et al. introduce an edge service provisioning model to minimize the total system cost by dynamically deploying and releasing applications on different ENs. The joint service placement and request routing in mobile edge computing (MEC) is investigated in [11] to minimize the workload to the cloud, considering the asymmetric bandwidth requirements of the services and the limited storage capacities of ENs. 
In \cite{lwan18}, Wang \textit{et al. }
address the service placement problem for social VR applications, aiming to minimize the total application deployment cost. However, none of these works consider the resource pricing design issue.
%examine the service placement problem for social VR applications to minimize the total application deployment cost. 
%, including the cloudlet activation, service placement, proximity, and colocation costs. 
%However, the resource pricing design issue is not considered in these works. 

Considerable research has also been conducted on pricing design in both cloud and edge networks.
%A substantial amount of research has also been carried out on pricing design in cloud and edge networks. 
In \cite{hxu13}, Xu \textit{et al.} propose a revenue maximization model and employ stochastic dynamic programming to address the dynamic pricing problem in an IaaS cloud. Similar research on joint virtual machine pricing, task scheduling, and server provisioning is explored in \cite{jzha14} using an online profit maximization algorithm.   Fisher market models are used in \cite{duongTCC,duongTON} to design fair pricing and allocation of edge resources. In \cite{ChengBandit}, a bandit algorithm is proposed to tackle the online pricing problem for heterogeneous edge resources.  However, these works primarily focus on centralized decision-making models. %did not consider the hierarchical decision making n

%Reference [28] studies the problem of resource pricing by thoroughly analyzing several dynamic pricing schemes based on auctions and fairness-seeking properties from the perspective of game-theory and existence of unique Nash or Stackelberg equilibrium. %In [29], the Lagrange multiplier method and a dynamic closed loop control scheme are integrated to solve the user perceived value-based dynamic pricing problem.

To address the hierarchical decision-making challenges related to resource allocation and pricing in both cloud and EC, several Stackelberg games and bi-level optimization models have been proposed.
%To address hierarchical decision-making problems for resource allocation and pricing in cloud and EC, various Stackelberg games and bi-level optimization models  have been proposed.
%-theoretical model has been applied in several studies to tackle various resource allocation and pricing problems.
%Stackelberg games and bi-level optimization have been proposed for studying resource allocation and pricing in cloud and edge computing. 
%In [30], the authors introduce a bi-level model and a heuristic algorithm to study the task allocation problem in a two-layer multi-community cloud/cloudlet social collaborative computational framework. 
% Reference \cite{mliu18} proposes a Stackelberg game between a single EN and multiple
% mobile users, in which the former seeks to maximize revenue
% within capacity constraints, while the latter seeks to minimize
% cost performing optimal task allocation. 
 Reference \cite{mliu18} models the interaction between edge clouds and users as a Stackelberg game, optimizing revenue and local offloading decisions.
Reference \cite{valerio20} addresses the problem of multi-service resource provisioning and pricing with fairness in a single-leader-multiple-followers Stackelberg game.  In \cite{yswang17}, a multi-leader-multi-follower two-stage Stackelberg game-based dynamic allocation algorithm is proposed for optimal bandwidth distribution to mobile terminals. \cite{wzhang20} formulates service pricing and offloading of IoT applications in a multi-MEC system are modeled as a Stackelberg game.   In \cite{phuang19}, the authors
present a bi-level optimization model, in which the upper-level
model represents the task allocation problem and the lower level captures the resource allocation problem, to minimize
energy consumption under delay constraints.

However, most existing models simplify follower problems by excluding integer variables or employing overly simplistic cost functions and constraints. These simplifications enable the derivation of closed-form solutions and facilitate the use of backward induction methods. \textit{In contrast, our model incorporates realistic costs and constraints, which introduces significant complexity. This results in a non-convex problem that lacks a closed-form solution. The need for discrete service placement decisions further complicates the follower problems and requires more sophisticated solution methods.} To the best of our knowledge, our model is the first to address both computing and storage resource pricing in EC while accounting for service placement costs in the follower problems.

\section{Conclusion}
\label{conc}
In this work, we propose a novel bi-level optimization model for joint service placement and edge resource pricing. The platform is responsible for optimizing the EN activation and resource pricing at different ENs in the upper level while 
the lower level involves each service making decisions regarding optimal service placement, resource procurement, and workload allocation.
%each service decides on optimal service placement, resource procurement, and workload allocation decisions in the lower-level.
%The presence of integer variables in the lower-level formulation makes solving the proposed \bmip ~ problem challenging. 
Due to the integer variables in the lower level, % of the formulation,
the proposed \bmip ~problem is challenging to solve.
The \bmip~ is converted into a constrained mathematical program through a reformulation procedure, which involves introducing duplicated lower-level variables, exploring potential values for lower-level integer variables to derive valid inequalities, and applying KKT conditions or LP duality, complemented by linearization techniques. We further developed %are inspired by the iterative decomposition algorithm CCG 
a decomposition-based iterative algorithm for solving the problem in a master-subproblem framework with finite convergence guarantee.  %that converges in finite iterations. 
Finally, extensive simulations are conducted to evaluate the performance of the proposed model and solutions. Our findings indicate that our proposed methodology outperforms existing approaches and provides an exact optimal solution to the joint optimal edge resource pricing and service placement problem.  Our results also indicate that the resource prices tend to decrease as the platform charges the service placement costs.  %This study contributes to the optimization of edge computing platforms and provides a promising approach to addressing the joint optimal pricing and service placement problem. 
%We present extensive numerical results to demonstrate the efficacy of our proposed model and technique. Our results reveal that the proposed methodology outperforms existing approaches and provides an exact optimal solution to the joint optimal resource pricing and service placement problem. This study contributes to the optimization of edge computing platforms and provides a promising approach to address the joint optimal resource pricing and service placement problem.

\bibliographystyle{IEEEtran}

\begin{IEEEbiography}[{\includegraphics[width=.95in,clip,keepaspectratio]{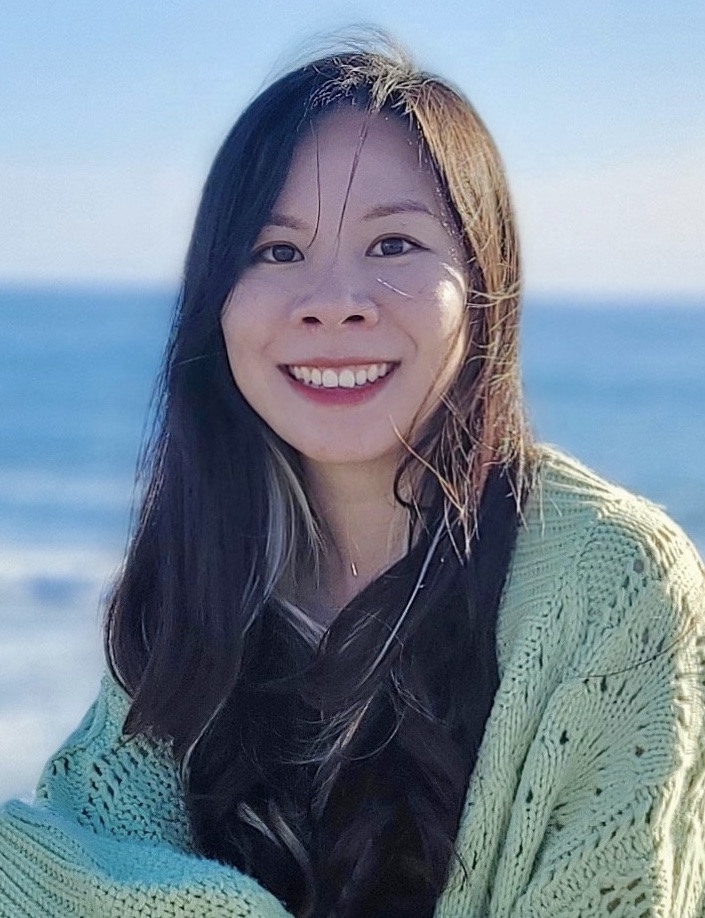}}]
{Duong Thuy Anh Nguyen} received the M.Sc. degree in Applied Mathematics from the University of Louisiana at Lafayette, LA, USA in 2019. She is currently pursuing the Ph.D. degree with the School of Electrical, Computer and Energy Engineering at Arizona State University, AZ, USA. Her current research interests include distributed optimization, operations research, and game theory, with a focus on developing mathematical models for decision-making under uncertainty and privacy-preserving %mechanism designs and 
distributed algorithms in multi-agent systems. 
	\end{IEEEbiography}

\vspace{-0.5cm}
\begin{IEEEbiography}[{\includegraphics[width=.95in,clip,keepaspectratio]{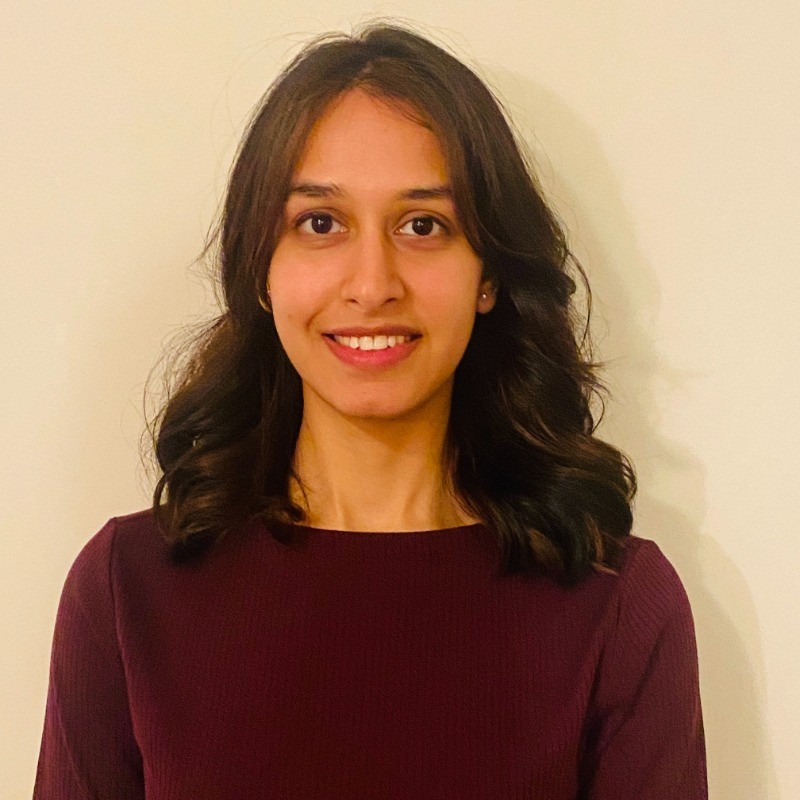}}]
{Tarannum Nisha} received the M.Sc. degree in electrical and computer engineering from the University of British Columbia, BC, Canada in 2021. She is currently a Data Scientist at a Vancouver-based startup, where she applies advanced quantitative techniques to real-world problems. Her research and professional interests span mathematical modeling, multi-agent systems, operations research, game theory, market design, and artificial intelligence/machine learning.
	\end{IEEEbiography}
    
\vspace{-0.5cm}
\begin{IEEEbiography}[{\includegraphics[width=1in,height=1.25in,clip,keepaspectratio]{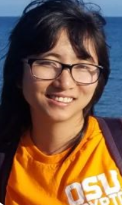}}]
    {Ni Trieu} received the B.Sc. degree in computer science from St. Petersburg State Polytechnic University, St. Petersburg, Russia, in 2013, and the master's and Ph.D. degrees from Oregon State University, Corvallis, OR, USA, 2017 and 202, respectively. She was a Postdoctoral Researcher with the University of California, Berkeley, CA, USA, in 2020. Since 2020, she has been with Arizona State University, Tempe, AZ, USA, where she is currently an Assistant Professor with the School of Computing, Informatics, and Decision Systems. Her research interests are in the area of cryptography and security, with a specific focus on secure computation and its applications. %, such as private database queries, contact tracing, biocomputing, privacy-preserving machine learning, data privacy, and secure optimization.
    \end{IEEEbiography}
    
\vspace{-0.5cm}
\begin{IEEEbiography}[{\includegraphics[width=1in,height=1.25in,clip,keepaspectratio]{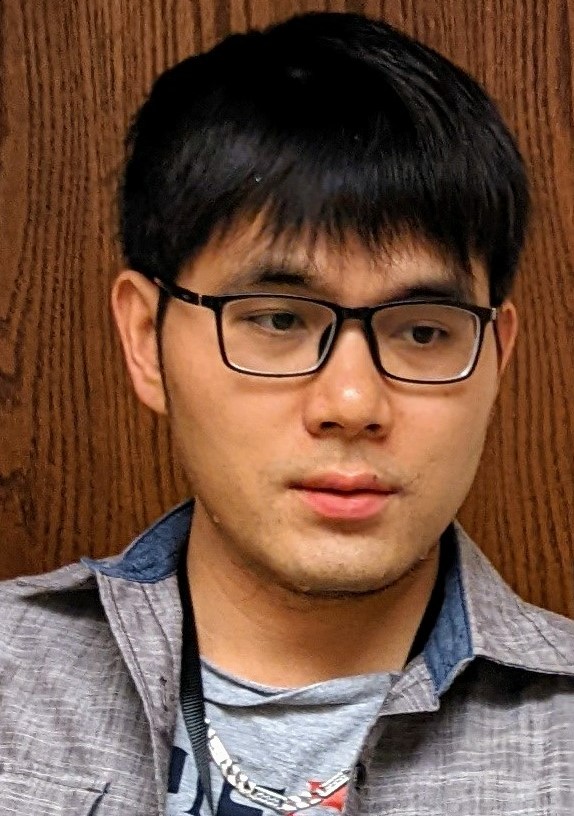}}]
{Duong Tung Nguyen} received the Ph.D. degree in electrical and computer engineering from the University of British Columbia, BC, Canada. He is currently an assistant professor in the School of Electrical, Computer and Energy Engineering at Arizona State University, AZ, USA. His research lies at the intersection of operations research, AI, economics, and engineering, focusing on developing mathematical models and techniques for decision-making and analysis of large-scale networked systems. 
\end{IEEEbiography}

\cleardoublepage

\input{appendix.tex}
\end{document}

%% file: appendix.tex
\appendix
\setcounter{page}{1}
\subsection{Feasible Set $\mathcal{S}^k(\bp,\bp^s,\bz,\bt^k)$} \label{app-FeasibleSetSk}
The feasible set $\mathcal{S}^k(\bp,\bp^s,\bz,\bt^k)$ associated with the second minimization problem in \eqref{sepenum} is provided below:
\begin{subequations}
\begin{alignat}{3}
\!\!p_0 y_0^k + \!\!\sum_j p_j y_j^k \leq B^k \!- \!\!\sum_j \left(\phi_j^k+s^kp_j^s\right) t_j^k&,&&(\mu_1^k)\\ 
t_j^k - z_j \leq 0 &, \forall j, &&(\nu_j^k)\!\!\\
\sum_i x_{i,j}^k - y_j^k \leq 0&, \forall j, &&(\Gamma_j^k)\!\!\\
\sum_i x_{i,0}^k - y_0^k \leq 0&, &&(\mu_2^k)\!\!\\
y_j^k \leq C_j t_j^k &, \forall j, &&(\sigma_j^k)\!\!\\
x_{i,0}^k d_{i,0} + \sum_j x_{i,j}^k d_{i,j} \le D^{k,m} R_i^k &,\forall i, &&(\xi_i^k) \\
\sum_{j} x_{i,j}^k +  x_{i,0}^k  + q_i^k = R_i^k &, \forall i, &&(\eta_i^k) \\
x_{i,j}^k  \leq a_{i,j}^k R_i^k &, \forall i,j,~ &&(\tau_{i,j}^k)\!\!
\end{alignat}
\end{subequations}
where $\mu_1^k$, $\nu_j^k$, $\Gamma_j^k$, $\mu_2^k$, $\sigma_j^k$, $\xi_i^k$ and $\tau_{i,j}^k$ are non-negative continuous dual variables, while $\eta_i^k \in \mathbb{R}$. It is important to note that $\bp$, $\bp^s$, $\bz$, and $\bt^k$ represent the parameters in this context.

\subsection{\textbf{BMIP} with Duplicated Lower-Level Variables}\label{app:bmipd}
The explicit formulation for \bmipd, as presented in \eqref{eq:bmipd}, is as follows:
\begin{align*}
&\textbf{\bmipd:}~~ \max_{\substack{\bp,\bp^s,\bz,\\\bx^\prime,\bq^\prime,\by^\prime,\bt^\prime}} ~ \sum_j p_j \sum_k y_j^{\prime k} + \sum_j\sum_k \left(\phi_j^k+s^kp_j^s\right) t_j^{\prime k} \nonumber\\
&\qquad\qquad\qquad- \sum_j \left(f_j z_j + c_j  \frac{\sum_k y_j^{\prime k}}{C_j} \right)  \\
&\!\!\text{subject to:}\\
% (\bp,\bp^s,\bz,\bx',\bq',\by',\bt') \in \mathcal{H}, \\
& \!\!\sum_k y_j^{\prime k} \leq z_j C_j; \sum_k s^k t_j^{\prime k} \leq z_jS_j; y_0^{\prime k} \geq \!\sum_i x_{i,0}^{\prime k}, \forall j,k,\!\\
&p_0 y_0^{\prime k} \!+ \!\!\sum_j p_j y_j^{\prime k} \!+ \!\!\sum_j \left(\phi_j^k+s^kp_j^s\right) t_j^{\prime k}  \leq \!B^k, ~~\forall k, \\
&t_j^{\prime k} \leq z_j; y_j^{\prime k} \geq \sum_i x_{i,j}^{\prime k};y_j^{\prime k} \leq C_j t_j^{\prime k};s^k t_j^{\prime k} \leq S_j,\forall j,k,\! \\
&x_{i,j}^{\prime k} \leq a_{i,j}^k R_i^k;~\sum_j x_{i,j}^{\prime k} + x_{i,0}^{\prime k} +  q_i^{\prime k} = R_i^k, ~\forall i,k, \\
&x_{i,0}^{\prime k} d_{i,0} +  \sum_j x_{i,j}^{\prime k} d_{i,j} \leq D^{k, \sf m} R_i^k, ~~\forall i,k, \\
& w^k \Big( \sum_i x_{i,0}^{\prime k} d_{i,0} + \sum_{i,j} x_{i,j}^{\prime k} d_{i,j} \Big) + p_0 y_0^{\prime k} + \sum_j p_j y_j^{\prime k}\nonumber\\
& +  \sum_i \psi_i^k q_i^{\prime k}+ \sum_j (\phi_j^k   +s^kp_j^s) t_j^{\prime k} \leq  \min_{\substack{(\bx^k, \bq^k, \by^k, \bt^k)\\ \in \mathcal{S}^k(\bp,\bp^s,\bz)}} \bigg\{p_0 y_0^k \nonumber\\
& + \sum_j p_j y_j^k +\sum_j \left(\phi_j^k+s^kp_j^s\right) t_j^k    + \sum_i \psi_i^k q_i^k  \nonumber\\
& +  w^k \Big( \sum_i x_{i,0}^k d_{i,0} + \sum_{i,j} x_{i,j}^k d_{i,j} \Big) \bigg\}, \forall k,\\
& z_j\in \{0,1\}, ~\forall j;~~ t_j^{\prime k} \in \{0,1\}, ~\forall j,k,\\
&x_{i,0}^{\prime k}, x_{i,j}^{\prime k}, q_{i}^{\prime k}, y_0^{\prime k}, y_{j}^{\prime k} \geq 0, \forall i,j,k,\\
&p_j = \sum_v p_j^v r_j^v;~~
\sum_v r_j^v = 1;~~ r_j^v\in\{0,1\}, ~~\forall j, v\\
&p^s_j = \sum_h p_j^{sh} r_j^{sh};~~
\sum_l r_j^{sh} = 1;~~ r_j^{sh}\in\{0,1\}, ~~\forall j, h.
\end{align*}

\subsection{The Karush–Kuhn–Tucker (KKT) Conditions}\label{app:KKTcond}
The Lagrangian function of the second maximization problem 
on the RHS of (\ref{sepenum}) is:
\beqn
\mathcal{L}(\bx,\!\bq,\by,\!\upsilon^k_{i,j},\!\upsilon^{0k}_{i}\!,\omega^{k}_{i},\!\gamma^k_{j},\!\gamma^{0k}\!,\mu_1^k,\mu_2^k,\nu^k_j,\!\Gamma^k_j,\!\sigma^k_j,\!\tau^k_{i,j},\!\eta^k_i,\xi^k_i) = \nonumber\\
-\!\sum_{i,j}\!\upsilon^k_{i,j} x^k_{i,j} \!-\!\! \sum_i\! \upsilon^{0k}_{i} x^k_{i,0} \!-\!\! \sum_i \!\omega^{k}_{i} q^k_i \!- \!\!\sum_{j} \!\gamma^k_{j} y^k_{j} \!-\! \gamma^{0k} y^k_{0} \!+\! p_0 y_0^k \nonumber\\
+ \! \sum_j p_j y_j^k \!+\! \sum_i \psi_i^k q_i^k + w^k \Big( \sum_i x_{i,0}^k d_{i,0}   +\!  \sum_{i,j} x_{i,j}^k d_{i,j} \Big)\nonumber\\
+ \mu_1^k \Big(\sum_j \left(\phi_j^k+s^kp_j^s\right)t_j^k+p_0 y_0^k + \sum_j p_j y_j^k - B^k\Big) \nonumber\\
+ \mu_2^k \Big(\!\sum_i x_{i,0}^k\!-\!y_0^k \!\Big) \!\!+\!\!\sum_j \Gamma^k_j \Big(\! \sum_{i}x^k_{i,j}\!-\! y_j^k\!\Big)\!\!+\!\! \sum_j \! \sigma^k_j \Big(\!y_j^k\!-\!C_j t_j^k\!\Big)\nonumber\\
+\! \sum_{i,j}\! \tau^k_{i,j} \Big(x^k_{i,j} \!-\! a^k_{i,j}R_i^k\Big) \!+\!\sum_i \eta^k_i \Big(x_{i,0}^k \!+ \!\!\sum_{j} x_{i,j}^k \!+\!  q_i^k \!- \!R_i^k\Big)\nonumber\\
+\nu_j^k(t_j^k - z_j)+\sum_i \xi^k_i \Big(x_{i,0}^k d_{i,0}+\sum_{j} x_{i,j}^k d_{i,j}-R_i^kD^{k,m}\Big).\nonumber
\eeqn
The KKT conditions give 
\begin{subequations}
\begin{alignat}{3}
\label{eq:kktxi0}
&\frac{\partial L}{\partial x^k_{i,0}} \!= \! -\upsilon^{0k}_{i}\!+\!w^kd_{i,0}\!+\!\mu_2^k\!+\!\eta^k_i \!+\!\xi^k_id_{i,0} \!=\! 0&&,\!\forall i,\!j,\!k \\
\label{eq:kktxij}
&\!\!\frac{\partial L}{\partial x^k_{i,j}} \!\!=\!  -\upsilon^k_{i,j}\!\!+\!w^k\!d_{i,j}\!\!+\!\Gamma^k_jv\!+\!\tau^k_{i,j} \!\!+\!\eta^k_i\!\!+\!\xi^k_id_{i,j} \!\!=\! 0&&,\!\forall i,\!j,\!k \\
\label{eq:kktqi}
&\frac{\partial L}{\partial q^k_i} = - \omega_i^k +  (\psi_i^k + \eta_i^k) = 0&&,\!\forall i,k \\
\label{eq:kkty0}
&\frac{\partial L}{\partial y^k_{0}} =  -\gamma^{0k}+ p_0+\mu_1^kp_0 -\mu_2^k = 0&&, \forall k\\
\label{eq:kktyj}
&\frac{\partial L}{\partial y^k_{j}} =  -\gamma^k_{j}+ p_j+\mu_1^kp_j -\Gamma^k_j+\sigma^k_j = 0&&,\forall j,k \\
\label{eq:kktCS1}
&\!\!0 \!\leq\! \mu_1^k \bot  B^k \!- \!\!\sum_j\! \left(\!\phi_j^k\!+\!s^kp_j^s\right)t_j^k\!-\!p_0 y_0^k \!- \!\!\sum_j \!p_j y_j^k &&\!\!\geq 0,\forall k \!\! \\
\label{eq:kktCS2}
&0 \leq \mu_2^k ~\bot~ y_0^k-\sum_i x_{i,0}^k \geq 0 &&,\forall k  \\
  \label{eq:kktCS3}
&0 \leq \Gamma^k_j ~\bot~ y_j^k-\sum_ix^k_{i,j} \geq 0&&, \forall j,k  \\
  \label{eq:kktCS4}
&0 \leq \sigma^k_j ~\bot~ C_j t_j^k-y_j^k \geq 0&&, \forall j,k  \\
  \label{eq:kktCS5}
&0 \leq \tau^k_{i,j} ~\bot~ a^k_{i,j}R_i^k- x^k_{i,j} \geq 0&&, \!\forall i,\!j,\!k  \\
  \label{eq:kktCS6}
&\sum_{j} x_{i,j}^k +  x_{i,0}^k + q_i^k = R_i^k&&, \forall i,k\\
  \label{eq:kktCS7}
&0 \leq \nu^k_j \!~\bot~\! z_j - t_j^k \geq 0&&, \forall j,k  \\
  \label{eq:kktCSe}
&0 \leq \xi^k_i \!~\bot~\! R_i^kD^{k,m}\!\!-x_{i,0}^k d_{i,0}-\!\!\sum_{j} x_{i,j}^k d_{i,j} \geq 0&&, \forall i,k  \\
\label{eq:kktCE1}
&0 \leq \upsilon^k_{i,j} ~\bot~ x^k_{i,j} \geq 0&&,\!\forall i,\!j,\!k\!\!\\
\label{eq:kktCE2}
&0 \leq \upsilon^{0k}_{i} ~\bot~ x^k_{i,0} \geq 0&&,\forall i,k\\
\label{eq:kktCE2b}
&0 \leq \omega^{k}_{i} ~\bot~ q^k_{i} \geq 0&&,\forall i,k\\
  \label{eq:kktCE3}
&0 \leq \gamma^k_{j} ~\bot~ y^k_{j} \geq 0&&,\forall j,k\\
\label{eq:kktCE4}
&0 \leq \gamma^{0k} ~\bot~ y^k_{0} \geq 0&&,\forall k
\end{alignat}
\end{subequations}
where (\ref{eq:kktxi0})-(\ref{eq:kktyj}) are the stationary conditions, (\ref{eq:kktCS1})-(\ref{eq:kktCE4}) are
the primal feasibility, dual feasibility and complementary slackness conditions.
From (\ref{eq:kktxij}), we have:
\beqn
 w^kd_{i,j}+\Gamma^k_j+\tau^k_{i,j}+\eta^k_i +\xi^k_id_{i,j}= \upsilon^k_{i,j},~\forall i,j,k. \label{eq:kktex1}
\eeqn
Moreover, from \eqref{eq:kktCE1}, if $x^k_{i,j} > 0$, then $\upsilon^k_{i,j}=0$ and \eqref{eq:kktex1} implies that 
$$ w^kd_{i,j}+\Gamma^k_j+\tau^k_{i,j} +\eta^k_i+\xi^k_id_{i,j} =0,~\forall i,j,k.$$
Thus, 
\beqn \label{eq:kktex2}
\Big( w^kd_{i,j}+\Gamma^k_j+\tau^k_{i,j} +\xi^k_id_{i,j}\Big)x_{i,j} =0,~\forall i,j,k.
\eeqn
Therefore, the conditions (\ref{eq:kktxij}), \eqref{eq:kktCE1} and \eqref{eq:kktex2} are equivalent to the following constraint:
\beqn
\label{eq:kktxijex}
0 \leq  w^kd_{i,j}\!+\Gamma^k_j\!+\tau^k_{i,j}\!+\!\eta^k_i \!+\!\xi^k_id_{i,j} ~\bot~ x^k_{i,j} \geq 0,\forall i,j,k.
\eeqn
Similarly, we have
\begin{align}
\label{eq:kktxi0ex}
0 \leq  w^kd_{i,0}+\mu_2^k +\eta^k_i+\xi^k_id_{i,0} ~\bot~ x^k_{i,0} \geq 0&,\forall i,k,\\
\label{eq:kktqiex}
0 \leq  (\psi_i^k + \eta_i^k) ~\bot~ q^k_i \geq 0&,\forall i,k,\\
\label{eq:kktyjex}
0 \leq  p_j+\mu_1^kp_j -\Gamma^k_j+\sigma^k_j ~\bot~ y^k_{j} \geq 0&,\forall j,k,\\
\label{eq:kkty0ex}
0 \leq  p_0+\mu_1^kp_0 -\mu_2^k ~\bot~ y_0^k \geq 0&,\forall k.
\end{align}
It is worth noting that we can directly use the  set of KKT conditions  (\ref{eq:kktxi0})-(\ref{eq:kktyj})
 to solve the problem, but it will involve more variables (i.e., $\upsilon$, $\omega$ and $\gamma$) compared to solving the subproblem with (\ref{eq:kktxijex})-(\ref{eq:kkty0ex}).
 
 In brief, based on  the KKT conditions above, we can infer that
the problem (\ref{sepenum}) is equivalent to the following problem with complementary constraints, for all $k$:
\begin{align}
\label{eq:sepenumKKT}
&\!\!\!w^k \Big( \sum_i x_{i,0}^{\prime k} d_{i,0} + \sum_{i,j} x_{i,j}^{\prime k} d_{i,j} \Big) + \! \sum_i \psi_i^k q_i^{\prime k}   +  \sum_j p_j y_j^{\prime k}\nonumber\\
& + p_0 y_0^{\prime k} + \sum_j (\phi_j^k+s^kp_j^s) t_j^{\prime k} \leq \min_{\substack{\bt^k \in  \bT^k\!,\bx^k,\bq^k,\by^k\!,\\\bu^k\!,\mu_1,\mu_2,\nu,\Gamma,\sigma,\tau,\xi}}    p_0 y_0^k \nonumber\\
&+ \sum_j p_j y_j^k + \sum_{j}(\phi_j^k+s^kp_j^s)t_j^k +  \sum_i \psi_i^k q_i^k\nonumber\\
&+ w^k \Big( \sum_i x_{i,0}^k d_{i,0}   +  \sum_{i,j} x_{i,j}^k d_{i,j} \Big)
\end{align}
subject to
\begin{align*}
% \label{eq:kktstart}
% &0 \leq  w^kd_{i,j}+\Gamma^k_j+\tau^k_{i,j}+\eta^k_i +\xi^k_id_{i,j} ~\bot~ x^k_{i,j} \geq 0,&&\forall i,j\\
% &0 \leq  w^kd_{i,0}+\mu_2^k +\eta^k_i+\xi^k_id_{i,0} ~\bot~ x^k_{i,0} \geq 0,&&\forall i\\
% &0 \leq  p_j+\mu_1^kp_j -\Gamma^k_j+\sigma^k_j ~\bot~ y^k_{j} \geq 0,&&\forall j\\
% &0 \leq  p_0+\mu_1^kp_0 -\mu_2^k ~\bot~ y_0^k \geq 0,&&\\
% & 0 \leq \mu_1^k \bot  B^k - \sum_j (\phi_j^k+s^kp_j^s)t_j^k-p_0 y_0^k - \sum_j p_j y_j^k &&\geq 0  \\
% &0 \leq \mu_2^k ~\bot~ y_0^k-\sum_i x_{i,0}^k \geq 0&&  \\
% &0 \leq \Gamma^k_j ~\bot~ y_j^k-\sum_ix^k_{i,j} \geq 0&&, \forall j  \\
% &0 \leq \sigma^k_j ~\bot~ C_j t_j^k-y_j^k \geq 0&&, \forall j  \\
% &0 \leq \tau^k_{i,j} ~\bot~ a^k_{i,j}R_i^k- x^k_{i,j} \geq 0&&, \forall i,j  \\
% &\sum_{j} x_{i,j}^k +  x_{i,0}^k = R_i^k&&, \forall i\\
% \label{eq:kktend}
% &0 \leq \xi^k_i ~\bot~ R_i^kD^{k,m}-x_{i,0}^k d_{i,0}-\sum_{j} x_{i,j}^k d_{i,j} \geq 0&&, \forall i.
\eqref{eq:kktxijex} - \eqref{eq:kkty0ex}, ~ \eqref{eq:kktCS1} - \eqref{eq:kktCSe}.
\end{align*}

Note that a complimentary constraint $0 \leq x \bot \pi \geq 0$ means $x \geq 0, \pi \geq$ and $x . \pi = 0$. Thus, it is a nonlinear constraint. However, this nonlinear complimentary constraint can be transformed into equivalent exact linear constraints by using the Fortuny-Amat transformation \cite{bigM}. Specifically, the complementarity  condition $0 \leq x \bot \pi \geq 0$ is  equivalent  to the following set of mixed-integer linear constraints:
\beqn
x \geq 0;~~ x \leq (1-u)M 
\pi \geq 0;~~ \pi \leq  uM,~~ u\in\{0;1\},
\eeqn
where $M$ is a sufficiently large constant.
By applying this transformation to all the complementary constraints listed in the KKT reformulation in \eqref{eq:sepenumKKT} above, we obtain an MILP that is equivalent to the subproblem (\ref{sepenum}).

\subsection{KKT-based Reformulation}
\label{KKTccg}
We present a single-level equivalent reformulation of \textbf{BMIP} in this subsection which serves as the foundation for our solution scheme. The primary goal of this reformulation is to expand (\ref{enumx}) by enumeration. We assume that the remaining lower-level problem has a finite optimal value for any possible ($p, z, t$). 
We now apply the classical reformulation method using KKT conditions to the second maximization problem on the RHS of (\ref{sepenum}). 
% Please refer to Appendix~\ref{app:KKTcond} for more details. 

The explicit form of the MILP that is equivalent to the subproblem (\ref{sepenum}) is given for all $k$ as follows:
\begin{align}
\label{sepenummilpMILP}
&w^k \Big(\! \sum_i\! x_{i,0}^{\prime k} d_{i,0} \!+\!\! \sum_{i,j}\! x_{i,j}^{\prime k} d_{i,j} \!\Big) \!+\! \!\sum_i\! \psi_i^k q_i^{\prime k} \!+\! p_0 y_0^{\prime k} \!+\!\!  \sum_j p_j y_j^{\prime k} \nonumber\\
&+ \!\sum_j \left(\phi_j^k+s^kp_j^s\right) t_j^{\prime k}\! \leq \!\min_{\substack{\bt^k \in  \bT^k\!,\bx^k,\bq^k,\by^k\!,\\\bu^k\!,\mu_1,\mu_2,\nu,\Gamma,\sigma,\tau,\xi}}  \sum_{j} \left(\phi_j^k+s^kp_j^s\right) t_j^k  \nonumber\\
&+ p_0 y_0^k+\!\sum_j p_j y_j^k+ w^k \Big( \sum_i x_{i,0}^k d_{i,0} +\! \sum_{i,j} x_{i,j}^k d_{i,j} \Big) \!\!\!
\end{align}
subject to
\begin{align*}
&0 \leq B^k - \!\sum_j \left(\phi_j^k+s^kp_j^s\right)t_j^k-p_0 y_0^k - \!\sum_j p_j y_j^k  \leq  u^{1,k}M^{1,k},\\
&0 \leq \mu_1^k \leq (1-u^{1,k}) M^{1,k}, \\
&0 \leq y_0^k-\sum_i x_{i,0}^k \leq u^{2,k} M^{2,k},  \\
&0 \leq \mu_2^k \leq (1-u^{2,k}) M^{2,k},  \\
&0 \leq y_j^k-\sum_ix^k_{i,j} \leq u_{j}^{3,k} M_{j}^{3,k} ,~~\forall j \\
&0 \leq \Gamma^k_j \leq (1-u_{j}^{3,k}) M_{j}^{3,k} ,~~\forall j \\
&0 \leq C_j t_j^k-y_j^k \leq u_{j}^{4,k} M_{j}^{4,k} ,~~\forall j \\
&0 \leq \sigma^k_j \leq (1-u_{j}^{4,k}) M_{j}^{4,k} ,~~\forall j \\
&0 \leq  a^k_{i,j}R_i^k- x^k_{i,j} \leq u_{i,j}^{5,k} M_{i,j}^{5,k} ,~~\forall i,j\\
&0 \leq \tau^k_{i,j} \leq (1-u_{i,j}^{5,k}) M_{i,j}^{5,k} ,~~\forall i,j\\
&\sum_{j} x_{i,j}^k +  x_{i,0}^k + q_i^k= R_i^k,~~ \forall i,k\\
&0 \leq z_j-t_j^k \leq u_{j}^{6,k} M_{j}^{6,k} ,~~\forall j \\
&0 \leq \nu^k_j \leq (1-u_{j}^{6,k}) M_{j}^{6,k} ,~~\forall j \\
&0 \leq R_i^kD^{k,m}-x_{i,0}^k d_{i,0}-\sum_{j} x_{i,j}^k d_{i,j} \leq u_{i}^{7,k} M_{i}^{7,k} ,~~\forall i,k\\
&0 \leq \xi^k_i \leq (1-u_{i}^{7,k}) M_{i}^{7,k} ,~~\forall i,k\\
&0 \leq  w^kd_{i,j}+\Gamma^k_j+\tau^k_{i,j}+\eta^k_i +\xi^k_id_{i,j} \leq u_{i,j}^{8,k} M_{i,j}^{8,k} ,~~\forall i,j\\
&0 \leq x^k_{i,j} \leq (1-u_{i,j}^{8,k}) M_{i,j}^{8,k} ,~~\forall i,j\\
&0 \leq w^kd_{i,0}+\mu_2^k +\eta^k_i+\xi^k_id_{i,0} \leq u_{i}^{9,k} M_{i}^{9,k} ,~~\forall i,k\\
&0 \leq x^k_{i,0} \leq (1-u_{i}^{9,k}) M_{i}^{9,k} ,~~\forall i,k\\
&0 \leq (\psi_i^k+\eta_i^k) \leq u_{i}^{10,k} M_{i}^{10,k} ,~~\forall i,k\\
&0 \leq q^k_i \leq (1-u_{i}^{10,k}) M_{i}^{10,k} ,~~\forall i,k\\
&0 \leq p_j+\mu_1^kp_j -\Gamma^k_j+\sigma^k_j \leq u_{j}^{11,k} M_{j}^{11,k} ,~~\forall j \\
&0 \leq y^k_{j} \leq (1-u_{j}^{11,k}) M_{j}^{11,k} ,~~\forall j \\
&0 \leq p_0+\mu_1^kp_0 -\mu_2^k \leq u^{12,k} M^{12,k}, \\
&0 \leq y_0^k \leq (1-u^{12,k}) M^{12,k}, \\
&\bu\in \{0,1\},
% &u^{1,k}\!,\!u^{2,k}\!,\!u_{j}^{3,k}\!,\!u_{j}^{4,k}\!,\!u_{i,j}^{5,k}\!,\!u_{i}^{7,k}\!,\!u_{i,j}^{8,k}\!,\!u_{i}^{9,k}\!,\!u_{j}^{11,k}\!,\!&&\!\!u^{12,k}\!\in \{0,1\}.
\end{align*}
where $\bu$ represents the set of binary variables $u^{1,k}$, $u^{2,k}$, $u_{j}^{3,k}$, $u_{j}^{4,k}$, $u_{i,j}^{5,k}$, $u_{j}^{6,k}$, $u_{i}^{7,k}$, $u_{i,j}^{8,k}$, $u_{i}^{9,k}$, $u_{i}^{10,k}$, $u_{j}^{11,k}$, $u^{12,k}$. Also,  $M^{1,k}$, $M^{2,k}$, $M_{j}^{3,k}$, $M_{j}^{4,k}$, $M_{i,j}^{5,k}$, $M_{j}^{6,k}$, $M_{i}^{7,k}$, $M_{i,j}^{8,k}$, $M_{i}^{9,k}$, $M_{i}^{10,k}$, $M_{j}^{11,k}$, $M^{12,k}$ are sufficiently large numbers. The value of each $M$ should be large enough to ensure feasibility of the associated constraint. On the other hand, the value of each $M$ should not be too large to enhance the computational speed of the solver. Indeed, the value of each $M$ should be tighten to the limits of parameters and variables in the corresponding constraint. For instance, $M_i^0$ needs to be larger or equal to the maximum value of $x_{i}$, which is $\lambda_i$. Thus,  $M_i^0$ can be set to be a small number greater than $\lambda_i$.

\subsection{Bilinear Linearization}\label{app-BilinearLinearization}
To handle the bilinear terms $p_j \mu_1^k$, we can utilize (\ref{ipc5}) and express it in a different form as follows:
\beqn
\label{pmu1linearccg}
p_j \mu_1^k = \sum_v p_j^v r_j^v \mu_1^k = \sum_v p_j^v \pi_j^{v,k},
\eeqn
where $\pi_j^{v,k} = r_j^v \mu_1^k$. Note that $\pi_j^{v,k}$ is a continuous variable and we have $\pi_j^{v,k} = \mu_1^k$ if $r_j^v = 1$ and $\pi_j^{v,k} = 0$, otherwise. Hence, using (\ref{pmu1linearccg}), the bilinear term $p_j \mu_1^k$ can be written as a linear function of $\pi_j^{k} = (\pi_j^{1,k}, \ldots, \pi_j^{V,k})$. 

Additionally, the constraints $\pi_j^{v,k} = r_j^v \mu_1^k, \forall j,k$ can be implemented through the following linear inequalities \cite{linear2}:
\begin{subequations}\label{pmulinearccg}
\begin{align}
\label{pmulinear1ccg}
\pi_j^{v,k} \leq M r_j^v,~ \forall j, k, v;~
\pi_j^{v,k} \leq \mu_1^k,~\forall j, k, v,\\
\label{pmulinear2ccg}
\pi_j^{v,k} \geq 0, ~ \forall j, k, v;~ \pi_j^{v,k} \geq \mu_1^k + M r_j^v - M, ~ \forall j, k, v,
\end{align}
\end{subequations}
where $M$ is a sufficiently large number. 

Similarly, the bilinear terms $\sum_{j,k} p_j y_j^{\prime k}$ in  (\ref{sepenummilpMILP}) is products of discrete and continuous variables. By expressing $p_j$ using (\ref{ipc5}), i.e.,
\beqn
\label{pyplinearccg}
p_j y_j^{\prime k} = \sum_v p_j^{v} r_j^{v} y_j^{\prime k} = \sum_v p_j^{v} \rho_j^{v,k},
\eeqn
we can linearize bilinear terms $\rho_j^{v,k}=r_j^v y_j^{\prime k}$ using the following linear inequalities: 
\begin{subequations}\label{pyprimelinearccg}
\begin{align}
\label{pyprimelinear1ccg}
\rho_j^{v,k} \leq M r_j^v,~ \forall j, k, v;~
\rho_j^{v,k} \leq y_j^{\prime k},~\forall j, k, v,\\
\label{pyprimelinear2ccg}
\rho_j^{v,k} \geq 0, ~ \forall j, k, v;~ \rho_j^{v,k} \geq y_j^{\prime k} + M r_j^v - M, ~ \forall j, k, v.
\end{align}
\end{subequations}

The bilinear terms $p_j^s t_j^{\prime k}$ in (\ref{sepenummilpMILP}), are discretized using (\ref{ipc6}) for $p_j^s$, as follows,
\beqn
\label{pstlinearccg}
p^s_j t_j^{\prime k} = \sum_h p_j^{sh} r_j^{sh} t_j^{\prime k} = \sum_h p_j^{sh} \zeta_j^{h,k},
\eeqn
where $\zeta_j^{h,k}=r_j^{sh} t_j^{\prime k}$ can be expressed as linear inequalities: \begin{subequations}\label{ptlinearccg}
\begin{align}
\label{ptlinear1ccg}
\zeta_j^{h,k} \leq r_j^{sh},~ \forall j, k, h;~
\zeta_j^{h,k} \leq t_j^{\prime k},~\forall j, k, h,\\
\label{ptlinear2ccg}
\zeta_j^{h,k} \geq 0, ~ \forall j, k, h;~ \zeta_j^{h,k} \geq r_j^{sh}+t_j^{\prime k}-1, ~ \forall j, k, h.
\end{align}
\end{subequations}

Following a similar procedure, the bilinear terms $\varrho_j^{k}=\nu_j^{k}z_j$ in \eqref{eq:SubProbStrongDuality3} can be implemented using linear inequalities:
\begin{subequations}\label{nuzlinearccg}
\begin{align}
\label{nuzlinear1ccg}
\varrho_j^{k} \leq M z_j,~ \forall j, k;~
\varrho_j^{k} \leq \nu_j^{k},~\forall j, k,\\
\label{nuzlinear2ccg}
\varrho_j^{k} \geq 0, ~ \forall j, k;~ \varrho_j^{k} \geq \nu_j^{k} + M z_j - M, ~ \forall j, k.
\end{align}
\end{subequations}

Based on the aforementioned linearization steps, we can represent the bi-level problem given by (\ref{eq:bmipd}) as an equivalent single-level MILP. It is worth noticing that this approach results in a large number of variables and constraints in the reformulation compared to the reformulation using LP duality as presented in the main manuscript.

% \subsection{Single-Level MILP Reformulation of \bmip~with KKT Reformulation }
% \label{app-singlemilp-KKT}

\subsection{Single-Level MILP Reformulation of \bmip }
\label{app-singlemilp-dual}
The single-level mixed-integer linear program \bmipl~representing the leader problem is presented as follows:
\begin{align}\label{SMIP-objmd}
\bmipl: \max_{\substack{\bp,\bp^s \!,\bz,\bx^\prime \!,\by^\prime \!,\\\bt^\prime \!, \br, \mu, \Gamma, \nu, \sigma,\\ \xi, \tau, \kappa, \pi, \rho, \varrho, \zeta}} \sum_{j,k}&\! \left(\!\sum_v p_j^{v} \rho_j^{v,k} \!+  \phi_j^kt_j^{\prime k}\!+s^k \!\sum_h p_j^{sh} \zeta_j^{h,k}\right) \nonumber\\
&- \sum_j \left(f_j z_j + c_j  \frac{\sum_k y_j^{\prime k}}{C_j} \!\!\right)  
\end{align}
subject to
\begin{align*}
&\sum_k y_j^{\prime k} \leq z_j C_j, \forall j;~ \sum_k s^k t_j^{\prime k} \leq z_jS_j, \forall j;~ y_0^{\prime k} \geq \sum_i x_{i,0}^{\prime k}, \forall j,k\\
&p_0 y_0^{\prime k} + \!\sum_{j}\!\Big(\sum_{v} p_j^{v} \rho_j^{v,k} + \phi_j^kt_j^{\prime k}+s^k\!\sum_h p_j^{sh} \zeta_j^{h,k}\Big)  \!\leq B^k\!,~ \forall k \\
&t_j^{\prime k} \leq z_j,~~ y_j^{\prime k} \geq \sum_i x_{i,j}^{\prime k}, ~~y_j^{\prime k} \leq C_j t_j^{\prime k}, ~~s^k t_j^{\prime k} \leq S_j,~\forall j,k \\
&x_{i,j}^{\prime k} \leq a_{i,j}^k R_i^k, ~\forall i,j,k~; ~~\sum_j x_{i,j}^{\prime k} + x_{i,0}^{\prime k} + q_{i}^{\prime k} = R_i^k, ~\forall i,k \\
&x_{i,0}^{\prime k} d_{i,0} +  \sum_j x_{i,j}^{\prime k} d_{i,j} \leq D^{k, \sf m} R_i^k, ~~\forall i,k \\
%\label{cut1}
&w^k \Big( \sum_i x_{i,0}^{\prime k} d_{i,0} + \sum_{i,j} x_{i,j}^{\prime k} d_{i,j} \Big) + \! \!\sum_i \psi_i^k p_i^{\prime k}  + \sum_{j,v} p_j^{v} \rho_j^{v,k} \nonumber\\
&+ p_0 y_0^{\prime k} +\!\sum_j\! \left(\!\phi_j^kt_j^{\prime k} \!+ \! s^k\sum_h p_j^{sh} \zeta_j^{h,k}\!\right) \!\leq\! \sum_j \!\left(\phi_j^k+s^kp_j^s\right) t_j^{k,l}  \nonumber\\ 
&-B^k\mu_1^{k,l} +\sum_j \left(\phi_j^k\mu_1^{k,l}+s^k\sum_h p_j^{sh}  \kappa_j^{h,k,l}\right) t_j^{k,l} \nonumber\\
&+\nu_j^{k,l}t_j^{k,l} - \varrho_j^{k,l} -\sum_i R_i^k D^{k,m} \xi_i^{k,l} + \sum_i R_i^k \eta_i^{k,l}\nonumber\\
&- \sum_j C_j t_j^{k,l} \sigma_j^{k,l} - \sum_{i,j} a_{i,j}^k \tau_{i,j}^{k,l} R_i^k, ~~\forall k, 1 \leq l \leq L\\
%\label{cut2}
&p_0 \left( 1 + \mu_1^{k,l} \right) - \mu_2^{k,l} \geq 0,~~\forall k, 1 \leq l \leq L\\
%\label{cut3}
&\sum_v p_j^v \pi_j^{v,k,l} -\Gamma_j^{k,l} + \sigma_j^{k,l} + p_j \geq 0,~~\forall j,k, 1 \leq l \leq L\\
%\label{cut4}
& \eta_i^{k,l} \leq \psi_i^k,~~\forall i,k, 1 \leq l \leq L\\
&\mu_2^{k,l} + d_{i,0} \xi_i^{k,l} - \eta_i^{k,l} \geq -w^k d_{i,0},~~\forall i,k, 1 \leq l \leq L\\
%\label{cut5}
&\Gamma_j^{k,l} + d_{i,j} \xi_i^{k,l} + \tau_{i,j}^{k,l} - \eta_i^{k,l} \geq -w^k d_{i,j},\forall i,j,k, 1 \leq l \leq L\\
&p_j = \sum_v p_j^v r_j^v,~\forall j;~
\sum_v r_j^v = 1,~\forall j;~ r_j^v\in\{0,1\}, ~\forall j, v\\
&p^s_j = \sum_h p_j^{sh} r_j^{sh},~\forall j;~
\sum_l r_j^{sh} = 1,~\forall j;~ r_j^{sh}\in\{0,1\}, ~\forall j, h\\
&\kappa_j^{h,k,l} \leq M r_j^{sh},~
\kappa_j^{h,k,l} \leq \mu_1^{k,l},~\kappa_j^{h,k,l} \geq \mu_1^{k,l} + M r_j^{sh} - M, ~ \forall j, k, h,l\\
&\pi_j^{v,k,l} \leq M r_j^{v},~ \pi_j^{v,k,l} \leq \mu_1^{k,l},~\pi_j^{v,k,l} \geq \mu_1^{k,l} + M r_j^{v} - M, ~ \forall j, k, v,l\\
&\rho_j^{v,k} \leq M r_j^v,~
\rho_j^{v,k} \leq y_j^{\prime k},~\rho_j^{v,k} \geq y_j^{\prime k} + M r_j^v - M, ~ \forall j, k, v\\
&\zeta_j^{h,k} \leq r_j^{sh},~ \zeta_j^{h,k} \leq t_j^{\prime k},~\zeta_j^{h,k} \geq r_j^{sh}+t_j^{\prime k}-1, ~ \forall j, k, h\\
&\varrho_j^{k,l} \leq M z_j,~ \varrho_j^{k,l} \leq \nu_j^{k,l},~\varrho_j^{k,l} \geq \nu_j^{k,l} + M z_j - M, ~ \forall j, k,l,\\
&x_{i,j}^{\prime k},x_{i,0}^{\prime k},y_j^{\prime k},y_0^{\prime k}\geq 0, ~~\forall i,j,k \\
&\mu_1^{k,l}, \Gamma_j^{k,l}, \nu_j^{k,l}, \mu_2^{k,l}, \sigma_j^{k,l}, \xi_i^{k,l}, \tau_{i,j}^{k,l} \geq 0, ~~\forall i,j,k, 1 \leq l \leq L\\
&\kappa_j^{h,k,l}, \pi_j^{v,k,l}, \rho_j^{v,k}, \zeta_j^{h,k}, \varrho_j^{k,l} \geq 0, ~~\forall i,j,k,v,h, 1 \leq l \leq L\\
&z_j\in \{0,1\};~~ t_j^{\prime k} \in \{0,1\}; ~~ \eta_i^{k,l} \in\mathbb{R}, ~\forall i,j,k, 1 \leq l \leq L.
\end{align*}